\documentclass[aap,MSNbibl,noautosecdot,citesort,dvips]{arximspdf}
\usepackage{subenv}
\usepackage{graphicx}

\doi{10.1214/11-AAP759}
\volume{22}
\issue{1}
\pubyear{2012}
\firstpage{70}
\lastpage{127}

\makeatletter

\newcommand{\tendsdown}{\alpha\searrow0}

\newcommand{\implies}{\Longrightarrow}

\newcommand{\third}{\frac{1}{3}}
\newcommand{\twothirds}{\frac{2}{3}}

\newcommand{\bZero}{\mathbf{0}}
\newcommand{\bOne}{\mathbf{1}}
\newcommand{\LIFT}{\Delta}
\newcommand{\LIFTW}{\Delta{W}}
\newcommand{\CLVR}{\Xi}
\newcommand{\bvarepsilon}{\bolds{\varepsilon}}
\newcommand{\bdelta}{\bolds{\delta}}
\newcommand{\tbdelta}{\vec{\bdelta}}
\newcommand{\bpi}{{\bolds{\pi}}}
\newcommand{\brho}{{\bolds{\rho}}}
\newcommand{\bsigma}{\bolds{\sigma}}
\newcommand{\blambda}{\bolds{\lambda}}
\newcommand{\dLambda}{\partial\Lambda}
\newcommand{\tlambda}{\vec{\lambda}}
\newcommand{\tblambda}{\vec{\blambda}}
\newcommand{\bxi}{{\bolds{\xi}}}
\newcommand{\bzeta}{{\bolds{\zeta}}}
\newcommand{\bA}{\mathbf{A}}
\newcommand{\fba}{\tilde{\mathbf{a}}}
\newcommand{\fa}{\tilde{a}}
\newcommand{\cA}{\mathcal{A}}
\newcommand{\hba}{\hat{\mathbf{a}}}
\newcommand{\bB}{\mathbf{B}}
\newcommand{\bc}{\mathbf{c}}
\newcommand{\cD}{\mathcal{D}}
\newcommand{\hi}{{\hat{\imath}}}
\newcommand{\hj}{{\hat{\jmath}}}
\newcommand{\cJ}{\mathcal{J}}
\newcommand{\cK}{\mathcal{K}}
\newcommand{\N}{{N}}
\newcommand{\bn}{\mathbf{n}}
\newcommand{\bP}{\mathbf{P}}
\newcommand{\tbP}{\vec{\bP}}
\newcommand{\bQ}{\mathbf{Q}}
\newcommand{\fbq}{\tilde{\mathbf{q}}}
\newcommand{\fq}{\tilde{q}}
\newcommand{\hbq}{\hat{\mathbf{q}}}
\newcommand{\tbq}{\vec{\mathbf{q}}}
\newcommand{\tq}{\vec{q}}
\newcommand{\tqsup}{\vec{q}}
\newcommand{\tbr}{\vec{\mathbf{r}}}
\newcommand{\br}{\mathbf{r}}
\newcommand{\hbr}{\hat{\mathbf{r}}}
\newcommand{\tR}{\vec{R}}
\newcommand{\bs}{\mathbf{s}}
\newcommand{\hs}{\hat{s}}
\newcommand{\fls}{\tilde{s}}
\newcommand{\sS}{\mathcal{S}}
\newcommand{\T}{\mathsf{T}}
\newcommand{\Tfluid}{T^{\mathrm{fluid}}}
\newcommand{\bu}{\mathbf{u}}
\newcommand{\bv}{\mathbf{v}}
\newcommand{\WORK}{W}
\newcommand{\hWORK}{\hat{W}}
\newcommand{\cX}{\mathcal{X}}
\newcommand{\bx}{\mathbf{x}}
\newcommand{\tbx}{\vec{\bx}}
\newcommand{\fx}{\tilde{x}}
\newcommand{\bY}{\mathbf{Y}}
\newcommand{\fby}{\tilde{\mathbf{y}}}
\newcommand{\fy}{\tilde{y}}
\newcommand{\hbz}{\hat{\mathbf{z}}}
\newcommand{\bZ}{\mathbf{Z}}
\newcommand{\PRIMAL}{\operatorname{PRIMAL}}
\newcommand{\DUAL}{\operatorname{DUAL}}
\newcommand{\ALGDUAL}{\operatorname{ALGD}}
\newcommand{\ALGD}{\operatorname{ALGD}}
\newcommand{\SSC}{\mathcal{W}}
\newcommand{\SSCm}{\mathcal{W}^{\max}}
\newcommand{\INV}{\mathcal{I}}
\newcommand{\FMS}{\operatorname{FMS}}
\newcommand{\FMSm}{\operatorname{FMSm}}
\newcommand{\hcW}{\hat{\mathcal{W}}}
\newcommand{\modcont}{\operatorname{mc}}
\newcommand{\clust}{\mathrm{CP}}
\newcommand{\hE}{\hat{E}}

\newcommand{\argmin}{\arg\min}
\newcommand{\tend}{\to}
\newcommand{\Prob}{\mathbb{P}}
\newcommand{\Reals}{\mathbb{R}}
\newcommand{\IntegersP}{\mathbb{Z}_+}
\newcommand{\Naturals}{\mathbb{N}}
\newcommand{\RealsP}{\Reals_+}
\newcommand{\mmin}{\wedge}
\newcommand{\mmax}{\vee}
\newcommand{\argmax}{\arg\max}
\newcommand{\union}{\cup}
\newcommand{\bdot}{\cdot}

\newproclaim{assumption}{Assumption}[section]
\newproclaim{definition}[assumption]{Definition}

\newtheorem{proposition}[assumption]{Proposition}
\newtheorem{theorem}[assumption]{Theorem}
\newtheorem{corollary}[assumption]{Corollary}
\newtheorem{lemma}[assumption]{Lemma}
\newtheorem{conjecture}[assumption]{Conjecture}

\makeatother

\begin{document}
\begin{frontmatter}

\title{Switched networks with maximum weight policies: Fluid approximation
and multiplicative state space collapse}
\runtitle{Network scheduling}

\begin{aug}
\author[A]{\fnms{Devavrat} \snm{Shah}\corref{}\thanksref{t1}\ead[label=e1]{devavrat@mit.edu}}
and
\author[B]{\fnms{Damon} \snm{Wischik}\thanksref{t2}\ead[label=e2]{D.Wischik@cs.ucl.ac.uk}}
\runauthor{D. Shah and D. Wischik}
\affiliation{Massachusetts Institute of Technology and University
College London}
\address[A]{Department of EECS\\
Massachusetts Institute of Technology\\
Cambridge, Massachusetts 02139\\
USA\\
\printead{e1}} 
\address[B]{Department of Computer Science\\
University College London\\
Gower St.\\
London WC1E 6BT\\
United Kingdom\\
\printead{e2}}
\end{aug}

\thankstext{t1}{Supported by NSF CAREER CNS-0546590.}
\thankstext{t2}{Supported by a Royal Society university research
fellowship. Collaboration supported
by the British Council Researcher Exchange program.}

\received{\smonth{8} \syear{2007}}
\revised{\smonth{1} \syear{2011}}

%
\begin{abstract}
We consider a queueing network in which there are constraints on which
queues may be served simultaneously; such networks may be used to
model input-queued switches and wireless networks. The scheduling
policy for such a network specifies which queues to serve at any point
in time. We consider a family of scheduling policies, related to the
maximum-weight policy of Tassiulas and Ephremides [\textit{IEEE Trans. Automat.
Control} \textbf{37} (1992) 1936--1948], for single-hop and
multihop networks. We specify a fluid model and show that
fluid-scaled performance processes can be approximated by fluid model
solutions. We study the behavior of fluid model solutions under critical
load, and characterize invariant states as those states which solve a
certain network-wide optimization problem. We use fluid model results
to prove multiplicative state space collapse. A notable feature of our
results is that they do not assume complete resource pooling.
\end{abstract}

%

\end{frontmatter}

\section{\texorpdfstring{Introduction.}{Introduction}}
\label{sec:intro}

A switched network consists of a collection of queues, operating
in discrete time. In each time slot, queues are offered
service according to a \textit{service schedule} chosen
from a specified finite set. For example, in a three-queue
system, the set of allowed schedules might consist of ``Serve~3 units
of work each from queues $A$ and $B$'' and ``Serve 1 unit of work each
from queues $A$ and $C$, and 2 units from queue $B$.''
The rule for choosing a~schedule is called the \textit{scheduling policy}.
New work may
arrive in each time slot; let each queue have a dedicated exogenous
arrival process, with specified mean arrival rates. Once work is
served, it may either rejoin one of the queues or leave the network.

Switched networks are special cases of what Harrison
\cite{harrisoncanonical,harrisoncanonicalcorr} calls
``stochastic
processing networks.'' We believe that switched networks are general
enough to model a variety of interesting applications. For example,
they have been used to model input-queued switches, the devices at the
heart of high-end internet routers, whose underlying silicon
architecture imposes constraints on which traffic streams can be
transmitted simultaneously \cite{daibala}. They have also been used to
model a multi-hop wireless network in which interference limits the
amount of service that can be given to each host \cite{tassiula1}.

The main result of this paper is Theorem \ref{thm:heavytraffic}, which
proves multiplicative state space collapse
(as defined in Bramson \cite{bramson})
for a switched network
running a generalized version of the maximum-weight scheduling policy
of Tassiulas and Ephremides \cite{tassiula1}, in the diffusion (or
heavy traffic)
limit. Whereas previous works on switched networks and
stochastic processing networks in the diffusion
limit \cite{stolyar,LD05,LD08} have assumed the ``complete resource
pooling'' condition, which roughly means that there is a single
bottleneck cut constraint, we do not make this assumption.
Section \ref{sec:discussion} discusses further the related work
and our reasons for being interested in the case without complete resource
pooling.

To prove multiplicative state space collapse we follow the general
method laid out by Bramson \cite{bramson}. In Section \ref{sec:model} we
specify a stochastic switched network model and describe the
generalized maximum-weight policy. In Section~\ref{sec:fluid} we
specify a fluid model and prove that
fluid-scaled performance processes of the switched network
are approximated by solutions of this fluid
model. Sections \ref{sec:fluid.behavior} and
\ref{sec:fluid.behavior.m} give properties of the solutions of the
fluid model for single-hop and multi-hop networks,
respectively. Specifically, for nonoverloaded fluid model solutions,
we characterize the invariant states and prove
that fluid model solutions
converge towards the
set of invariant states.
In Section~\ref{sec:ht} we use these properties to prove
multiplicative state space collapse.

We use the cluster-point method of Bramson~\cite{bramson} to prove the
fluid model approximation in Section \ref{sec:fluid}, rather than
following an approach based on weak convergence. The former
yields uniform bounds on the error of fluid model approximations, and
these uniform bounds are needed in proving multiplicative state space
collapse. However, the assumptions we make on the arrival process are
not the same as those of Bramson \cite{bramson}.

In Section \ref{sec:optimal} we give results concerning the fluid
model behavior of a~general single-hop switched network in critical
load, and the set of invariant states for the input-queued
switch, under a condition that we call ``complete loading.'' Motivated
by these results, we define a scheduling policy which we conjecture is
optimal in the diffusion limit.

\subsection*{\texorpdfstring{Notation.}{Notation}}

Let $\Naturals$ be the set of natural numbers $\{1,2,\ldots\}$,
let $\IntegersP=\{0,1,\allowbreak 2,\ldots\}$, let $\Reals$ be the set of real
numbers and let $\RealsP=\{x\in\Reals\dvtx x\geq0\}$.
Let $1_{\{\cdot\}}$ be the indicator function, where
$1_{\mathrm{true}}=1$ and $1_{\mathrm{false}}=0$. Let $x\mmin
y=\min(x,y)$, $x\mmax y=\max(x,y)$ and $[x]^+=x\mmax0$. When $x$
is a vector, the maximum is taken componentwise.

We will reserve bold letters for vectors in
$\Reals^\N$, where $N$ is the number of queues, for example,
$\bx=[x_n]_{1\leq n\leq N}$.
Superscripts on vectors are used to denote labels, not exponents,
except where otherwise noted;
thus, for example, $(\bx^0,\bx^1,\bx^2)$ refers to three arbitrary
vectors.
Let $\bZero$ be the vector of all 0s, and $\bOne$ be the
vector of all 1s. Use the norm
$|\bx|=\max_n|x_n|$. For vectors $\bu$ and $\bv$ and
functions $f\dvtx\Reals\to\Reals$, let
\[
\bu\bdot\bv= \sum_{{n}=1}^\N u_{n}v_{n},\qquad
\bu\bv= [u_n v_n]_{1\leq n\leq N}
\quad\mbox{and}\quad
f(\bu) = [ f(u_{n}) ]_{1\leq{n}\leq\N},
\]
and let matrix multiplication take precedence over dot product so that
\[
\bu\bdot A\bv= \sum_{n=1}^N u_n \Biggl(\sum_{m=1}^N A_{n m}
v_m\Biggr).
\]
Let $A^\T$ be the transpose of matrix $A$.
For a set $\sS\subset\Reals^N$, denote its convex hull by~$\langle
{\sS}\rangle$.

For a fixed $T>0$, and $I\in\Naturals$,
let $C^I(T)$ be the set of continuous functions $[0,T]\to\Reals^I$,
where $\Reals^I$ is equipped with
the norm $|x|={\max_i }|x_i|$.
Equip $C^I(T)$ with the norm
\[
\|f\| = \sup_{t\in[0,T]} |f(t)|.
\]
Let $d(f,g)=\|f-g\|$ be the metric induced by this norm.
For $E \subset C^I(T)$ and $f \in C^I(T)$, let
$d(f,E) = \inf\{ d(f,g) \dvtx g \in E\}$.
Define the modulus of continuity $\modcont_\delta(\cdot)$ by
\[
\modcont_\delta(f) = \sup_{|s-t|<\delta} | f(s)-f(t) |,
\]
where $s,t\in[0,T]$. Since $[0,T]$ is compact, each $f\in C^I(T)$
is uniformly continuous, therefore $\modcont_\delta(f)\tend0$
as $\delta\tend0$.\vspace*{-3pt}

\section{\texorpdfstring{Switched network model.}{Switched network
model}}
\label{sec:model}

We now introduce the switched network model. Section
\ref{sec:model.queue} describes the general system model, Section
\ref{sec:model.alg} specifies the class of scheduling policies we are
interested in and Section \ref{sec:model.stochastic} lists the
probabilistic assumptions about the arrival process that are needed
for the main theorems.\vspace*{-3pt}

\subsection{\texorpdfstring{Queueing dynamics.}{Queueing dynamics}}
\label{sec:model.queue}

Consider a collection of $N$ queues. Let time be discrete, indexed by
$\tau\in\{0,1,\ldots\}$. Let $Q_n(\tau)$ be the amount of work in
queue $n\in\{1,\ldots,N\}$ at time slot $\tau$.
Following our general notation for vectors, we write
$\bQ(\tau)$ for $[Q_n(\tau)]_{1\leq n\leq N}$. The initial queue sizes
are $\bQ(0)$.
Let $A_n(\tau)$ be the total amount of work arriving to
queue $n$, and $B_n(\tau)$ be the cumulative
potential service to queue $n$, up to time $\tau$,
with $\bA(0)=\bB(0)=\bZero$.\vadjust{\goodbreak}

We first define the queueing dynamics for a single-hop switched network.
Defining
$d\bA(\tau)=\bA(\tau+1)-\bA(\tau)$ and
$d\bB(\tau)=\bB(\tau+1)-\bB(\tau)$, the basic Lindley recursion
that we
will consider is
%
%
\begin{equation}
\label{eq:lindley1}
\bQ(\tau+1) = [ \bQ(\tau)-d\bB(\tau)]^+ + d\bA(\tau),
\end{equation}
where the \mbox{$[\cdot]^+$} is taken componentwise. The fundamental
``switched network'' constraint is that there is some finite set
$\sS\subset\RealsP^N$ such that
%
%
\begin{equation}
\label{eq:scheduling}
d\bB(\tau)\in\sS\qquad\mbox{for all $\tau$}.
\end{equation}
We will refer to $\bpi\in\sS$ as a schedule and $\sS$ as the set of
allowed schedules. In the applications in this paper, the schedule is
chosen based on current queue sizes,
which is why it is
natural to write the basic Lindley recursion as (\ref{eq:lindley1})
rather than the more standard $[\bQ(\tau)+d\bA(\tau)-d\bB(\tau)]^+$.

For the analyses in this paper it is useful to keep track of two other
quantities.
Let $Y_n(\tau)$ be the cumulative amount of idling at
queue $n$, defined by $\bY(0)=\bZero$ and
%
%
\begin{equation}
\label{eq:def.idling}
d\bY(\tau) =
[d\bB(\tau)-\bQ(\tau)]^+,
\end{equation}
where $d\bY(\tau)=\bY(\tau+1)-\bY(\tau)$. Then (\ref
{eq:lindley1}) can be rewritten
%
%
\begin{equation}
\label{eq:discrete.queue.singlehop}
\bQ(\tau) = \bQ(0) + \bA(\tau) - \bB(\tau) + \bY(\tau).
\end{equation}
Also, let $S_\bpi(\tau)$ be the cumulative time spent on
schedule $\bpi$ up to time $\tau$, so that
%
%
\begin{equation}
\label{eq:service}
\bB(\tau)=\sum_{\bpi\in\sS} S_\bpi(\tau) \bpi.
\end{equation}

For a multi-hop switched network, let $R\in\{0,1\}^{N\times N}$ be the
routing matrix, $R_{m n}=1$ if work served from queue $m$ is sent to
queue $n$ and $R_{m n}=0$ otherwise; if $R_{m n}=0$ for all $n$, then
work served from queue $m$ departs the network. For each $m$ we
require $R_{m n}=1$ for at most one $n$. (Tassiulas and Ephremides
\cite{tassiula1} described
a network model with routing choice, whereas we have restricted
ourselves to fixed routing for the sake of simplicity.)
The scheduling constraint (\ref{eq:scheduling}) is as before, the
definition of idling (\ref{eq:def.idling}) is as before and the
queueing dynamics are now defined by
\[
Q_n(\tau\!+\!1) = Q_n(\tau) + dA_n(\tau)
- \bigl(dB_n(\tau)-dY_n(\tau)\bigr)
+ \sum_m\!R_{m n}\bigl(dB_m(\tau)-dY_m(\tau)\bigr).
\]
Equivalently,
%
%
\begin{equation}
\label{eq:discrete.queue.multihop}
\bQ(\tau) = \bQ(0)
+ \bA(\tau)
- (I-R^\T)\bigl(\bB(\tau)-\bY(\tau)\bigr).
\end{equation}
Note that $\bA$ includes only exogenous arrivals to the network, not
internally routed traffic.
We will assume that routing is acyclic, that is, that work served
from
some queue $n$ never returns to queue $n$.
For example, Border Gateway Protocol (BGP)
utilized for routing in the internet is designed to be acyclic
\cite{BGP}. This implies that the
inverse $\tR=(I-R^\T)^{-1}$ exists; by considering the expansion\vadjust{\goodbreak}
$\tR=I+R^\T+(R^\T)^2+\cdots$ it is clear that $\tR_{m n}\in\{0,1\}$
for all $m$, $n$ and that $\tR_{m n}=1$ if work injected at queue $n$
eventually passes through $m$, and $\tR_{m n}=0$ otherwise. When $R=0$
we obtain a single-hop switched network.

A straightforward bound we shall need is
%
%
\begin{equation}
\label{eq:multihop.discrete.queue.bound}
Q_n(\tau) \leq Q_n(\tau')
+ A_n(\tau)-A_n(\tau')
+ \sum_m R_{m n}\bigl(B_m(\tau)-B_m(\tau')\bigr)
\end{equation}
for $\tau'\leq\tau$.

\subsection{\texorpdfstring{Scheduling policy.}{Scheduling policy}}
\label{sec:model.alg}

A policy that decides which
schedule to choose at each time slot $\tau\in\IntegersP$ is
called a \textit{scheduling policy}.
In this paper we will be interested in the max-weight scheduling
policy, introduced by Tassiulas and Ephremides \cite{tassiula1}. We
will refer to it as MW.

\subsubsection{\texorpdfstring{Max-weight policy for single-hop
network.}{Max-weight policy for single-hop
network}}
We describe the policy first for a single-hop network.
Let $\bQ(\tau)$ be the vector of queue sizes at
time~$\tau$.
Define the weight of a schedule $\bpi\in\sS$ to
be $\bpi\bdot\bQ(\tau)$. The MW policy then chooses\setcounter{footnote}{2}\footnote{%
There may be several schedules which jointly have the greatest
weight. To be concrete, we might specify some fixed numbering of
schedules and choose the highest-numbered maximum-weight
schedule. Alternatively, we might treat MW not as a policy per
se but as a constraint on the set of allowed sample paths.
For example, in a stochastic setting, we might allow $d\bB(\tau)$ to be
a random variable, measurable with respect to the underlying
probability space, satisfying (\ref{eq:maxweight}) for every
randomness. This permits ``break ties at random.''
For the analyses in this paper, it makes no difference
which of these two options is used.}
for time slot
$\tau$ a schedule $d\bB(\tau)$ with the greatest weight,
%
%
\begin{equation}
\label{eq:maxweight}
d\bB(\tau) \in\mathop{\argmax}_{\bpi\in\sS} \bpi\bdot\bQ
(\tau).
\end{equation}
This policy can be
generalized to choose a schedule which maximizes
$\bpi\bdot\bQ(\tau)^\alpha$, where the exponent is
taken componentwise for some $\alpha>0$; call this the MW-$\alpha$
policy.
More generally, one could choose a schedule
such that
%
%
\begin{equation}
\label{eq:discrete.mwm}
d\bB(\tau) \in\mathop{\argmax}_{\bpi\in\sS} \bpi\bdot f(\bQ
(\tau))
\end{equation}
for some function $f\dvtx\RealsP\to\RealsP$; call this the MW-$f$ policy.
It is assumed that~$f$ satisfies the following scale-invariance property:
\begin{assumption}
\label{cond.f}
Assume $f$ is differentiable and strictly increasing with
$f(0)=0$. Assume also that for any $\mathbf{q}\in\RealsP^\N$ and
$\bpi\in\sS$, with $m(\mathbf{q})=\break\max_{\brho\in\sS}\brho\bdot
f(\mathbf{q})$,
\[
\bpi\bdot f(\mathbf{q}) = m(\mathbf{q})
\quad\implies\quad
\bpi\bdot f(\kappa\mathbf{q}) = m(\kappa\mathbf{q})\qquad
\mbox{for all $\kappa\in\RealsP$}.
\]
\end{assumption}

This is satisfied by $f(x)=x^\alpha$, $\alpha>0$,
but it is not satisfied, for example, for an input-queued switch with
$f(x)=\log(1+x)$.

\subsubsection{\texorpdfstring{Max-weight policy for multi-hop
network.}{Max-weight policy for multi-hop
network}}

Now we define the multi-hop version of the MW-$f$ scheduling policy.
This policy chooses a
schedule~$d\bB(\tau)$ at time $\tau$ such that
\[
d\bB(\tau) \in\mathop{\argmax}_{\bpi\in\sS}
\bpi\cdot(I-R)f(\bQ(\tau)).
\]
Recall that matrix multiplication takes precedence over the $\cdot$
operator, so the $\argmax$ is of
$\bpi\cdot\{(I-R)f(\bQ(\tau))\}$;
note also that
\[
[R f(\bQ)]_n = \sum_m R_{n m} f(Q_m) = f([R\bQ]_n),
\]
where $[R\bQ]_n$ is the queue size at the first queue downstream of
$n$ (or $0$ if there is no queue downstream).
Thus
%
%
\begin{equation}
\label{eq:backpressure}
[(I-R)f(\bQ)]_n = f(Q_n) - f([R\bQ]_n).
\end{equation}
The difference $f(Q_n)-f([R\bQ]_n)$ is interpreted as the pressure to
send work from queue $n$ to the queue downstream of $n$; if the downstream
queue has more work in it than the upstream queue, then there is no pressure
to send work downstream. For this reason, it is also known as
\textit{backpressure} policy.

As before we will assume that $f$ satisfies a scale-invariance property,
the multi-hop equivalent of Assumption \ref{cond.f}:

\begin{assumption}
\label{multihop.cond.f}
Assume $f$ is differentiable and strictly increasing with
$f(0)=0$. Assume also that for any $\mathbf{q}\in\RealsP^\N$ and
$\bpi\in\sS$, with $m(\mathbf{q})=\break\max_{\brho\in\sS}\brho\bdot
(I-R)f(\mathbf{q})$,
\[
\bpi\bdot(I-R)f(\mathbf{q}) = m(\mathbf{q})
\quad\implies\quad
\bpi\bdot(I-R)f(\kappa\mathbf{q}) = m(\kappa\mathbf{q})
\qquad\mbox{for all $\kappa\in\RealsP$}.
\]
\end{assumption}

We further require that the scheduler always have the option of not
sending work downstream at any individual queue. Our Lyapunov
function, and indeed our whole fluid analysis in Section
\ref{sec:fluid.behavior.m},
rely on this assumption.

\begin{assumption}\label{cond.monotone}
For the multi-hop setting, assume that $\sS$ satisfies the following:
if $\bpi\in\sS$ is an allowed schedule,
and $\brho\in\RealsP^N$ is some other
vector with $\rho_n\in\{0,\pi_n\}$ for all $n$, then
$\brho\in\sS$.
\end{assumption}

In the rest of this paper, whenever we refer to a network running the
\mbox{MW-$f$} back-pressure policy, we mean that Assumptions
\ref{multihop.cond.f} and \ref{cond.monotone} are satisfied.

\subsection{\texorpdfstring{Stochastic model.}{Stochastic model}}
\label{sec:model.stochastic}

Some of the results in this paper are about fluid-scaled processes,
and others are about multiplicative state space collapse in the
diffusion scaling, and the different
results make different assumptions about the arrival process.
\begin{assumption}[(Assumptions for the fluid scale)]
\label{cond:fluid.stoch}
Let $\bA(\cdot)$ be a random process with stationary increments.
Assume it has a well-defined mean\vadjust{\goodbreak} arrival rate vector $\blambda$,
that is, assume $\lim_{\tau\tend\infty} A_n(\tau)/\tau$ exists almost
surely and is deterministic for every queue $1\leq n\leq N$, and define
%
%
\begin{equation}
\label{cond.arr}
\lambda_n = \lim_{\tau\tend\infty} \frac{1}{\tau} A_n(\tau).
\end{equation}
Assume there is a sequence of deviation terms $\delta_r\in\RealsP$,
$r\in\Naturals$, such that $\delta_r\tend0$ as $r\tend\infty$ and
\[
\Prob \biggl(\sup_{\tau\leq r} \frac{1}{r}| \bA(\tau
)-\blambda\tau| \geq
\delta_r \biggr)
\tend0 \qquad\mbox{as $r\tend\infty$}.
\]
\end{assumption}
\begin{assumption}[(Assumptions for multiplicative state space collapse)]
\label{cond:mssc.stoch}
Let $\bA^r(\cdot)$ be a sequence of random processes indexed by
$r\in\Naturals$. For each $r$, assume that $\bA^r$ has stationary
increments, and a well-defined mean arrival rate vector
$\blambda^r$, and that there is some limiting arrival rate vector
$\blambda$ such that\looseness=-1
\[
\blambda^r\tend\blambda\qquad\mbox{as $r\tend\infty$}.
\]\looseness=0
Assume there is a sequence of deviation terms
$\delta_z\in\RealsP$, $z\in\Naturals$, such that
$\delta_z\tend0$ as $z\tend\infty$ and
\[
z (\log z)^2 \Prob \biggl(\sup_{\tau\leq z} \frac{1}{z}
| \bA^r(\tau)-\blambda^r\tau| \geq\delta_z \biggr)
\tend0\qquad
\mbox{as $z\tend\infty$, uniformly in $r$}.
\]
\end{assumption}

If the arrival process is the same for all $r$, say
$\bA^r=\bA$ where $\bA$ has a~well-defined mean arrival
rate vector, then Assumption \ref{cond:mssc.stoch} reduces to
%
%
\begin{equation}
\label{eq:cond.arr}
\Prob \biggl(\sup_{\tau\leq r} \frac{1}{r}|
\bA(\tau)-\blambda\tau| \geq\delta_r \biggr)
=
o\biggl(\frac{1}{r(\log r)^2}\biggr),
\end{equation}
and it implies Assumption \ref{cond:fluid.stoch}.
For any arrival process with i.i.d. increments that are uniformly
bounded, that is, such that there is an $A^{\max}$ for which
\[
A_n(\tau+1)-A_n(\tau)
\in
[ 0, A^{\max}]\qquad
\mbox{for all $n$, $\tau$},
\]
equation (\ref{eq:cond.arr}) holds with $\delta_r = C\sqrt{\log
r}/\sqrt{r}$,
with choice of an appropriate constant $C$ that depends on $A^{\max}$,
by an application of concentration inequality by Azuma \cite
{azumainequality} and
Hoeffding \cite{hoeffdinginequality}.
More generally, it holds when the increments are not uniformly bounded but
instead satisfy a reasonable moment bound. For example, an application of
Doob's maximal inequality~\cite{doobmaximal} with bounded fourth moment
and $\delta_r = r^{-1/6}$ yields a stronger result than
(\ref{eq:cond.arr}); this can be used to show that a Poisson process
satisfies that equation. Furthermore
(\ref{eq:cond.arr}) holds for a much larger class of stationary
arrival processes beyond processes with i.i.d. increments, for
example, Markov modulated processes (see Dembo and Zeitouni \cite{DZ}).

\subsection{\texorpdfstring{Motivating example.}{Motivating example}}
\label{ssec:iq}

An internet router has several input ports and output ports. A data
transmission cable is attached to each of these ports. Packets arrive
at the input ports. The function of the router is to work out which\vadjust{\goodbreak}
output port each packet should go to, and to transfer packets to the
correct output ports. This last function is called \textit{switching}.
There are a~number of possible switch architectures; we will consider
the commercially popular input-queued switch architecture.

\begin{figure}

\includegraphics{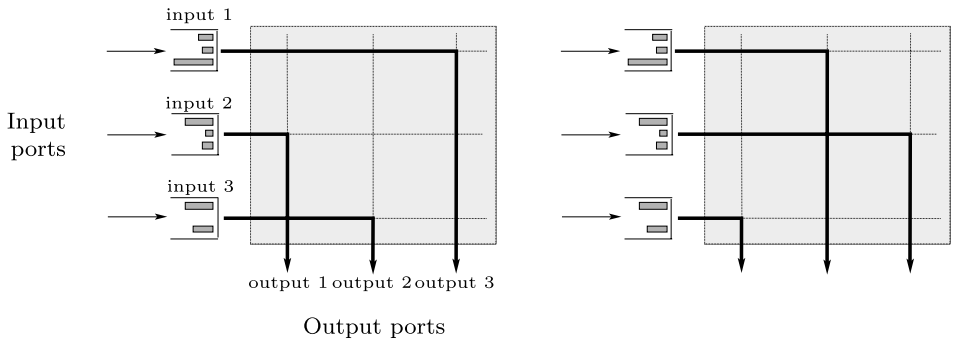}

\caption{An input-queued switch, and two example matching of inputs
to outputs.}
\label{fig:switch-iq}
\end{figure}

Figure \ref{fig:switch-iq} illustrates an input-queued switch with
three input ports and three output ports. Packets arriving at input
$i$ destined for output $j$ are stored at input port~$i$, in queue
$Q_{i,j}$; thus there are $N=9$ queues in total. (For this example
it is more natural to use
double indexing, e.g., $Q_{3,2}$, whereas for general switched networks
it is more natural to use single indexing, e.g., $Q_n$ for $1\leq n\leq N$.)

%
%

The switch operates in discrete time. In each time slot, the switch
fabric can transmit a number of packets from input ports to output
ports, subject to the two constraints that each input can transmit
at most one packet and that each output can receive at most one
packet. In other words, at each time slot the switch can choose a
\textit{matching} from inputs to outputs. The schedule
$\bpi\in\RealsP^{3\times3}$ is given by $\pi_{i,j}=1$ if
input port $i$ is matched to output port $j$ in a~given time slot,
and $\pi_{i,j}=0$ otherwise. Clearly $\bpi$ is a permutation
matrix, and the set $\sS$ of allowed schedules is the set of
$3\times3$ permutation matrices.

Figure \ref{fig:switch-iq} shows two possible matchings. In the
left-hand figure, the matching allows a packet to be transmitted
from input port 3 to output port 2, but since $Q_{3,2}$ is
empty, no packet is actually transmitted.

\section{\texorpdfstring{Related work.}{Related work}}
\label{sec:discussion}

Keslassy and McKeown \cite{KM} found from extensive simulations of an
input-queued switch that the average queueing delay is
different under MW-$\alpha$ policies for different values
of $\alpha> 0$. They conjecture:

\begin{conjecture}
\label{conj:keslassy}
For an input-queued switch running the MW-$\alpha$ policy,
the average queueing delay decreases as $\alpha$ decreases.
\end{conjecture}

Though our work is motivated by the desire to establish
Conjecture \ref{conj:keslassy}, we have not been able to prove it.
But whereas the two main analytic approaches that have been employed
in the literature yield results for the\vadjust{\goodbreak} input-queued switch that are
insensitive to $\alpha$, our result about multiplicative state
space
collapse \textit{is} sensitive, as shown in Section \ref{sec:optimal}.
We speculate that our result might
eventually form part of a proof of the conjecture.

The two main analytic approaches that have been employed in the
literature are stability analysis and heavy traffic analysis. In
stability analysis, one calculates the set of arrival rates for which
a policy is stable (in the sense of
\cite{tassiula1,MAW,daibala,KM,dev,andrewsetal}). All the prior work in this context leads
to the conclusion that MW-$\alpha$ has the optimal stability region,
regardless of $\alpha$.

In heavy traffic analysis, one looks at queue size behavior under
a diffusion (or heavy traffic) scaling. This regime was
first described by Kingman~\cite{kingmanht}; since then a
substantial body of theory has developed, and modern treatments
can be found in \cite{mike2,bramson,williams,whittspl}.
Stolyar has studied MW-$\alpha$ for a~generalized
switch model in the diffusion scaling, and
obtained a complete characterization of the diffusion approximation
for the queue size process, \textit{under a condition known
as} ``\textit{complete resource pooling}.'' This condition effectively
requires that a clever scheduling policy be able to balance
work between all the heavily loaded queues. Stolyar \cite{stolyar}
showed in a remarkable paper that the limiting queue size lives in a
one-dimensional state space. Operationally, this means that all one
needs to keep track of is the one-dimensional total amount of work in
the system (called the \textit{workload}), and at any point in time one
can assume that the individual queues have all been balanced.
Dai and Lin \cite{LD05,LD08} have established that similar result holds
in the more general setting of a stochastic processing network.

Under the complete resource pooling condition, the results in
\cite{stolyar,LD05,LD08} imply that the performance of MW-$\alpha$ in
an input-queued switch is always optimal (in the diffusion scaling)
regardless of the value of $\alpha>0$. Therefore these results do not
help in addressing Conjecture \ref{conj:keslassy}. This is our
motivation for studying switched networks in the absence of complete
resource pooling. Technically, the lifting map for a
critically-loaded input-queued switch is degenerate and insensitive to
$\alpha$ under complete resource pooling, but it is sensitive to
$\alpha$ otherwise.

We prove multiplicative state space collapse, following the method of
Bramson~\cite{bramson}.
The complement of Bramson's work is by Williams \cite{williams},
and consists of proving a diffusion approximation, using an appropriate
invariance principle along with the multiplicative state space
collapse. We do not carry out this complementary aspect.
Stolyar \cite{stolyar} and Dai and Lin \cite{LD05,LD08}  have proved
the diffusion
approximation under the complete resource pooling condition, and
Kang and Williams \cite{KW} have made progress toward it in the case
without complete
resource pooling, for an input-queued switch under the MW-$1$ policy.

Whereas in heavy traffic models of other systems
\cite{mike2,bramson,williams,stolyar} the lifting
map from workloads to queue sizes is linear,
we find instead that it is nonlinear---in fact it can be
expressed as the solution to an optimization problem.
The objective function of the problem is a natural\vadjust{\goodbreak} generalization of
the Lyapunov function introduced by Tassiulas and Ephremides \cite{tassiula1}
for proving stability of the MW-$1$ policy; the constraints
of the problem are closely
linked to the canonical representation of workload
identified by Harrison \cite{harrisoncanonical}.
The objective function for MW-$\alpha$ depends on $\alpha$, and this
hints that the performance measures might also depend on $\alpha$.

Finally, we take note of two related results. First, in \cite
{shahwischikinfocom}
we have reported some results about a critically loaded
input-queued switch without a complete resource pooling condition.
Second, a sequence of works by Kelly and Williams \cite{KW04} and Kang
et al. \cite{kellywilliamsssc} has
resulted in a diffusion approximation for a bandwidth sharing
network model operating under proportionally fair rate allocation,
assuming a technical ``local traffic'' condition,
but without assuming complete resource pooling.
They show
that the resulting diffusion approximation model has a product form
stationary distribution.

\section{\texorpdfstring{The fluid approximation.}{The fluid
approximation}}
\label{sec:fluid}

This section introduces the fluid model and establishes it as
an approximation to a fluid-scaled descriptor
of the switched network. Intuitively, the fluid model describes
the dynamics of the system at the ``rate'' level rather than at
finer granularity. The reader is referred to a recent
monography by Bramson \cite{Bramson08} and lecture notes by
Dai \cite{dailec} for a~detailed account of fluid approximation
for multiclass queueing networks.
In Section \ref{sec:fluid.model} we specify the fluid model,
in Section \ref{sec:fluid.result} we state the main result
and in Section \ref{sec:fluid.proof} we prove it.

\subsection{\texorpdfstring{Definition of fluid model.}{Definition of fluid
model}}
\label{sec:fluid.model}
Let time be measured by $t\in[0,T]$ for some fixed $T>0$.
Let $\mathbf{q}$, $\mathbf{a}$ and $\mathbf{y}$ all be continuous
functions mapping
$[0,T]$ into~$\RealsP^\N$, and let $s=(s_\bpi)_{\bpi\in\sS}$ be a
collection of continuous functions mapping $[0,T]$ into~$\RealsP$.
Let $x(\cdot)=(\mathbf{q}(\cdot),\mathbf{a}(\cdot),\mathbf
{y}(\cdot),s(\cdot))$. This lies
in $C^I(T)$ where $I=3N+|\sS|$.
The definition below requires these functions to be absolutely
continuous; such functions are differentiable almost everywhere, and
the time instants where they are differentiable are called ``regular
times.''
Any equations we write involving derivatives are taken to apply only at
regular times.\looseness=-1

\begin{definition}[(Fluid model solution for single-hop network)]
\label{def:fms}
Let $f\dvtx\allowbreak\RealsP\to\RealsP$ satisfy Assumption \ref{cond.f}.
Say that $x(\cdot) $ is a fluid model solution for a
single-hop switched network with arrival rate $\blambda\in\RealsP^N$
operating under the MW-$f$ policy if
it satisfies equations
(\ref{eq:fluid.queue})--(\ref{eq:fluid.mwm}) below. Write $\FMS$
for the set of all such $x \in C^I(T)$. Additionally, define
\begin{eqnarray*}
\FMS_K &=& \{ x\in\FMS\dvtx |\mathbf{q}(0)|\leq K \},\\
\FMS(\mathbf{q}_0) &=& \{ x\in\FMS\dvtx \mathbf{q}(0)=\mathbf
{q}_0 \}.
\end{eqnarray*}
The equations are:
%
%
\begin{eqnarray}\label{eq:fluid.queue}
\mathbf{q}(t) &=& \mathbf{q}(0) + \mathbf{a}(t) - \sum_{\bpi}
s_\bpi(t)\bpi+ \mathbf{y}(t);
\\
\label{eq:fluid.arrivals}
\mathbf{a}(t) &=& \blambda t;
\\
\label{eq:fluid.busy}
\sum_{\bpi\in\sS} s_{\bpi}(t) &=& t;
\\
\label{ineq:fluid.littleidle}
\mathbf{y}(t)&\leq&\sum_{\bpi\in\sS}s_{\bpi}(t)\bpi;
\end{eqnarray}
\begin{eqnarray}\label{eq:fluid.increasing}\qquad
&&\mbox{each $s_\bpi(\cdot)$ and $y_{n}(\cdot)$ is increasing
(not necessarily strictly increasing)};
\\
\label{eq:fluid.abscont}
&&\mbox{all the components of $x(\cdot)$ are absolutely continuous};
\\
\label{eq:fluid.idling}
&&\mbox{for regular times $t$, all }n\qquad \dot{y}_{n}(t) = 0
\qquad\mbox{if $q_{n}(t)>0$};
\\
\label{eq:fluid.mwm}
&&\mbox{for regular times $t$, all }\bpi\in\sS\nonumber\\[-8pt]\\[-8pt]
&&\qquad\dot{s}_\bpi
(t) = 0\qquad\mbox{if }
\bpi\bdot f(\mathbf{q}(t))<\max_{\brho\in\sS} \brho\bdot
f(\mathbf{q}(t)).\nonumber
\end{eqnarray}
\end{definition}

Here, $\mathbf{q}(t)$ represents the vector of queue sizes at time $t$,
$\mathbf{a}(t)$ represents the cumulative arrivals up to time $t$,
$\mathbf{y}(t)$ represents the cumulative idleness up to time $t$
and $s_\bpi(t)$ represents
the total amount of time spent on schedule~$\bpi$ up to time $t$.
The equation in (\ref{eq:fluid.queue}) is the continuous analog of
(\ref{eq:discrete.queue.singlehop}) combined with
(\ref{eq:service}), and the inequality is the analog of
the single-hop version of (\ref{eq:multihop.discrete.queue.bound}).
Equation (\ref{eq:fluid.arrivals}) represents an assumption about the
arrival process, related to (\ref{cond.arr}).
Equation
(\ref{eq:fluid.busy}) says that the scheduling policy must choose
some schedule at every timestep.
Both (\ref{ineq:fluid.littleidle}) and (\ref{eq:fluid.idling})
derive from the definition of idling, (\ref{eq:def.idling}).
Equation (\ref{eq:fluid.mwm}) is the continuous analog of (\ref
{eq:discrete.mwm}).
\begin{definition}[(Fluid model solution for multi-hop network)]
Let $f\dvtx\allowbreak\RealsP\to\RealsP$ satisfy Assumption \ref{multihop.cond.f},
and let $\sS$ satisfy Assumption \ref{cond.monotone}.
Say that $x(\cdot)$ is a fluid model solution for a multi-hop switched
network operating under the MW-$f$ policy if
it satisfies equations (\ref{eq:fluid.arrivals})--(\ref{eq:fluid.idling}),
and additionally~(\ref{eq:multihop.fluid.queue}) and~(\ref{eq:multihop.fluid.mwm})
below. Let $\FMSm$ be the set of all such $x\in C^I(T)$.
Also, let $\FMSm_K$ and $\FMSm(\mathbf{q}_0)$ be defined analogously
to the
single-hop case.
The extra equations are
%
%
\begin{equation}
\label{eq:multihop.fluid.queue}
\mathbf{q}(t) = \mathbf{q}(0) + \mathbf{a}(t) - (I-R^\T)
\biggl(\sum_{\bpi} s_\bpi(t)\bpi- \mathbf{y}(t)\biggr)
\end{equation}
and
%
%
\begin{eqnarray}
\label{eq:multihop.fluid.mwm}
&&\mbox{for all regular times $t$, all $\bpi\in\sS$}
\nonumber\\[-8pt]\\[-8pt]
&&\qquad \dot{s}_\bpi(t)=0\qquad\mbox{if }
\bpi\bdot(I-R)f(\mathbf{q}(t))<\max_{\brho\in\sS} \brho\bdot
(I-R)f(\mathbf{q}(t)).\nonumber
\end{eqnarray}
\end{definition}

When we refer to ``fluid model solutions for any scheduling policy,'' we
mean processes $x(\cdot)\in C^I(T)$ satisfying
(\ref{eq:fluid.queue}) to (\ref{eq:fluid.idling}) in the single-hop
case,
or satisfying (\ref{eq:fluid.arrivals}) to (\ref{eq:fluid.idling})
plus (\ref{eq:multihop.fluid.queue}) in the multi-hop case.

\subsection{\texorpdfstring{Main fluid model result.}{Main fluid model
result}}
\label{sec:fluid.result}

The development in this section follows the general pattern of
Bramson \cite{bramson}. There is, however, a difference in
presentation that is
worth noting. The main result of this section, Theorem
\ref{thm:fluid}, is a~general purpose sample path-wise result: it does
not make any probabilistic claim nor does it depend on any
probabilistic assumptions. It can be applied to a switched network
with stochastic arrivals in two ways: to obtain a result about fluid
approximations (Corollary \ref{cor:fluid}), and to obtain a~result
about multiplicative state space collapse (Section \ref{sec:ht}).

We start by defining the fluid scaling. Consider a switched network of
the type described in Section \ref{sec:model.queue} running a
scheduling policy of the type described in Section
\ref{sec:model.alg}. Write
$X(\tau)=(\bQ(\tau),\bA(\tau),\bY(\tau),S(\tau))$, $\tau\in
\IntegersP$,
to denote its sample path. Given a scaling parameter $z\geq1$,
define the fluid-scaled sample path
$\fx(t)=(\fbq(t),\fba(t),\fby(t),\fls(t))$ for $t\in\RealsP$ by
\begin{eqnarray*}
\fbq(t) &=& \bQ(z t)/z,\qquad
\fba(t) = \bA(zt)/z,\\
\fby(t) &=& \bY(zt)/z,\qquad
\fls_\bpi(t) = S_\bpi(zt)/z
\end{eqnarray*}
after extending the domain of $X(\cdot)$ to
$\RealsP$ by linear interpolation in each interval $(\tau,\tau+1)$.
In this section we are interested in the evolution of $\fx(t)$ over
$t\in[0,T]$ for some fixed $T>0$, therefore we take $\fx$ to lie in
$C^I(T)$ with $I=3\N+|\sS|$.

The following theorem concerns uniform convergence of a set of
fluid-scaled sample paths. Every fluid-scaled sample path is assumed
to relate to some (unscaled) switched network, and all the switched
networks are assumed to have the same network data, that is, the same
number of queues~$N$, the same set of allowed schedules $\sS$, the
same routing matrix $R$, and the same scheduling policy.

The convergence is indexed by a parameter $j$ lying in some totally
ordered countable set. For Corollary \ref{cor:fluid} we will use
$j\in\Naturals$, and for Section \ref{sec:ht} we will use a subset of
$\Naturals\times\Naturals$ as the index set. We are purposefully using
the symbol~$j$ here as an index, rather than the $r$ used elsewhere,
to remind the reader that the index set is interpreted differently in
different results.

\begin{theorem}
\label{thm:fluid}
Let $\cX$ be the set of all possible sample paths for single-hop
switched networks with the network data specified above,
running the MW-$f$\vadjust{\goodbreak} scheduling policy, where $f$ satisfies Assumption
\ref{cond.f}.
Fix $K>0$ and $\blambda\in\RealsP^N$. Let there be sequences
$\varepsilon_j \in\RealsP$ and $\blambda^j\in\RealsP^N$,
indexed by $j$ in some totally ordered countable set, such
that
%
%
\begin{equation}
\label{cond.rateconv}
\varepsilon_j \tend0 \quad\mbox{and}\quad \blambda^j
\tend\blambda,\qquad
\mbox{as $j\tend\infty$}.
\end{equation}
Consider a sequence of subsets
$G_j\subset C^I(T)\times[1,\infty)$ which
satisfy the following: for every $(\fx,z)\in G_j$ there is some
unscaled sample path $X\in\cX$ such that $\fx$ is the fluid-scaled
version\vadjust{\goodbreak}
of $X$ with scaling parameter $z$
(here $z$ is permitted to be a function of $X$);
and furthermore
%
%
\begin{eqnarray}
\label{eq:fluid.scale}
\inf\{z \dvtx (\fx,z)\in G_j\} &\tend&\infty
\qquad\mbox{as $j\tend\infty$,}\\
\label{eq:fluid.arr}
\sup_{(\fx,z)\in G_j} \sup_{t\in[0,T]} | \fba(t)-\blambda
^jt |
&\leq&\varepsilon_j \qquad\mbox{for all $j$}
\end{eqnarray}
and
\begin{equation}
\label{eq:fluid.qbound}
\sup_{(\fx,z)\in G_j} | \fbq(0)|\leq K
\qquad\mbox{for all $j$.}
\end{equation}
Then
%
%
\begin{equation}\label{eq:fluid1}
\sup_{(\fx,z)\in G_j} d(\fx,\FMS_K) \tend0 \qquad\mbox
{as $j\tend\infty$}.
\end{equation}
Furthermore, fix $\mathbf{q}_0\in\RealsP^N$
and a sequence $\varepsilon'_j\tend0$,
and assume that the sets~$G_j$ also satisfy
%
%
\begin{equation}
\label{eq:fluid.initq}
\sup_{(\fx,z)\in G_j} |\fbq(0)-\mathbf{q}_0| \leq
\varepsilon'_j\qquad
\mbox{for all $j$}.
\end{equation}
Then
%
%
\begin{equation}\label{eq:fluid2}
\sup_{(\fx,z)\in G_j} d(\fx,\FMS(\mathbf{q}_0)) \tend0
\qquad\mbox{as
$j\tend\infty$}.
\end{equation}
Equivalent results to (\ref{eq:fluid1}) and (\ref{eq:fluid2}) apply to
multi-hop switched networks, with references to $\FMS$ replaced by
$\FMSm$ and the set $\cX$ modified to refer to multi-hop networks
running the MW-$f$ scheduling policy where $f$ satisfies Assumption
\ref{multihop.cond.f} and $\sS$ satisfies Assumption \ref{cond.monotone}.
\end{theorem}

The above theorem as stated applies to the MW-$f$ scheduling policy,
but it is clear from the proof that a corresponding limit result
holds, relating sample paths of \textit{any} scheduling policy to fluid
models defined by equations~(\ref{eq:fluid.queue})--(\ref{eq:fluid.idling}).

The following corollary is a straightforward application of Theorem
\ref{thm:fluid}. It specializes the theorem to the case
of a single random system $X$, and the sequence of fluid-scaled versions
indexed by $r\in\Naturals$ where the $r$th version uses scaling
parameter $r$.
The arrival process is assumed to satisfy certain stochastic assumptions.
This corollary is useful when studying
the behavior of a single switched network with random arrivals, over
long timescales.
\begin{corollary}
\label{cor:fluid}
Consider a single-hop switched network as described in Section
\ref{sec:model.queue}, running the MW-$f$ policy as described in
Section \ref{sec:model.alg} where~$f$ satisfies Assumption \ref{cond.f}.
Let the arrival process $\bA(\cdot)$ satisfy Assumption~\ref
{cond:fluid.stoch},
and let the initial queue size $\bQ(0)\in\RealsP^N$ be random.
For $r\in\Naturals$, let
\begin{eqnarray*}
\fbq^r(t) &=& \bQ(rt)/r, \qquad\fba^r(t) = \bA(rt)/r,\\
\fby^r(t) &=& \bY(rt)/r, \qquad\fls_\bpi^r(t) = S_\bpi(rt)/r,
\end{eqnarray*}
and let $\fx^r(t)=(\fbq^r(t),\fba^r(t),\fby^r(t),\fls^r(t))$, for
$t\in[0,T]$ where $T>0$ is some fixed time horizon.
Then for any $\delta>0$
\[
\Prob\bigl(d(\fx^r(\cdot),\FMS(\bZero))<\delta
\bigr)\tend1
\qquad\mbox{as $r\tend\infty$}.
\]
The same conclusion holds for a multi-hop switched network running the
MW-$f$ back-pressure policy where $f$ satisfies Assumption
\ref{multihop.cond.f} and $\sS$ satisfies Assumption~\ref{cond.monotone},
with $\FMS$ replaced by $\FMSm$.
\end{corollary}
\begin{pf}
First define the event
$E_r$ by
\[
E_r = \biggl\{ \sup_{\tau\leq r}\frac{1}{r}
| \bA(\tau)-\blambda\tau|<\delta_r
\quad\mbox{and}\quad
|\bQ(0)|\leq\sqrt{r}
\biggr\},
\]
where $\blambda$ and $\delta_r$ are as in Assumption
\ref{cond:fluid.stoch}. By this we mean that $E_r$ is a~subset of the
probability sample space, and we write $X(\cdot)(\omega)$ etc. for
$\omega\in E_r$ to emphasize the dependence on $E_r$.

We will apply
Theorem \ref{thm:fluid} with index set $j\equiv r\in\Naturals$ to the
sequence of sets
\[
G_j \equiv G_r = \{ (\fx^r(\cdot)(\omega),r)
\dvtx
\omega\in E_r \}.
\]
In order to apply the theorem we will pick constants as follows. Let
$K=1$, let $\blambda$ be as in Assumption \ref{cond:fluid.stoch},
$\blambda^j=\blambda$ for all $j$,
$\varepsilon_j\equiv\varepsilon_r=T\delta_r$ where $\delta_r$ is
as in
Assumption~\ref{cond:fluid.stoch}, $\mathbf{q}_0=\bZero$ and
$\varepsilon'_r=1/\sqrt{r}$.
We now need to verify the conditions of Theorem
\ref{thm:fluid}.
Equation (\ref{cond.rateconv}) holds by the choice of
$\blambda^j$ and by Assumption \ref{cond:fluid.stoch}.
Equation (\ref{eq:fluid.scale}) holds automatically by choice of
$G_j$.
To see that (\ref{eq:fluid.arr}) holds, rewrite event $E_r$ in terms
of the fluid scaled arrival process~$\fba^r$ to see
\[
\sup_{t\in[0,T]} | \fba^r(t)(\omega)-\blambda t| <
T\delta_r\qquad
\mbox{for all $r$ and $\omega\in E_r$},
\]
which implies (\ref{eq:fluid.arr}); likewise for (\ref{eq:fluid.qbound})
and (\ref{eq:fluid.initq}). We conclude that (\ref{eq:fluid2}) holds.
It may be rewritten in terms of $E_r$ as
%
%
\begin{equation}
\label{eq:fluid.corr}
\sup_{\omega\in E_r} d(\fx^r(\cdot)(\omega),\FMS(\bZero
))\tend0\qquad
\mbox{as $r\tend\infty$.}
\end{equation}

We next argue that $\Prob(E_r)\tend1$ as $r\tend\infty$. The event
$E_r$ is
the intersection of two events, one concerning arrivals and the other
concerning initial queue size. The probability of the former $\tend1$
as $r\tend\infty$ by Assumption \ref{cond:fluid.stoch}.
For the latter, $\Prob(|\bQ(0)|\leq\sqrt{r})\tend1$ as $r\tend
\infty$
since $\bQ(0)$ is assumed not to be infinite. Therefore
$\Prob(E_r)\tend1$. Combining this with (\ref{eq:fluid.corr}) gives
the desired result for single-hop networks.
The multi-hop version follows similarly.
\end{pf}

\subsection{\texorpdfstring{Proof of Theorem \protect\ref
{thm:fluid}.}{Proof of Theorem 4.3}}
\label{sec:fluid.proof}

We shall present the proof of Theorem \ref{thm:fluid} for a single-hop
network in detail followed by main ideas required to extend it to
multi-hop networks.

\subsubsection{\texorpdfstring{Cluster points.}{Cluster points}}

Here we are interested in convergence in $C^I(T)$, where $I=3N+|\sS|$
and $T>0$ is fixed, equipped with the norm \mbox{$\|\cdot\|$}.
The appropriate concept for proving convergence is
\textit{cluster points}. Consider any metric space $E$ with metric $d$
and a sequence $(E_1,E_2,\ldots)$ of subsets of $E$. Say that $x\in E$
is a \textit{cluster point} of the sequence if $\liminf_{j\tend\infty}
d(x,E_j)=0$ where $d(x, E_j) = \inf\{d(x,y) \dvtx y \in E_j\}$.
\begin{proposition}[{[Cluster points in $C^I(T)$]}\footnote{Taken from
Bramson \cite{bramson}, Proposition 4.1.}]
\label{prop:tight}
Given $K>0$, $A>0$ and a sequence $B_j\tend0$, let
\[
K_j = \{ x\in C^I(T) \dvtx |x(0)|\leq K \mbox{ and }
\modcont_\delta(x)\leq A\delta+B_j \mbox{ for all $\delta>0$}
\},
\]
and consider a sequence $(E_1,E_2,\ldots)$ of subsets of $C^I(T)$ for
which $E_j\subset K_j$. Then $\sup_{y\in E_j} d(y,\clust)\tend0$ as
$j\tend\infty$, where $\clust$ is the set of cluster points of
$(E_1,E_2,\ldots)$.
\end{proposition}

\subsubsection{\texorpdfstring{Proof of Theorem \protect\ref
{thm:fluid}.}{Proof of Theorem 4.3}}

Let $E_j=\{\fx\dvtx (\fx,z)\in G_j\}$. Lemma \ref{lem:fluid.tight} below
shows that $E_j\subset K_j$, with $K_j$ as defined in Proposition
\ref{prop:tight} for appropriate constants $K$, $A$ and
$B_j$. By applying that proposition,
\[
\sup_{\fx\in E_j} d(\fx,\clust)\tend0\qquad
\mbox{as $j\tend\infty$},
\]
where $\clust$ is the set of cluster points of the $E_j$ sequence.
Lemma \ref{lem:limitpoints} below shows that all cluster points of
the $E_j$ sequence satisfy the fluid model equations. Every cluster
point $x$
must also satisfy $|\mathbf{q}(0)|\leq K$, by
(\ref{eq:fluid.qbound}). Therefore
\[
\sup_{\fx\in E_j} d(\fx,\FMS_K)\tend0
\qquad\mbox{as $j\tend\infty$}.
\]
If in addition (\ref{eq:fluid.initq}) holds, then every cluster point
$x$ must also satisfy $\mathbf{q}(0)=\mathbf{q}_0$. Therefore
\[
\sup_{\fx\in E_j} d(\fx,\FMS(\mathbf{q}_0))\tend0
\qquad\mbox{as $j\tend\infty$}.
\]

\begin{lemma}[(Tightness of fluid scaling)]
\label{lem:fluid.tight}
Let $K$ and $G_j$ be as in Theorem~\ref{thm:fluid}. Then
there exist a constant $A>0$ and a sequence $B_j\tend0$ such that
for every $(\fx,z)\in G_j$, $|\fx(0)|\leq K$ and
\[
| \fx(u)-\fx(t)|
\leq A |u-t| + B_j\qquad
\mbox{for all $0\leq t,u\leq T$.}
\]
\end{lemma}
\begin{pf}
Consider $(\fx,z)\in G_j$,
where $\fx=(\fbq,\fba,\fby,\fls)$.
As per the definitions in Section
\ref{sec:model.queue}, the only nonzero component
of $\fx(0)$ is $\fbq(0)$ and $|\fbq(0)|\leq K$ by choice of $G_j$,
hence $|\fx(0)|\leq K$. For the second inequality,
without loss of generality pick any $0\leq t<u\leq T$, and let us now
look at each component of $|\fx(u)-\fx(t)|$ in turn.\vadjust{\goodbreak}

For arrivals, let $\lambda^{\max}=\sup_j|\blambda^j|$; this is
finite by the assumption that $\blambda^j\tend\blambda$ in Theorem
\ref{thm:fluid}.
Then for $(\fx,z)\in G_j$,
\begin{eqnarray*}
| \fba(u)-\fba(t) |
&\leq&
| \fba(u)-\blambda^j u |
+ | \blambda^j (u-t)|
+ | \fba(t)-\blambda^j t | \\[-2pt]
&\leq&
2 \varepsilon_j + |\blambda^j|(u-t) \qquad\mbox{by (\ref
{eq:fluid.arr})}\\[-2pt]
&\leq&
2 \varepsilon_j + \lambda^{\max} (u-t).
\end{eqnarray*}
For idling, let $S^{\max}=\max_{\bpi\in\sS} \max_{n}
\pi_{n}$. This is the maximum amount of service that can be offered to
any queue per unit time, and it must be finite since $|\sS|$ is
finite. Then, based on (\ref{eq:def.idling}),
\begin{eqnarray*}
| \fy_n(u)-\fy_n(t) |
&\leq&
(u-t) S^{\max} + 2 S^{\max}/z \\[-2pt]
&\leq&(u-t) S^{\max} + 2 S^{\max}/z_j^{\min},
\end{eqnarray*}
where $z_j^{\min}=\inf\{z \dvtx (\fx,z)\in G_j\}$.
For service, let $S_\bpi(\cdot)$ be the unscaled process that
corresponds to $\fls_\bpi(\cdot)$;
since $S_\bpi(\cdot)$ is increasing and since a schedule
must be chosen not more than once every time slot,
\[
| \fls_\bpi(u) - \fls_\bpi(t) |
\leq
\frac{1}{z} \bigl( S_\bpi(\lceil z u\rceil) - S_\bpi(\lfloor z
t\rfloor)
\bigr)
\leq(u-t) + 2/z
\leq(u-t) + 2/z_j^{\min}.
\]
For queue size, note that (\ref{eq:discrete.queue.singlehop}) carries
through to
the fluid model scaling, that is,
\[
\fbq(t) = \fbq(0) + \fba(t) - \sum_\bpi\fls_\bpi(t)\bpi+ \fby(t),
\]
thus
\begin{eqnarray*}
&&| \fq_n(u)-\fq_n(t) |\\[-2pt]
&&\qquad\leq
| \fa_n(u)-\fa_n(t) |
+ \sum_\bpi\pi_n | \fls_\bpi(u) - \fls_\bpi(t) |
+ | \fy_n(u)-\fy_n(t) |\\[-2pt]
&&\qquad\leq
(u-t) ( \lambda^{\max} + |\sS| S^{\max} + S^{\max} ) +
(2|\sS| S^{\max} + 2S^{\max} )/z_j^{\min} + 2
\varepsilon_j.
\end{eqnarray*}
Putting all these together,
%
%
\begin{equation}\label{eq:tight.cont.bound0}
|\fx(u)-\fx(t)| \leq A (u-t) + B_j,\vadjust{\goodbreak}
\end{equation}
where the constants are
\[
A = (N+1)\lambda^{\max} + 2N S^{\max} + |\sS| + N|\sS|
S^{\max}
\]
and
\[
B_j = (4N S^{\max} + 2|\sS| + 2N|\sS|
S^{\max}) \frac{1}{z_j^{\min}}
+ (2+2N)\varepsilon_j.
\]
By the assumptions of Theorem \ref{thm:fluid}, $\varepsilon_j\tend0$
and $z_j^{\min}\tend\infty$ as $j\tend\infty$, thus $B_j\tend0$
as required.
\end{pf}

\begin{lemma}[(Dynamics at cluster points)]
\label{lem:limitpoints}
Make the same assumptions as Theorem \ref{thm:fluid}, and let
$E_j = \{\fx\dvtx (\fx, z) \in G_j\}$. Then $x \in\FMS_K$
if $x$ is a~cluster point of the $E_j$ sequence.\vadjust{\goodbreak}
\end{lemma}
\begin{pf}
From Lemma \ref{lem:fluid.tight} and Proposition
\ref{prop:tight}, it follows that\vspace*{1pt}\break
$\limsup_{\fx\in E_j} d(\fx$, $\clust) \tend0$ as $j\tend\infty$
where $\clust$ is the set of cluster points of the sequence $E_j$. Let
$x$ be
one such cluster point. That is,\vspace*{1pt} there exists a~subsequence $j_k$
and a collection $\fx^{j_k}\in E_{j_k}$ such that $\fx^{j_k}\tend x$.
It easily\vspace*{1pt} follows that $|x(0)|\leq K$ since $|\fx^{j_k}(0)| \leq K$
for all $\fx^{j_k} \in E_{j_k}$ as argued in Lemma~\ref{lem:fluid.tight}.
Using this, we wish to establish that $x$ satisfies all the
fluid model equations to conclude $x \in\FMS_K$. For convenience,
we shall omit the subscript $k$ in the rest of the proof; that
is, we shall use $j$ in place of $j_k$ and $j\tend\infty$.\vspace*{8pt}

\textit{Proof of} (\ref{eq:fluid.queue}), (\ref{eq:fluid.busy}),
(\ref{eq:fluid.increasing}).\quad
The discrete (unscaled) system satisfies these properties; therefore
the scaled systems $\fx^j$ do too. Taking the limit yields the fluid
equations.\vspace*{8pt}

\textit{Proof of} (\ref{ineq:fluid.littleidle}).\quad
In (\ref{eq:def.idling}), $d\bB(\tau)$ and $\bQ(\tau)$ are both
nonnegative (component-wise), hence $d\bY(\tau)\leq d\bB(\tau)$
for all $\tau$. Summing up
over $\tau$, we see the discrete (unscaled) system satisfies
the equivalent of (\ref{ineq:fluid.littleidle}), so as above we obtain
the fluid equation.\vspace*{8pt}

\textit{Proof of} (\ref{eq:fluid.arrivals}).\quad
Observe that
\[
\sup_{t\in[0,T]} | \mathbf{a}(t)-\blambda t |
\leq
\sup_{t\in[0,T]} | \mathbf{a}(t)-\fba^j(t) |
+
\sup_{t\in[0,T]} | \fba^j(t)-\blambda^j t |
+
T|\blambda^j-\blambda|.
\]
Each term converges to 0 as $j\tend\infty$: the first because
$\fx^j\tend x$, the second because $\fx^j\in E_j$ so
the deviation in $\fba^j$ is bounded by $\varepsilon_j$ and
$\varepsilon_j\to0$
and the third because $\blambda^j\tend\blambda$. Since the left-hand side
does not depend on $j$, it must be that $\mathbf{a}(t)=\blambda
t$.\vspace*{8pt}

\textit{Proof of} (\ref{eq:fluid.abscont}).\quad
In Lemma \ref{lem:fluid.tight} we found constants $A$ and $B_j$
such that for all $\fx\in E_j$
\[
| \fx(u) - \fx(t) | \leq A |u-t| + B_j
\]
with $B_j \to0$ as $j\to\infty$. Taking the limit
of $|\fx^j(u)-\fx^j(t)|$
as $j\tend\infty$,
we find that $|x(u)-x(t)|\leq A|u-t|$; that is, $x$ is (globally)
Lipschitz continuous (of order $1$ with respect to
the appropriate metric as defined earlier). This
immediately implies that $x$ is absolutely continuous.\vspace*{8pt}

\textit{Proof of} (\ref{eq:fluid.idling}).\quad
Since $x = (\mathbf{q},\mathbf{a}, \mathbf{y}, s)$ is absolutely continuous,
each component is too, which means that $y_{n}$ is
differentiable for almost all $t$. Pick some such
$t$, and suppose that $q_{n}(t)>0$. Consider some
small interval $I=[t,t+\delta]$ about $t$. Since
$q_{n}$ is continuous, we can choose $\delta$
sufficiently small that $\inf_{s\in I} q_{n}(s)>0$.
Since $\|\fbq^j-\mathbf{q}\|\tend0$, we can find $c>0$
such that $\inf_{s\in I} \fq^j_{n}(s) > c$ for all
$j$ sufficiently large. Since $\fx^j \in E_j$, there exists
a corresponding unscaled version of the system, say $X^j$,
and scaling parameter, say~$z_j$, so that $\fx^j(\cdot) = X^j(z_j
\cdot)/z_j$.
Therefore, it must be that the corresponding unscaled queue satisfies
$\inf_{s\in I} Q^j_{n}(z_j s)>z_j c$. That is, the\vadjust{\goodbreak}
queue size in the entire interval never vanishes
to $0$ and hence idling in the entire interval is not possible.
Therefore after rescaling we find $\fy^j_{n}(t+\delta/2)= \fy^j_{n}(t)$.
(The switch from $\delta$ to $\delta/2$ sidesteps any
discretization problems.) Therefore the same holds for
$y_{n}$ in the limit. We assumed $y_{n}$ to be differentiable
at $t$; the derivative must be $0$.\vspace*{8pt}

\textit{Proof of} (\ref{eq:fluid.mwm}).\quad
Pick a $t$ at which $s_\bpi$ is differentiable, and suppose that
$\bpi\bdot f(\mathbf{q}(t)) < \max_\brho\brho\bdot f(\mathbf
{q}(t))$. As above,
pick some small interval $I=[t,t+\delta]$ and $j$ sufficiently
large that
\[
\bpi\bdot f(\fbq^j(s))
<
\max_{\brho\in\sS} \brho\bdot f(\fbq^j(s))
\qquad\mbox{for $s\in I$}.
\]
Writing this in terms of the unscaled system and applying
Assumption \ref{cond.f},
\[
\bpi\bdot f(\bQ^j(z_j s))
<
\max_{\brho\in\sS} \brho\bdot f(\bQ^j(z_j s))\qquad
\mbox{for $s\in I$}.
\]
The MW-$f$ policy ensures by (\ref{eq:discrete.mwm}) that $\bpi$
will not be chosen throughout this entire interval,
so after rescaling we find $\fls^j_\bpi(t+\delta/2)-\fls^j_\bpi
(t)= 0$,
and taking the limit gives $s_\bpi(t+\delta/2)=s_\bpi(t)$.
Since $s_\bpi$ is assumed to be differentiable at $t$; the
derivative must be $0$.
\end{pf}

\subsubsection{\texorpdfstring{Proof of Theorem \protect\ref{thm:fluid} for
multi-hop networks.}{Proof of Theorem 4.3 for multi-hop networks}}

The proof of Theorem~\ref{thm:fluid} for single-hop network
applies verbatim, except that the two lemmas need to be replaced.
\begin{eqnarray*}
\mbox{Lemma \ref{lem:fluid.tight} (Tightness of fluid scaling)}
&\longrightarrow&
\mbox{Lemma \ref{lem:multihop.fluid.tight}.}
\\
\mbox{Lemma \ref{lem:limitpoints} (Dynamics at cluster points)}
&\longrightarrow&
\mbox{Lemma \ref{lem:multihop.limitpoints}.}
\end{eqnarray*}

\begin{lemma}[(Tightness of fluid scaling)]
\label{lem:multihop.fluid.tight}
Make the same assumptions as Theorem \ref{thm:fluid}, multi-hop case,
and use the same definition of $G_j$.
Then there exist a constant $A>0$
and a sequence $B_j\tend0$ such that for every\vadjust{\goodbreak}
$(\fx,z)\in G_j$, $|\fx(0)|\leq K$ and
\[
| \fx(u)-\fx(t)|
\leq A |u-t| + B_j\qquad
\mbox{for all $0\leq t,u\leq T$.}
\]
\end{lemma}
\begin{pf}
Consider $(\fx,z) \in G_j$, $\fx=(\fbq,\fba,\fby,\fls)$.
The bound $|\fx(0)| \leq K$ follows from an
argument similar to that in the single-hop case. The bounds on
the arrival process, the idleness and service allocation are as
in the single-hop case: for any $0\leq t < u \leq T$,
\begin{eqnarray*}
| \fba(u)-\fba(t) |
&\leq&(u-t) \lambda^{\max} + 2 \varepsilon_j,\\
| \fy_{n}(u)-\fy_{n}(t) |
&\leq&(u-t) S^{\max} + 2 S^{\max}/z_j^{\min},\\
| \fls_\bpi(t) - \fls_\bpi(s) |
&\leq&(u-t) + 2/z_j^{\min},
\end{eqnarray*}
where $z_j^{\min}=\inf\{z \dvtx (\fx,z)\in G_j\}$.
The bound on queue size is a little different. Note that
(\ref{eq:discrete.queue.multihop}) carries through to
the fluid-scaling, that is,
\[
\fbq(t) = \fbq(0) + \fba(t) - (I-R^\T) \sum_\bpi\fls_\bpi
(t)\bpi+
\fby(t),
\]
thus
\begin{eqnarray*}
| \fq_{n}(u)-\fq_{n}(t) |
&\leq&
| \fa_{n}(u)-\fa_{n}(t) |
+ \sum_\bpi|[(I-R^\T)\bpi]_n| | \fls_\bpi(u)
- \fls_\bpi(t) |\\[-2pt]
&&{} + | \fy_{n}(u)-\fy_{n}(t) |\\[-2pt]
&\leq&
(u-t)
\bigl( \lambda^{\max} + |\sS| (N S^{\max})S^{\max} + S^{\max}
\bigr) \\[-2pt]
&&{} +
\bigl(2|\sS| (N S^{\max})S^{\max} + 2 S^{\max} \bigr)/z_j^{\min}
+ 2\varepsilon_j.
\end{eqnarray*}
Putting all these together, for any $(\fx, z) \in G_j$,
\[
|\fx(u)-\fx(t)| \leq A (t-s) + B_j,
\]
where the constants $A$ and $B_j$ are
%
%
\begin{eqnarray}
\label{eq:multihop.tight.cont.bound}
A &=&
(1+N)\lambda^{\max} +2NS^{\max} + |\sS| +
|\sS|(NS^{\max})^2,\nonumber\\[-9pt]\\[-9pt]
B_j &=&
\bigl( 4N S^{\max} + 2|\sS| + 2|\sS|(N S^{\max})^2\bigr)
\frac{1}{z_j^{\min}}
+ (2N+2)\varepsilon_j.
\nonumber
\end{eqnarray}
Here $B_j\tend0$ as $j\tend\infty$ since
$\varepsilon_j \tend0$ by (\ref{cond.rateconv}) and $z_j^{\min}
\tend\infty$
as $j\tend\infty$ by~(\ref{eq:fluid.scale}).
\end{pf}

\begin{lemma}[(Dynamics at cluster points)]
\label{lem:multihop.limitpoints}
Under the setup of Theorem~\ref{thm:fluid} for a multihop network,
let $E_j = \{\fx\dvtx(\fx, z) \in G_j\}$. Then $x \in\FMSm_K$
if $x$ is a cluster point of the $E_j$ sequence.\vspace*{-2pt}
\end{lemma}

\begin{pf}
Given a cluster point $x=(\mathbf{q},\mathbf{a},\mathbf{y},s)$,
let there be $(\fx^j,z_j) \in G_j$ so that $\fx^j \to x$,
as in the proof of Lemma \ref{lem:limitpoints}.
Now the bound $|x(0)|\leq K$
and equations\vadjust{\goodbreak} (\ref{eq:fluid.arrivals})--(\ref{eq:fluid.idling})
all work exactly as in the single-hop case, as does the queue size
equation~(\ref{eq:multihop.fluid.queue}). The only equation that needs further
argument is the MW-$f$ backpressure equation
(\ref{eq:multihop.fluid.mwm}).\vspace*{8pt}

\textit{Proof of} (\ref{eq:multihop.fluid.mwm}).\quad
Pick a $t$ at which $s_\bpi$ is differentiable, and suppose that
$\bpi\bdot(I-R)f(\mathbf{q}(t)) < \max_\brho\brho\bdot
(I-R)f(\mathbf{q}(t))$. As in Lemma
\ref{lem:limitpoints}, proof of (\ref{eq:fluid.mwm}), it must be that
there is some small interval $I=[z_jt,z_jt+z_j\delta]$ such that
$\bpi$ is not chosen for any $\tau\in I$, therefore $\dot{s}_\bpi(t)=0$.
\end{pf}

\section{\texorpdfstring{Fluid model behavior (single-hop case).}{Fluid
model behavior (single-hop
case)}}
\label{sec:fluid.behavior}

In this section we prove certain properties of fluid model solutions,
which will be needed for the main result of this paper,
multiplicative state space
collapse. In order to state these properties, we first need some
definitions. We then state a portmanteau theorem listing all the
properties, and give an example to illustrate the definitions.
The rest of the section is given over to proofs and supplementary lemmas.\vadjust{\goodbreak}

This section deals with a single-hop switched network;
in the next section we give corresponding results for multi-hop.
Our reason for giving separate single-hop
and multi-hop proofs, rather than just treating single-hop as a~special
case of multi-hop, is that our multi-hop results place additional
restrictions on the set of allowed schedules (Assumption
\ref{cond.monotone}) beyond what is required for single-hop
networks. This mainly affects the proof; the portmanteau theorem for
multi-hop is nearly identical to that for single-hop.

\begin{definition}[(Admissible region)]
Let $\sS\subset\RealsP^N$ be the set of allowed schedules. Let
$\langle{\sS}\rangle$ be the convex hull of $\sS$,
\[
\langle{\sS}\rangle=
\biggl\{ \sum_{\bpi\in\sS} \alpha_\bpi\bpi\dvtx
\sum_{\bpi\in\sS} \alpha_\bpi=1\mbox{, and }
\alpha_\bpi\geq0 \mbox{ for all $\bpi$} \biggr\}.
\]
Define the \textit{admissible
region} $\Lambda$ to be
\[
\Lambda= \{\blambda\in\RealsP^N \dvtx
\blambda\leq\bsigma\mbox{ componentwise,
for some $\bsigma\in\langle{\sS}\rangle$} \}.
\]
\end{definition}
\begin{definition}[(Static planning problems and virtual resources)]
Define the optimization problem $\PRIMAL(\blambda)$ for
$\blambda\in\RealsP^N$ to be
\begin{eqnarray*}
&&\mbox{minimize }
\sum_{\bpi\in\sS} \alpha_\bpi
\quad\mbox{over}\quad
\alpha_\bpi\in\RealsP\qquad\mbox{for all $\bpi\in\sS$}\\
&&\mbox{such that }
\blambda\leq\sum_{\bpi\in\sS} \alpha_\bpi\bpi
\mbox{ componentwise}.
\end{eqnarray*}
Let $\DUAL(\blambda)$ be the dual to this: it is
\begin{eqnarray*}
&&\mbox{maximize }
\bxi\bdot\blambda
\quad\mbox{over}\quad
\bxi\in\RealsP^\N\qquad
\mbox{such that }
\max_{\bpi\in\sS} \bxi\bdot\bpi\leq1.
\end{eqnarray*}
Let $E$ be the set of extreme points of the feasible region of the
dual problem; the feasible region\vadjust{\goodbreak} is a finite convex polytope so $E$ is
finite. Define the set of \textit{virtual resources}
$\sS^*\subset\RealsP^N$ to be the set of maximal extreme points,
\[
\sS^* = \{\bxi\in E \dvtx \mbox{for all }\bzeta\in E,
\bxi\leq\bzeta\implies\bxi=\bzeta\}.
\]
Define the set of \textit{critically loaded virtual resources}
$\CLVR(\blambda)$ to be
\[
\Xi(\blambda) = \{\bxi\in\sS^* \dvtx \bxi\bdot\blambda= 1
\}.
\]
\end{definition}

Both problems are clearly feasible, and the optimum is attained in
each. By Slater's condition there is strong duality, that is,
$\PRIMAL(\blambda)=\DUAL(\blambda)$. [When we write
$\PRIMAL(\blambda)$ or $\DUAL(\blambda)$ in mathematical expressions,
we mean the optimum value, not the optimizer.]
Clearly, $\PRIMAL(\blambda)\leq1$ if and only if $\blambda$ is
feasible.

Laws \cite{lawsphd,lawsaap} and Kelly and Laws \cite
{kellylawsrespool} used
primal and dual problems of this general sort for describing
multi-hop queueing networks with routing choice.
Harrison \cite{harrisoncanonical}
extended the problems for stochastic processing networks.
\begin{definition}[(Lyapunov function and lifting map)]
\label{def:lyapunov}
Let the scheduling policy be MW-$f$, where $f$ satisfies Assumption
\ref{cond.f}.
Define the function $L\dvtx\RealsP^N\to\RealsP$ by
\[
L(\mathbf{q}) = F(\mathbf{q})\bdot\bOne,
\]
where $F(x)=\int_0^xf(y) \,dy$ for $x\in\Reals$, and
$F(\mathbf{q})=[F(q_n)]_{1\leq n\leq N}$ as per the notation in
Section \ref{sec:intro}.
Define the optimization problem $\ALGD(\mathbf{q})$ to be
\begin{eqnarray*}
&&\mbox{minimize }
L(\br)
\quad\mbox{over}\quad
\br\in\RealsP^N\\
&&\mbox{such that }
\bxi\bdot\br\geq\bxi\bdot\mathbf{q}
\qquad\mbox{for all } \bxi\in\CLVR(\blambda)
\quad\mbox{and}\\
&&\hphantom{\mbox{such that }}r_n\leq q_n\qquad\mbox{for all $n$ such that $\lambda_n=0$}.
\end{eqnarray*}
Note that $F$ is strictly convex and increasing, and the feasible
region is convex; hence this problem has a unique optimizer.
Define the \textit{lifting map} $\LIFTW\dvtx\RealsP^N\to\RealsP^N$ by
setting $\LIFTW(\mathbf{q})$ to be the optimizer.
\end{definition}

Note that $\ALGD$ and
$\LIFTW$ both depend on $\blambda$ and $f$, but we will surpress this
dependency when the context makes it clear which $\blambda$ and $f$
are meant.

The results in this section apply to any
$\blambda\in\Lambda$. However, if \mbox{$\PRIMAL(\blambda)<1$}, then
$\CLVR(\blambda)$ is empty, so $\LIFTW(\mathbf{q})=\bZero$ for all
$\mathbf{q}$. The
results are only interesting when $\PRIMAL(\blambda)=1$, so we define
\[
\dLambda= \{ \blambda\in\Lambda\dvtx\PRIMAL(\blambda)=1
\}.
\]

We can now state the main result of this section.
\begin{theorem}[(Portmanteau theorem, single-hop version)]
\label{thm:fp}
Let $\blambda\in\Lambda$.
Consider a single-hop switched network running MW-$f$, where $f$
satisfies Assumption~\ref{cond.f}.
\begin{longlist}
\item
For any $K<\infty$, $\{\mathbf{q}\in\RealsP^N\dvtx L(\mathbf{q})\leq K\}
$ is compact.
Also, for any fluid model solution with arrival rate
$\blambda$, $L(\mathbf{q}(t))\leq L(\mathbf{q}(0))$ for all $t\geq0$.
\item
$\LIFTW$ is continuous.
\item
If $\mathbf{q}=\LIFTW(\mathbf{q})$, then
$\LIFT\WORK(\kappa\mathbf{q})=\kappa\LIFTW(\mathbf{q})$ for
all $\kappa>0$.
\item
Say that $\mathbf{q}^0$ is an \textit{invariant state} if all fluid
model solutions
$\mathbf{q}(\cdot)$ with arrival rate $\blambda$,
starting at $\mathbf{q}(0)=\mathbf{q}^0$, satisfy $\mathbf
{q}(t)=\mathbf{q}^0$ for all $t\geq0$.
Then
$\mathbf{q}^0$ is an invariant $\mbox{state} \iff \mathbf{q}^0 = \LIFTW
(\mathbf{q}^0)$.
\item
For any $\varepsilon>0$ there exists some $H_\varepsilon<\infty$
such that, if $\mathbf{q}(\cdot)$ is a fluid model solution with arrival
rate $\blambda$, and $|\mathbf{q}(0)|\leq1$,
then $|\mathbf{q}(t)-\LIFTW(\mathbf{q}(t))|<\varepsilon$ for all
$t\geq H_\varepsilon$.
\end{longlist}
\end{theorem}

A loose interpretation of these results is that the MW-$f$ scheduling
policy seeks always to reduce $L(\mathbf{q})$ [part (i)],
but it is constrained from
reducing it too much, because it is not permitted to reduce the
workload at
any of the critically loaded virtual resource (the constraints of
$\ALGDUAL$).
However, it can choose how to allocate work between queues, subject to
those constraints.
It heads towards\vadjust{\goodbreak} a state where it is impossible to reduce
$L(\mathbf{q})$ any further [parts (iv) and
(v)].
In all the examples we have looked at, the fluid model solutions
reach an invariant state in finite time,
that is, (v) holds also for
$\varepsilon=0$, but we have not been able to prove this in general.

\subsection{\texorpdfstring{Example to illustrate $\Lambda$,
$\dLambda$, $\sS^*$ and $\CLVR$.}{Example to illustrate Lambda,
partial Lambda, S^* and Xi}}
Consider a system with \mbox{$N=2$} queues, $A$ and $B$. Suppose the set
$\sS$ of possible schedules consists of ``serve three packets from
queue $A$'' and ``serve one packet each from $A$ and~$B$.'' Write these two schedules as $\bpi^1=(3,0)$ and $\bpi^2=(1,1)$,
respectively.
Let~$\lambda_A$ and~$\lambda_B$ be the arrival
rates at the two queues, measured in packets per time slot.

\subsubsection*{\texorpdfstring{Determining $\Lambda$ and $\dLambda
$.}{Determining $\Lambda$ and
$\dLambda$}}
The arrival rate vector $\blambda=(\lambda_A,\lambda_B)$ is feasible
if there is some $\bsigma=(1-x)\bpi^1+x\bpi^2$ with $0\leq x\leq1$
such that $\blambda\leq\bsigma$. In words, the arrival rates are
feasible if the switch can divide its time between the two possible
schedules in such a way that the service rates at the two queues are
at least as big as the arrival rates.
Schedule $\bpi^2$ is the only schedule which serves queue $B$, so we
would need $x\geq\lambda_B$. If $\lambda_B>1$, then it is impossible
to serve all the work that arrives at queue $B$. Otherwise, we may as
well set
$x=\lambda_B$.
The total
amount of service given to queue $A$ is then $3(1-x)+x=3-2\lambda_B$;
if $\lambda_A\leq3-2\lambda_B$, then it is possible to serve all the
work arriving at queue $A$. We have concluded that
\[
\Lambda= \bigl\{ (\lambda_A,\lambda_B) \dvtx
\lambda_B\leq1 \mbox{ and } \tfrac{1}{3}\lambda_A+\tfrac
{2}{3}\lambda_B\leq1
\bigr\}.
\]
Further algebra tells us that
\[
\PRIMAL(\blambda) =
\max\bigl(\lambda_B, \tfrac{1}{3}\lambda_A+\tfrac{2}{3}\lambda_B
\bigr).
\]
Hence
\[
\dLambda= \bigl\{ (\lambda_A,\lambda_B) \in\Lambda\dvtx
\lambda_B=1 \mbox{ or }
\tfrac{1}{3}\lambda_A+\tfrac{2}{3}\lambda_B=1 \bigr\}.
\]

\subsubsection*{\texorpdfstring{Determining $\sS^*$ and $\CLVR
$.}{Determining $\sS^*$ and
$\CLVR$}}
The feasible region of $\DUAL(\blambda)$ is
\[
\{ (\xi_1,\xi_2)\in\RealsP^2 \dvtx 3\xi_1\leq1 \mbox{ and }
\xi_1+\xi_2\leq1 \}.
\]
The extreme points may be found by sketching the feasible region; they are
$(0,0)$, $(\third,0)$, $(\third,\twothirds)$ and $(0,1)$.
Clearly the maximal extreme points, that is, the virtual resources, are
\[
\sS^* = \bigl\{ \bigl(\tfrac{1}{3},\tfrac{2}{3}\bigr), (0,1) \bigr\}.
\]
The set of critically loaded virtual resources
depends on $\lambda_A$ and $\lambda_B\dvtx(0,1)\in\CLVR(\blambda)$ iff
$\lambda_B=1$, and $(\third,\twothirds)\in\CLVR(\blambda)$ iff
$\lambda_A/3+2\lambda_B/3=1$.

\subsubsection*{\texorpdfstring{Interpretation of virtual
resources.\footnote{cf. Laws \cite{lawsphd}, Example
4.4.3.}}{Interpretation of virtual
resources$^1$}}

Each virtual resource $\bxi\in\sS^*$ may be interpreted as a virtual
queue. For example, take $\bxi=(\third,\twothirds)$, and define
the virtual queue size to be $\bxi\bdot\bQ= Q_A/3+2Q_B/3$.
Think of the virtual\vadjust{\goodbreak} queue as consisting of tokens:
every time a packet arrives to
queue $A$ put $\third$ tokens into the virtual queue, and every time
a packet arrives to queue $B$ put in $\twothirds$ tokens. The schedule
$\bpi^1$ can remove at most $3\times\third=1$ token, and schedule~%
$\bpi^2$ can remove at most $\third+\twothirds=1$ token.
In order that the total rate at which tokens arrive should be no more
than the maximum rate at which we can remove tokens, we need
\[
\lambda_A/3+2\lambda_B/3 \leq1,
\]
that is, $\blambda\bdot\bxi\leq1$. If
$\DUAL(\blambda)=\PRIMAL(\blambda)>1$, then there is some
$\bxi\in\sS^*$ such that $\blambda\bdot\bxi>1$, which means that the
corresponding virtual queue is unstable; hence the original system is
unstable.

\subsection{\texorpdfstring{Proofs for the portmanteau theorem.}{Proofs
for the portmanteau
theorem}}

Throughout this subsection we consider a single-hop switched network
running MW-$f$ with arrival rates $\blambda\in\Lambda$.

The first claim of Theorem \ref{thm:fp}(i), that
$\{\mathbf{q}\in\RealsP^N\dvtx L(\mathbf{q})\leq K\}$ is compact for any
$K<\infty$,
follows straightforwardly from the facts that $L(\mathbf{q})\tend
\infty$ as
$|\mathbf{q}|\tend\infty$, and $L(\cdot)$ is continuous.
The second claim follows from a standard result (first given by
Dai and Prabhakar \cite{daibala}, for an input-queued switch), which
we include here
for the sake of completeness.
\begin{lemma}
\label{lem:lyapunov}
For all $\mathbf{q}\in\RealsP^N$,
%
%
\begin{equation}
\label{eq:lyapunov.neg}
\blambda\bdot f(\mathbf{q}) - \max_{\bpi\in\sS} \bpi\bdot
f(\mathbf{q}) \leq0.\vadjust{\goodbreak}
\end{equation}
Also, every fluid model solution satisfies
\[
\frac{d}{dt} L(\mathbf{q}(t)) = \blambda\bdot f(\mathbf{q}(t)) -
\max_{\bpi\in\sS} \bpi\bdot
f(\mathbf{q}(t))
\leq0.
\]
\end{lemma}
\begin{pf}
Since $\blambda\in\Lambda$,
we can write $\blambda\leq\bsigma$ componentwise for some
$\bsigma=\sum_\bpi\alpha_\bpi\bpi$
with $\alpha_\bpi\geq0$ and
$\sum\alpha_\bpi=1$. Hence
\begin{eqnarray*}
\blambda\bdot f(\mathbf{q}) - \max_\brho\brho\bdot f(\mathbf{q})
&=&
\sum_\bpi\alpha_\bpi\bpi\bdot f(\mathbf{q}) - \max_\brho\brho
\bdot
f(\mathbf{q})\\
&\leq&
\biggl(\sum_\bpi\alpha_\bpi-1\biggr)\max_\brho\brho\bdot
f(\mathbf{q})
\leq0.
\end{eqnarray*}
For the claim about fluid model solutions,
\begin{eqnarray*}
\frac{d}{dt} L(\mathbf{q}(t)) &=& \dot{\mathbf{q}}(t) \bdot
f(\mathbf{q}(t))\\
&=&
\biggl(\blambda-\sum_{\bpi\in\sS} \dot{s}_\bpi(t)\bpi+\dot
{\mathbf{y}}(t) \biggr)
\bdot f(\mathbf{q}(t))\qquad
\mbox{by differentiating (\ref{eq:fluid.queue})}\\
&=&
\biggl(\blambda-\sum_\bpi\dot{s}_\bpi(t)\bpi\biggr) \bdot
f(\mathbf{q}(t))
\qquad\mbox{by (\ref{eq:fluid.idling}), using $f(0)=0$}\\
&=&
\blambda\bdot f(\mathbf{q}(t)) - \max_\brho\brho\bdot f(\mathbf
{q}(t)) \sum_\bpi
\dot{s}_\bpi(t)
\qquad\mbox{by (\ref{eq:fluid.mwm})}\\
&=&
\blambda\bdot f(\mathbf{q}(t)) - \max_\brho\brho\bdot f(\mathbf{q}(t))
\qquad\mbox{by (\ref{eq:fluid.busy})}\\
&\leq&
0
\qquad\mbox{by (\ref{eq:lyapunov.neg})}.
\end{eqnarray*}
\upqed\end{pf}

To prove Theorem \ref{thm:fp}(ii), it is useful to
work with a ``fuller'' representation of the lifting map. Let $E$ be the
set of extreme feasible solutions of $\DUAL(\blambda)$, and define
%
%
\begin{equation}
\label{eq:clvrplus}
\CLVR^+(\blambda) = \{\bxi\in E \dvtx \bxi\bdot\blambda=1
\}.
\end{equation}
This includes nonmaximal extreme points, whereas $\CLVR(\blambda)$
only includes
maximal extreme points.
\begin{lemma}
\label{lem:liftwplus}
The lifting map $\LIFTW(\mathbf{q})$ is the unique solution to the
optimization problem $\ALGDUAL^+(\mathbf{q})$,
\begin{eqnarray*}
&&\mbox{minimize }
L(\br)
\quad\mbox{over}\quad
\br\in\RealsP^N\\
&&\mbox{such that }
\bxi\bdot\br\geq\bxi\bdot\mathbf{q}
\qquad\mbox{for all } \bxi\in\CLVR^+(\blambda).
\end{eqnarray*}
\end{lemma}
\begin{pf}
$\ALGDUAL^+(\mathbf{q})$ has a unique minimum for the same reason
that $\ALGDUAL(\mathbf{q})$ has a unique minimum.\vadjust{\goodbreak}

Next we claim that if $\br$ is feasible for $\ALGDUAL(\mathbf{q})$
then it
is feasible for $\ALGDUAL^+(\mathbf{q})$. Pick any
$\bxi\in\CLVR^+(\blambda)$. By definition,
$\bxi$ is an extreme feasible solution of $\DUAL(\blambda)$ and
$\bxi\bdot\blambda=1$. Since it is an extreme feasible solution,
$\bxi\leq\bzeta$ for some virtual resource $\bzeta\in\sS^*$. Since
$\bxi\bdot\blambda=1$ we know $\bzeta\bdot\blambda\geq1$, but by
assumption $\blambda\in\Lambda$; hence $\bzeta\bdot\blambda=1$
and furthermore
$\xi_n<\zeta_n$ only
for~$n$ where $\lambda_n=0$. Now,
\[
\bxi\bdot\br- \bxi\bdot\mathbf{q}
=
(\bzeta\bdot\br- \bzeta\bdot\mathbf{q})
+ (\bxi-\bzeta)\bdot(\br-\mathbf{q}).
\]
We assumed that $\br$ is feasible for $\ALGDUAL(\mathbf{q})$; by the first
constraint of $\ALGDUAL(\mathbf{q})$ the first term in the preceding equation
is positive; by the second constraint the second term is positive.
We have shown that $\bxi\bdot\br\geq\bxi\bdot\mathbf{q}$ for all
$\bxi\in\CLVR^+(\blambda)$; hence $\br$ is feasible for $\ALGDUAL
^+(\mathbf{q})$.

Next we claim that if $\br$ is optimal for $\ALGDUAL^+(\mathbf{q})$,
then it is
feasible for $\ALGDUAL(\mathbf{q})$.
Clearly it satisfies the first constraint of $\ALGDUAL(\mathbf{q})$. Suppose
it does not satisfy the second constraint, that is, that
$r_n>q_n$
for some $n$ where $\lambda_n=0$, and define $\br'$ by $r'_m=r_m$ if
$m\neq n$ and $r'_n=q_n$. Then $\br'<\br$, hence $L(\br')<L(\br)$.
Also, $\br'$ is feasible for $\ALGDUAL^+(\blambda)$. To see this,
pick any $\bzeta\in\CLVR^+(\blambda)$, and let $\bxi\in\CLVR
^+(\blambda)$ be such that
$\zeta_m=\xi_m$ if $m\neq n$ and $\xi_n=0$. Then
\[
\bzeta\bdot\br' = \bxi\bdot\br' + \zeta_n r'_n
=
\bxi\bdot\br+ \zeta_n r'_n
\geq
\bxi\bdot\mathbf{q}+ \zeta_n r'_n
=
\bzeta\bdot\mathbf{q}.\vadjust{\goodbreak}
\]
The inequality is because $\br$ is feasible for
$\ALGDUAL^+(\mathbf{q})$. This contradicts optimality of $\br$.

Putting these two claims together completes the proof.
\end{pf}

With this representation,
the lifting map $\LIFTW$ can be split into two parts. Let
$\CLVR^+(\blambda)=\{\bxi^1,\ldots,\bxi^V\}$ and define the \textit
{workload map}
$W\dvtx\RealsP^N\to\RealsP^V$ by $W(\mathbf{q}) = [\bxi^v\bdot\mathbf
{q}]_{1\leq v\leq
V}$. Also define $\LIFT\dvtx\RealsP^V\to\RealsP^N$ by
%
%
\begin{equation}
\label{eq:liftwplus}
\LIFT(w) = \argmin\{ L(\br) \dvtx\br\in\RealsP^N
\mbox{ and }
\bxi^v\bdot\br\geq w_v
\mbox{ for }
1\leq v\leq V
\}.
\end{equation}
(This has a unique optimum for the same reason that $\ALGDUAL^+$ and
$\ALGDUAL$ have.)
Then the lifting map is simply the composition of $\LIFT$ and $W$. It
is clear that $W$ is continuous; to prove Theorem
\ref{thm:fp}(ii) we just need to prove that $\LIFT$ is
continuous.
\begin{lemma}
\label{lem:lift.continuous}
$\LIFT$ is continuous.
\end{lemma}
\begin{pf}
If $\CLVR^+(\blambda)$ is empty, then $\LIFT$ is trivial and the
result is trivial. In what follows, we shall assume that
$\CLVR^+(\blambda)$ is nonempty, and we will abbreviate it to $\CLVR^+$.
Furthermore note that for every $\bxi\in\CLVR^+$ there is some queue~$n$
such that $\xi_n>0$; this is because $\bxi\bdot\blambda=1$ by
definition of $\CLVR^+$.

Pick any sequence $w^k\tend w\in\RealsP^V$, and let $\br^k=\LIFT(w^k)$
and $\br=\LIFT(w)$. We want to prove that $\br^k\tend\br$. We shall
first prove that there is a compact set $[0,h]^N$ such that
$\br^k\in[0,h]^N$ for all $k$. We shall then prove that any convergent
subsequence of $\br^k$ converges to $\br$; this establishes
continuity of~$\LIFT$.\vadjust{\goodbreak}

First---compactness. A suitable value for $h$ is
\[
h = \max_{1\leq v\leq V} \max_{n\dvtx\xi_n>0} \sup_k \frac{w^k_v}{\xi^v_n}.
\]
Note than the maximums are over a nonempty set, as noted at the
beginning of the proof. Note also that $h$ is finite because $w$ is
finite. Now, suppose that $\br^k\notin[0,h]^N$ for some $k$,
that is, that there is some queue $n$ for which $r^k_n>h$, and let
$\br'=\br^k$ in each component except for $r'_n=h$. We claim that
$\br'$ satisfies the constraints
of the optimization problem for $\LIFT(w^k)$.
To see this, pick
any $\bxi^v\in\CLVR^+$; either $\xi^v_n=0$ in which case
$\bxi^v\bdot\br'=\bxi^v\bdot\br^k\geq w^k_v$, or $\xi^v_n>0$ in which
case $\bxi^v\bdot\br' \geq\xi^v_n h \geq w^k_v$ by construction of
$h$. Applying this repeatedly, if $\br^k\notin[0,h]^N$,
then we can reduce it to a queue size
vector in $[0,h]^N$, thereby improving on $L(\br^k)$, yet still
meeting the constraints of the optimization problem for $\LIFT(w^k)$;
this contradicts the
optimality of $\br^k$. Hence $\br^k\in[0,h]^N$.

Next---convergence on subsequences. With a slight abuse of notation,
let $\LIFT(w^k)=\br^k\tend\bs$ be a convergent subsequence, and recall
that $\LIFT(w)=\br$ and $w^k\tend w$. By continuity of the constraints,
$\bs$ is feasible for the optimization problem for $\LIFT(w)$;
we shall next show that
$L(\bs)\leq L(\br)$. Since $\br$ is the unique optimum, it must be
that $\bs=\br$.\vadjust{\goodbreak}

It remains to show that $L(\bs)\leq L(\br)$. Consider the sequence
$\br+\varepsilon^k\bOne$ as candidate solutions to the problem
$\LIFT(w^k)$
where
\[
\varepsilon^k = \max_{1\leq v\leq V}
\frac{w^k_v-w_v}{\bxi^v\bdot\bOne}.
\]
This choice ensures that the candidates are feasible, since
\[
\bxi^v\bdot(\br+\varepsilon^k\bOne)
=
\bxi^v\bdot\br+ \varepsilon^k \bxi^v\bdot\bOne
\geq
\bxi^v\bdot\br+ w^k_v-w_v
\geq w^k_v.
\]
(If we had used $\bxi\in\CLVR$ rather than $\bxi\in\CLVR^+$,
it would not necessarily be true that the candidates are feasible;
this is why we introduced Lemma \ref{lem:liftwplus}.)
Since the candidates are feasible solutions to the problem $\LIFT
(w^k)$, and
$\br^k$ is an optimal solution, it must be that
\[
L(\br^k) \leq L(\br+\varepsilon^k\bOne).
\]
Taking the limit as $k\tend\infty$, and noting that $L$ is continuous
and $\varepsilon^k\tend0$, we find
\[
L(\bs) \leq L(\br)
\]
as required. This completes the proof.
\end{pf}

For the proof of Theorem \ref{thm:fp}(iii), it is useful
to work with a different representation of the constraint of
$\ALGDUAL$, provided by the following lemma.

\begin{lemma}
\label{lem:equiv}
\begin{longlist}
\item
$\LIFTW(\mathbf{q})=[\mathbf{q}+t(\blambda-\bsigma)]^+$
for some $t\geq0$ and $\bsigma\in\langle{\sS}\rangle$.
\item
$[\mathbf{q}+t(\blambda-\bsigma)]^+$ is feasible for $\ALGDUAL
(\mathbf{q})$ for all
$t\geq0$ and $\bsigma\in\langle{\sS}\rangle$.\vadjust{\goodbreak}
\end{longlist}
\end{lemma}
\begin{pf}
{(i)}
We will shortly prove that the following are equivalent, for all
$\mathbf{q}$
and $\br\in\RealsP^N$:
%
%
\begin{eqnarray}
\label{eq:equiv1}
\br&\geq&\mathbf{q}+t(\blambda-\bsigma) \qquad\mbox{for some }
t\geq0, \bsigma\in\langle{\sS}\rangle;\\
\label{eq:equiv2}
\bxi\bdot\br&\geq&\bxi\bdot\mathbf{q} \qquad\mbox{for all }
\bxi\in\CLVR^+(\blambda).
\end{eqnarray}
We use this equivalence as follows. From Lemma \ref{lem:liftwplus} we
know that $\LIFTW(\mathbf{q})$ is the solution of $\ALGDUAL
^+(\mathbf{q})$. That is,
letting $\mathbf{q}'=\LIFTW(\mathbf{q})$, equation (\ref
{eq:equiv2}) holds with
$\mathbf{q}'$ in the place of $\br$. Hence (\ref{eq:equiv1}) holds
for some
$t\geq0$ and $\bsigma\in\langle{\sS}\rangle$; moreover since
$\mathbf
{q}'\geq\bZero$ it must
be that
\[
\mathbf{q}' \geq[\mathbf{q}+ t (\blambda-\bsigma)]^+.
\]
We claim that this inequality is in fact an equality. To see
this, note that $\br=[\mathbf{q}+t(\blambda-\bsigma)]^+$ satisfies
(\ref{eq:equiv1}); hence it satisfies (\ref{eq:equiv2}); hence it is a
feasible solution of $\ALGDUAL^+(\mathbf{q})$. Note also
that $L(\cdot)$ is increasing componentwise, hence $L(\mathbf
{q}')\geq
L(\br)$. But $\ALGDUAL^+(\mathbf{q})$ has a unique minimum, hence
$\mathbf{q}'=\br$
as required. This completes the proof of Lemma
\ref{lem:equiv}(i), once we have proved the equivalence
between (\ref{eq:equiv1}) and~(\ref{eq:equiv2}).\vadjust{\goodbreak}

\textit{Proof that} (\ref{eq:equiv1})${}\implies{}$(\ref{eq:equiv2}).\quad
Pick any $\bxi\in\CLVR^+(\blambda)$. By definition of
$\CLVR^+(\blambda)$, we know: $\bxi\geq\bZero$; $\bxi\bdot\bpi
\leq1$ for all
$\bpi\in\sS$, hence $\bxi\bdot\bsigma\leq1$ for all
$\bsigma\in\langle{\sS}\rangle$; and $\bxi\bdot\blambda=1$. Hence
\begin{eqnarray*}
\bxi\bdot\br
&\geq&
\bxi\bdot\mathbf{q}+ t (\bxi\bdot\blambda- \bxi\bdot
\bsigma)
\qquad\mbox{assuming $\mathbf{q}$, $\br\in\RealsP^N$ satisfying
(\ref{eq:equiv1})}\\
&=&
\bxi\bdot\mathbf{q}+ t (1 - \bxi\bdot\bsigma)
\geq
\bxi\bdot\mathbf{q}+ t(1-1)
=
\bxi\bdot\mathbf{q}.
\end{eqnarray*}

\textit{Proof that}
(\ref{eq:equiv1})${}\Longleftarrow{}$(\ref{eq:equiv2}).\quad
Let $\mathbf{q}$ and $\br$ satisfy (\ref{eq:equiv2}), and let
$\bsigma' = \blambda- (\br-\mathbf{q})/t$
for some sufficiently large $t\in\RealsP$. We shortly show that the
value of $\DUAL(\bsigma')$ at its optimum is $\leq1$. By strong
duality the
value of $\PRIMAL(\bsigma')$ at its optimum is likewise $\leq1$, and so
by definition of $\PRIMAL(\bsigma')$ we can find some $\bsigma\in
\langle{\sS}\rangle$ such
that $\bsigma'\leq\bsigma$ componentwise. Then
\[
\br= \mathbf{q}+ t(\blambda-\bsigma')
\geq\mathbf{q}+ t(\blambda-\bsigma),
\]
that is, $\br$ satisfies (\ref{eq:equiv1}).

It remains to show that the value of $\DUAL(\bsigma')$ at its optimum
is $\leq1$,
that is, that $\bzeta\bdot\bsigma'\leq1$ for all dual-feasible
$\bzeta$.
We have assumed that $\blambda\in\Lambda$, hence
$\bzeta\bdot\blambda\leq1$.
On one hand, if $\bzeta\bdot\blambda=1$, then
it follows from the definition of $\CLVR^+$ that
$\bzeta\in\langle{\CLVR^+(\blambda)}\rangle$,
hence
\begin{eqnarray*}
\bzeta\bdot\bsigma'
&=&
\bzeta\bdot\blambda- \bzeta\bdot(\br-\mathbf{q})/t\\
&=&
1 - \bzeta\bdot(\br-\mathbf{q})/t \qquad\mbox{since $\bzeta\bdot
\blambda=1$}\\
&\leq&1 \qquad\mbox{by (\ref{eq:equiv2})}.
\end{eqnarray*}
On the other hand, if $\bzeta\bdot\blambda<1$, then
\[
\bzeta\bdot\bsigma'
<
1-\bzeta\bdot(\br-\mathbf{q})/t,\vadjust{\goodbreak}
\]
and this is $<1$ for $t$ sufficiently large.
Either way,
$\bzeta\cdot\bsigma'\leq1$.
Therefore the value of $\DUAL(\bsigma')$ at its optimum is $\leq1$.

(ii)
For this, we need to check two feasibility conditions of
$\ALGDUAL(\mathbf{q})$.
The first feasibility condition is
\[
\bxi\bdot[\mathbf{q}+t(\blambda-\bsigma)]^+ \geq\bxi\bdot
\mathbf{q}\qquad
\mbox{for all $\bxi\in\CLVR(\blambda)$.}
\]
Pick any $\bxi\in\CLVR(\blambda)$. By definition of $\CLVR
(\blambda)$, $\bxi\geq\bZero$,
$\bxi\bdot\bpi\leq1$ for all $\bpi\in\sS$ hence
$\bxi\bdot\bsigma\leq1$ for all $\bsigma\in\langle{\sS}\rangle
$, and
$\bxi\bdot\blambda=1$,
thus
\[
\bxi\bdot[\mathbf{q}+t(\blambda-\bsigma)]^+
\geq
\bxi\bdot[\mathbf{q}+t(\blambda-\bsigma)]
=
\bxi\bdot\mathbf{q}+ t (1 - \bxi\bdot\bsigma)
\geq
\bxi\bdot\mathbf{q}
\]
as required. The second feasibility condition is that if $\lambda_n=0$
for some $n$, then
\[
[q_n + t (\lambda_n - \sigma_n)]^+
=
[q_n - t \sigma_n]^+
\leq
q_n.
\]
This is true because $\bsigma\geq\bZero$ componentwise for all
$\bsigma\in\langle{\sS}\rangle$.
\end{pf}

Theorem \ref{thm:fp}(iii)
is a corollary of the following lemma.
\begin{lemma}[(Scale-invariance of $\LIFTW$)]
\label{lem:lift.scaleinvariant}
Let $\mathbf{q}\in\RealsP^\N$.
Then $\LIFT\WORK(\kappa\mathbf{q})=\kappa\times\LIFT\WORK(\mathbf
{q})$ for all $\kappa>0$.\vadjust{\goodbreak}
\end{lemma}
\begin{pf}
We will first establish three
preliminary properties of $\LIFTW$. Preliminary 1 is used to prove 2,
and 2 and 3 are used in the main proof.\vspace*{7pt}

\textit{Preliminary} 1.\quad
If
$\mathbf{q}=\LIFTW(\mathbf{q}')$ for some
$\mathbf{q}'\in\RealsP^N$, then
%
%
\begin{equation}
\label{eq:opt.cond}
\blambda\bdot f(\mathbf{q}) = \max_{\bpi\in\sS} \bpi\bdot
f(\mathbf{q}).
\end{equation}
To see this,
suppose $\bpi\in\sS$ has maximal weight and consider
$\br=[\mathbf{q}+t(\blambda-\bpi)]^{+}$.
This is feasible for $\ALGDUAL(\mathbf{q})$ by Lemma \ref{lem:equiv}.
Now, using the fact that $f(0)=0$,
\[
\frac{d}{dt}
L\bigl([\mathbf{q}+t(\blambda-\bpi)]^+\bigr)\big|_{ t=0}
= (\blambda-\bpi)\bdot f(\mathbf{q}).
\]
Since $\mathbf{q}$ is optimal for $\ALGDUAL(\mathbf{q}')$ it is
optimal
for $\ALGDUAL(\mathbf{q})$, hence
$\blambda\bdot f(\mathbf{q})\geq\bpi\bdot f(\mathbf{q})$.
On the other hand, $\blambda\in\Lambda$ so $\blambda\leq\bsigma$
for some
$\bsigma\in\langle{\sS}\rangle$, hence $\blambda\bdot f(\mathbf
{q})\leq
\bsigma\bdot
f(\mathbf{q})\leq\bpi\bdot f(\mathbf{q})$.
Hence the result follows.\vspace*{7pt}

\textit{Preliminary} 2.\quad Suppose
that $\br=\LIFT\WORK(\mathbf
{q})$. From
Lemma \ref{lem:equiv},
$\br=[\mathbf{q}+t(\blambda-\bsigma)]^+$ for some
$t\geq0$ and
$\bsigma\in\langle{\sS}\rangle$. Then either $t=0$ or
%
%
\begin{equation}
\label{eq:optservice}
\bsigma\bdot f(\br) = \max_{\bpi\in\sS} \bpi\bdot f(\br).
\end{equation}
This is because $t$ is an optimal choice, so either $t$ is constrained
to be $0$ or
\[
\frac{d}{du} L\bigl([\mathbf{q}+u(\blambda-\bsigma)
]^+\bigr)\big|_{u=t}
=
(\blambda-\bsigma)\bdot f(\br)
=
0.
\]
In this second case, $\blambda\bdot f(\br)=\max_\bpi\bpi\bdot
f(\br)$ by (\ref{eq:opt.cond})
so the same is true for~$\bsigma$.\vadjust{\goodbreak}

\textit{Preliminary} 3.\quad
Suppose that $\br=\LIFT\WORK(\mathbf{q})$. From Lemma \ref
{lem:equiv}, we can write it as
$\br=[\mathbf{q}+t(\blambda-\bsigma)]^+$ for some
$\bsigma\in\langle{\sS}\rangle$. In fact, for any $T\geq t$ we can
write it as
%
%
\begin{equation}
\label{eq:extend}
\br= [\mathbf{q}+T(\blambda-\brho)]^+
\qquad\mbox{for some $\brho\in\langle{\sS}\rangle$}.
\end{equation}
To see this, recall that $\PRIMAL(\blambda)\leq1$, so we can
pick some $\bar{\blambda}\in\langle{\sS}\rangle$ such that
$\blambda\leq\bar{\blambda}$, whence
\begin{eqnarray*}
\br
&\geq&
[\mathbf{q}+t(\blambda-\bsigma)+(T-t)(\blambda-\bar{\blambda
})]^+\\[-2pt]
&=&
[\mathbf{q}+ T(\blambda-\brho)]^+
\qquad\mbox{where }\brho=\frac{t}{T}\bsigma+\frac{T-t}{T}\bar
{\blambda}\in\langle{\sS}\rangle.
\end{eqnarray*}
This last expression is feasible for $\ALGDUAL(\mathbf{q})$
by Lemma \ref{lem:equiv}.
Since
$\br$ is optimal for $\ALGDUAL(\mathbf{q})$, and the objective
function is
increasing pointwise, $\br=[\mathbf{q}+T(\blambda-\brho)]^+$
as claimed.\vspace*{8pt}

\textit{Main proof.}\quad
Let $\br=\LIFTW(\mathbf{q})$ and $\kappa\br'=\LIFTW(\kappa
\mathbf{q})$.
We know that $\kappa\br$ is feasible for $\ALGDUAL(\kappa\mathbf{q})$
because the constraints are linear;
we will
now show that $L(\kappa\br)\leq L(\kappa\br')$; hence $\kappa\br$ is
also optimal for $\ALGDUAL(\kappa\mathbf{q})$. By uniqueness of the optimum,
$\kappa\br=\kappa\br'$ as required.

It remains to prove that $L(\kappa\br)\leq L(\kappa\br')$.
Since $\br$ solves $\ALGDUAL(\mathbf{q})$ and $\kappa\br'$ solves
$\ALGDUAL(\kappa\mathbf{q})$, we can use Lemma \ref{lem:equiv} to write
\[
\br= [\mathbf{q}+ t(\blambda-\bsigma)]^+
,\qquad
\kappa\br' = [\kappa\mathbf{q}+ \kappa t'(\blambda-\bsigma
')]^+
\]
for $t,t'\!\in\!\RealsP$ and $\bsigma,\bsigma'\!\in\!\langle{\sS}\rangle$.
Indeed, for $T\!>\!\max(t,t')$ we can use (\ref{eq:extend}) to~write
\begin{eqnarray}
\br&=& \mathbf{q}+ T(\blambda-\brho+\mathbf{y})
\qquad\mbox{for $\brho\in\langle{\sS}\rangle$, $\mathbf{y}\in
\RealsP^\N
$}\nonumber\\[-2pt]
&&\eqntext{\mbox{where $y_{n}=0$ if
$r_{n}>0$},}\\[-2pt]
\br' &=& \mathbf{q}+ T(\blambda-\brho'+\mathbf{y}')
\qquad\mbox{for $\brho'\in\langle{\sS}\rangle$, $\mathbf{y}'\in
\RealsP
^\N$}\nonumber\\[-2pt]
&&\eqntext{\mbox{where $y'_{n}=0$ if
$r'_{n}>0$}.}
\end{eqnarray}
Now consider the value of $L(\cdot)$ along the trajectory from $\kappa
\br$ to
$\kappa\br'$.
Along this trajectory,
\begin{eqnarray*}
\frac{d}{du} L\bigl(\kappa\br+ (\br'-\br)u/T\bigr)\bigg|_{ u=0}
&=&
(\br'-\br) \bdot f(\kappa\br)/T\\[2pt]
&=&
(\brho-\brho' -\mathbf{y}+\mathbf{y}')\bdot f(\kappa
\br)\\[2pt]
&\geq&
(\brho-\brho' -\mathbf{y})\bdot f(\kappa\br)
\qquad\mbox{since $\mathbf{y}'\geq\bZero$}\\[2pt]
&=&
(\brho-\brho')\bdot f(\kappa\br)
\qquad\mbox{since $y_{n}=0$ if $r_{n}>0$}\\[2pt]
&\geq&
\brho\bdot f(\kappa\br) - \max_{\bpi\in\sS}\bpi\bdot f(\kappa
\br)
\qquad\mbox{for any $\brho'\in\langle{\sS}\rangle$}\\[2pt]
&=&0.
\end{eqnarray*}
The final equality is because $\brho\bdot f(\br)=\max_\bpi\bpi
\bdot
f(\br)$ by (\ref{eq:optservice}), so
$\brho\bdot f(\kappa\br) = \max_\bpi\bpi\bdot f(\kappa\br)$ by
Assumption \ref{cond.f}.
Since $L(\cdot)$ is convex, it follows that $L(\kappa\br')\geq
L(\kappa\br)$.
This completes the proof.\vadjust{\goodbreak}
\end{pf}

The proof of Theorem \ref{thm:fp}(iv) relies on the
following lemma.
\begin{lemma}[(Fluid model trajectories preserve $\ALGDUAL$ feasibility)]
\label{lem:feasible}
Consider any fluid model solution, for any scheduling policy, with
initial queue size
$\mathbf{q}(0)$. Then $\mathbf{q}(t)$ is feasible for
$\ALGDUAL(\mathbf{q}(0))$ for all $t\geq0$.
\end{lemma}
\begin{pf}
Pick any critically loaded virtual resource $\bxi\in\CLVR(\blambda)$.
By (\ref{eq:fluid.queue}),
\begin{eqnarray*}
\bxi\bdot\mathbf{q}(t)
&=&
\bxi\bdot\mathbf{q}(0) + t\bigl(\bxi\bdot\blambda-\bxi\bdot
\bsigma(t)\bigr)+\bxi\bdot\mathbf{y}(t)
\qquad\mbox{where }\bsigma(t)=\sum\bpi s_\bpi(t)/t\\
&\geq&
\bxi\bdot\mathbf{q}(0) + t\bigl(\bxi\bdot\blambda-\bxi\bdot
\bsigma(t)\bigr)
\qquad\mbox{since }\mathbf{y}(t)\geq\bZero\\
&\geq&
\bxi\bdot\mathbf{q}(0)+t(1-1)
= \bxi\bdot\mathbf{q}(0).
\end{eqnarray*}
The last inequality is because $\bxi\in\CLVR(\blambda)$; so
$\bxi\bdot\blambda=1$,
and $\bxi\bdot\bpi\leq1$ for all $\bpi\in\sS$ hence
$\bxi\bdot\bsigma\leq1$ for all $\bsigma\in\langle{\sS}\rangle$.
Finally, $\mathbf{q}(t)\leq\mathbf{q}(0)+t\blambda$ by (\ref
{eq:fluid.queue}) and~(\ref{ineq:fluid.littleidle}),
and this yields the second constraint of
$\ALGDUAL(\blambda)$ for queues $n$ with 0 arrival rate.\vadjust{\goodbreak}
\end{pf}

Theorem \ref{thm:fp}(iv) is implied by
parts (i) and (ii) of the following lemma.\vspace*{-2pt}

\begin{lemma}[(Characterization of invariant states of MW-$f$)]
\label{lem:fixedpoint}
The following are equivalent, for $\mathbf{q}^0\in\RealsP^N$:
\begin{longlist}
\item
$\mathbf{q}^0 = \LIFTW(\mathbf{q}^0)$;\vspace*{1pt}
\item
$\mathbf{q}^0$ is an invariant state;
\item
there exists a fluid model solution with $\mathbf{q}(t)=\mathbf{q}^0$
for all $t$;
\item
$\blambda\bdot f(\mathbf{q}^0) = \max_{\bpi\in\sS} \bpi\bdot
f(\mathbf{q}^0)$.\vspace*{-2pt}
\end{longlist}
\end{lemma}
\begin{pf}

\textit{Proof that} (i)${}\implies{}$(ii).\quad
Suppose that $\mathbf{q}^0=\LIFTW(\mathbf{q}^0)$, that is, that
$\mathbf{q}^0$ is optimal for
$\ALGDUAL(\mathbf{q}^0)$, and
consider any fluid model
solution which starts with $\mathbf{q}(0)=\mathbf{q}^0$.
On one hand, Lemma \ref{lem:lyapunov} says that
$L(\mathbf{q}(t))\leq L(\mathbf{q}^0)$.
On the other hand, Lemma~\ref{lem:feasible} says that $\mathbf{q}(t)$ is
feasible for $\ALGDUAL(\mathbf{q}^0)$.
Since $\ALGDUAL(\mathbf{q}^0)$ has a unique solution, it
must be that $\mathbf{q}(t)=\mathbf{q}^0$.\vspace*{8pt}

\textit{Proof that} (ii)${}\implies{}$(iii).\quad
It is easy to find a fluid model solution which starts at
$\mathbf{q}(0)=\mathbf{q}^0$: a limit point of the stochastic model
from Theorem
\ref{thm:fluid} will do. By (ii), the queue size vector
is constant.\vspace*{8pt}

\textit{Proof that} (iii)${}\implies{}$(iv).\quad
Suppose
there is a fluid model
solution with $\mathbf{q}(t)\,{=}\,\mathbf{q}^0$. Since $\mathbf{q}(\cdot
)$ is constant,
$\dot{L}(\mathbf{q}(t))\,{=}\,0$. Lemma \ref{lem:lyapunov} says that
\mbox{$\dot{L}(\mathbf{q}(t))\,{\leq}\,0$}, so the inequality in the proof
must be tight for all $t$, that is,
%
%
\begin{equation}
\label{eq:fixed.opt}
\blambda\bdot f(\mathbf{q}^0) = \max_{\bpi\in\sS} \bpi\bdot
f(\mathbf{q}^0).
\end{equation}

\textit{Proof that} (iv)${}\implies{}$(i).\quad
If $\mathbf{q}^0=\bZero$ then\vspace*{1pt} the result is trivial.
Otherwise, let $\br=\LIFTW(\mathbf{q}^0)$. By Lemma \ref{lem:equiv},
$\br=[\mathbf{q}^0+t(\blambda-\bsigma)]^+$ for some\vadjust{\goodbreak} $t\geq0$ and
$\bsigma\in\langle{\sS}\rangle$.
Consider the value of $L(\cdot)$ along the trajectory from $\mathbf
{q}^0$ to
$\br$
\begin{eqnarray*}
\frac{d}{du} L\bigl([\mathbf{q}^0+(\blambda-\bsigma
)u]^+\bigr) \bigg|_{ u=0}
&=&
(\blambda-\bsigma)\bdot f(\mathbf{q}^0)
\qquad\mbox{relying on }f(0)=0\\
&=&
\Bigl(\max_{\bpi\in\sS}\bpi\bdot f(\mathbf{q}^0)\Bigr)
- \bsigma\bdot f(\mathbf{q}^0)
\qquad\mbox{by part (iv)}\\
&\geq&0\qquad
\mbox{because }\bsigma\in\langle{\sS}\rangle.
\end{eqnarray*}
By convexity of $L$, $L(\br)\geq L(\mathbf{q}^0)$,
and $\mathbf{q}^0$ is obviously feasible for $\ALGDUAL(\mathbf{q}^0)$,
but we chose $\br$ to be optimal for $\ALGDUAL(\mathbf{q}^0)$,
and the optimum is unique. Therefore
$\mathbf{q}^0=\LIFTW(\mathbf{q}^0)$.\vspace*{-2pt}
\end{pf}

Theorem \ref{thm:fp}(v) is given by the
following lemma. Recall that we are using the norm $|\bx|={\max_n}
|x_n|$.\vspace*{-2pt}

\begin{lemma}
\label{lem:criticalfluid.time}
Given $\blambda\in\Lambda$,
for any $\varepsilon>0$ there exists an $H_\varepsilon$ such that
for every fluid model solution with arrival rate $\blambda$,
for which $|\mathbf{q}(0)|\leq1$, $|\mathbf{q}(t)-\LIFTW(\mathbf
{q}(t))|<\varepsilon$
for all $t\geq H_\varepsilon$.\vadjust{\goodbreak}
\end{lemma}
\begin{pf} The proof is inspired by Kelly and Williams \cite{KW04},
Theorem 5.2, Lemma 6.3.
We start with some definitions. Let
\begin{eqnarray*}
\cD&=& \{ \mathbf{q}\in\RealsP^\N\dvtx L(\mathbf{q})\leq
L(\bOne) \}\qquad
\mbox{for $L(\cdot)$ as in Definition \ref{def:lyapunov}};\\[-2pt]
\INV&=& \{ \mathbf{q}\in\cD\dvtx\LIFT\WORK(\mathbf
{q})=\mathbf{q}\};\\[-2pt]
\INV_\delta&=& \{ \mathbf{q}\in\cD\dvtx|\mathbf{q}-\br
|<\delta\mbox{ for some }
\br\in\INV\};\\[-2pt]
\cJ_\varepsilon&=& \{ \mathbf{q}\in\RealsP\dvtx |
\mathbf{q}-\LIFT\WORK(\mathbf{q}) | < \varepsilon\};\\[-2pt]
\cK_\delta&=& \Bigl\{ \mathbf{q}\in\cD\dvtx L(\mathbf{q})-L
(\LIFT\WORK(\mathbf{q}))
< \inf_{\br\in\cD\setminus\INV_\delta} L(\br)-L(\LIFT
\WORK(\br)) \Bigr\}.\vspace*{-2pt}
\end{eqnarray*}
We will argue that the function
$K(\mathbf{q})=L(\mathbf{q})-L(\LIFT\WORK(\mathbf{q}))$
is decreasing along fluid model trajectories, so once
you hit $\cK_\delta$ you stay there.
We will then argue that
$\INV\subset\cK_\delta\subset\INV_\delta\subset\cJ_\varepsilon
$ for sufficiently
small $\delta$.
Finally, we will bound the time it takes
to hit $\cK_\delta$.\vspace*{7pt}

\textit{$K$ is decreasing.}\quad
Lemma \ref{lem:lyapunov} says that for any fluid model solution,
$L(\mathbf{q}(\cdot))$ is decreasing.
From Lemma \ref{lem:feasible}, the feasible set for
$\ALGDUAL(\mathbf{q}(u))$ is a subset of the feasible set for
$\ALGDUAL(\mathbf{q}(t))$ for any $u\geq t\geq0$, hence
$\LIFTW(\mathbf{q}(u))\geq\LIFTW(\mathbf{q}(t))$, that is, $\LIFTW
(\mathbf{q}(\cdot))$ is increasing.
Therefore $K$ is decreasing (not necessarily strictly).

\textit{$\INV\subset\cK_\delta\subset\INV_\delta
\subset\cJ
_\varepsilon$.}\quad
To show $\INV\subset\cK_\delta$:
the map $\LIFTW$ is continuous by Theorem~\ref{thm:fp}(ii),
and $L(\cdot)$ is clearly
continuous, so $K(\cdot)$ is continuous;
also the set $\cD$ is compact by Theorem
\ref{thm:fp}(i),
and $\INV_\delta$ is open,
so $\cD\setminus\INV_\delta$ is compact; so
the infimum in the definition of $\cK_\delta$ is attained at some
$\hbr\in\cD\setminus\INV_\delta$. Now, $K(\mathbf{q})>0$ for
$\mathbf{q}\in\cD\setminus\INV$,
so $K(\hbr)>0$. Yet $K(\mathbf{q})=0$ for $\mathbf{q}\in\INV$.
Thus $\INV\subset\cK_\delta$.

It is clear by construction that $\cK_\delta\subset\INV_\delta$.\vadjust{\goodbreak}

To show $\INV_\delta\subset\cJ_\varepsilon$:
the map $\LIFT\WORK(\cdot)$ is continuous,
hence it is uniformly~con\-tinuous on the compact set $\cD$, so for any
$\varepsilon>0$ there exists a $\delta>0$ such that
\[
|\mathbf{q}-\br|<\delta
\quad\implies\quad
| \LIFT\WORK(\mathbf{q})-\LIFT\WORK(\br) | <
\varepsilon/2\qquad
\mbox{for $\mathbf{q}, \br\in\cD$}.\vspace*{-2pt}
\]
If $\mathbf{q}\in\INV_\delta$, then it is within $\delta$ of some
$\br\in\INV$, hence
\begin{eqnarray*}
| \mathbf{q}- \LIFT\WORK(\mathbf{q}) |
&\leq&
|\mathbf{q}-\br| + | \br-\LIFT\WORK(\br) | +
|
\LIFT\WORK(\br)-\LIFT\WORK(\mathbf{q})|\\[-2pt]
&<&
\delta+ 0 + \varepsilon/2\\[-2pt]
&<&
\varepsilon
\qquad\mbox{for $\delta$ sufficiently small}.\vspace*{-2pt}
\end{eqnarray*}

\textit{Time to hit $\cK_\delta$.}\quad
Consider first the rate of change of $K(\cdot)$ while the process is
in $\cD\setminus\cK_\delta$
%
%
\begin{eqnarray}\label{eq:thissup}
\dot{K}(\mathbf{q}(t))
&\leq&
\dot{L}(\mathbf{q}(t))
=
\blambda\bdot f(\mathbf{q}(t)) - \max_{\bpi\in\sS} \bpi\bdot
f(\mathbf{q}(t))\nonumber\\[-2pt]
&\leq&
\sup_{\br\in\cD\setminus\cK_\delta} \Bigl[\blambda\bdot f(\br
) -
\max_{\bpi\in\sS} \bpi\bdot f(\br) \Bigr]\\[-2pt]
&\leq&0\qquad
\mbox{by Lemma \ref{lem:lyapunov}.}
\nonumber\vspace*{-2pt}
\end{eqnarray}
The supremum in (\ref{eq:thissup}) is of a continuous function of $\br$,
taken over a compact set;
hence the\vadjust{\goodbreak} supremum is attained at some
$\hbr\in\cD\setminus\cK_\delta$.
If the supremum were equal to 0, then
$\blambda\bdot f(\hbr)=\max_\bpi\bpi\bdot f(\hbr)$,
so
$\hbr\in\INV$
by Lemma \ref{lem:fixedpoint};
but $\hbr\in\cD\setminus\cK_\delta$, and we just proved that
$\INV\subset\cK_\delta$;
hence the supremum is some $-\eta_\delta<0$.

Now consider any fluid model solution starting at
$\mathbf{q}(0)$ with $|\mathbf{q}(0)|\leq1$. If \mbox{$\mathbf{q}(0)\in
\cK_\delta$}, then
it remains in $\cK_\delta$, so the theorem holds trivially. If not,
then $\mathbf{q}(0)\leq\bOne$ componentwise, so $L(\mathbf
{q}(0))\leq L(\bOne)$, so\vspace*{1pt}
$\mathbf{q}(0)\in\cD$; also $L(\mathbf{q}(t))$ is decreasing so
$\mathbf{q}(t)\in\cD$ for
all $t\geq0$. Now, $\dot{K}(\mathbf{q}(t))\leq-\eta_\delta$ all
the time
that $\mathbf{q}(t)\in\cD\setminus\cK_\delta$, and this cannot go
on for longer than
$H_\varepsilon= K(\mathbf{q}(0))/\eta_\delta\leq L(\bOne)/\eta
_\delta$.\vspace*{-3pt}
\end{pf}

\section{\texorpdfstring{Fluid model behavior (multi-hop case).}{Fluid
model behavior (multi-hop
case)}}\vspace*{-3pt}
\label{sec:fluid.behavior.m}

In this section we describe properties of fluid model solutions for a
multi-hop switched network running MW-$f$ back-pressure, as described in
Section \ref{sec:model}.

Let $R$ be the routing matrix and
$\tR=(I-R^\T)^{-1}$; recall that $\tR_{m n}=1$ if work injected at
queue $n$ eventually passes through $m$, and $0$ otherwise.
For a vector $\bx\in\Reals^N$, let $\tbx=\tR\bx$:
for arrival rate vector $\blambda$, $\tlambda_n$ is the total arrival
rate of work destined to pass through queue $n$; for a queue size
vector $\mathbf{q}$, $\tq_n$ is the total amount of work at queue $n$ and
queues upstream of $n$.

The set $\Lambda$, the $\PRIMAL(\cdot)$ and $\DUAL(\cdot)$ problems,
the set $\sS^*$ of virtual resources, and $\CLVR(\cdot)$ are
defined as in the single-hop case. The difference is that we will
require $\tblambda\in\Lambda$, and we will define the set of
critically loaded virtual resources to be $\CLVR(\tblambda)$. We also
need to modify the definition of $\ALGDUAL$ and the lifting map.\vspace*{-3pt}

\begin{definition}[(Lifting map)]
With $L\dvtx\RealsP^N\to\RealsP$ as in the single-hop case, define the
optimization problem $\ALGDUAL(\mathbf{q})$ to be
\begin{eqnarray*}
&&\mbox{minimize }
L(\br)\quad
\mbox{over}\quad\br\in\RealsP^N\\[-3pt]
&&\mbox{such that }
\bxi\bdot\tbr\geq\bxi\bdot\tbq\qquad\mbox{for all
$\bxi\in\CLVR(\tblambda)$}\quad\mbox{and}\\[-3pt]
&&\hphantom{\mbox{such that }}\vec{r}_n\leq\tq_n \qquad\mbox{for all $n$ such
that $\tlambda_n=0$}.
\end{eqnarray*}
Note that $L$ is strictly convex and increasing componentwise, and the
feasible region is convex; hence this problem has a unique
optimizer. Define the \textit{lifting map}
$\LIFTW\dvtx\RealsP^N\to\RealsP^N$ by setting $\LIFTW(\mathbf{q})$ to
be the optimizer.\vspace*{-3pt}
\end{definition}

The main result of this section is the following. Throughout this
section we are considering a multi-hop network with arrival rate
vector $\blambda\geq\bZero$ such that $\tblambda\in\Lambda$, running
MW-$f$ back-pressure.\vspace*{-3pt}
\begin{theorem}[(Portmanteau theorem, multi-hop version)]
\label{thm:fp.m}
The statements of Theorem \ref{thm:fp} parts \textup{(i)--(v)} hold,
for multi-hop fluid model solutions and
using the multi-hop definition of $\LIFTW$.\vspace*{-3pt}
\end{theorem}

Some of the proofs for the single-hop case carry through to the
multi-hop case.
Other proofs rely\vadjust{\goodbreak} on the fact that
for single-hop networks, $\blambda\in\Lambda\implies\blambda\leq
\bsigma$
for some $\bsigma\in\langle{\sS}\rangle$, and these proofs require
modification. We will modify them to use the following result.\vspace*{-2pt}

\begin{lemma}\label{lem:monotone}
Under Assumption \ref{cond.monotone}, if
$\bsigma\in\langle{\sS}\rangle$ and $\bsigma'\in\RealsP^N$ is
such that
$\bsigma'\leq\bsigma$, then $\bsigma'\in\langle{\sS}\rangle$.\vspace*{-2pt}
\end{lemma}
\begin{pf}
It is sufficient to establish the result for the
case when $\bsigma'$ differs from $\bsigma$ in only one
component, as the repeated application of this will yield the full
result. Without loss of generality, assume the queues are numbered
such that $0\leq\sigma'_1<\sigma_1$ and
$\sigma'_n=\sigma_n$ for $n\geq2$. Since $\bsigma\in\langle{\sS
}\rangle$
there is a collection of positive constants $(a_\bpi)_{\bpi\in\sS}$
such that $\sum_\bpi a_\bpi=1$ and
$\bsigma=\sum_\bpi a_\bpi\bpi$.
By Assumption \ref{cond.monotone}, if $\bpi\in\sS$, then $\bpi'\in
\sS$
where
\[
\pi'_n = \cases{
0, &\quad if $n=1$,\cr
\pi_n, &\quad otherwise,}
\]
thus
$\bsigma''\in\langle{\sS}\rangle$ where $\bsigma''=\sum_\bpi
a_\bpi\bpi'$.
By construction, $\sigma''_1=0$ and $\sigma''_n=\sigma_n$ for
$n\geq2$. By choosing the appropriate convex combination
\[
\bsigma' = (1-x)\bsigma'' + x\bsigma
\qquad\mbox{with }
x=\sigma'_1/\sigma_1\in[0,1]
\]
we see $\bsigma'\in\langle{\sS}\rangle$.\vspace*{-2pt}
\end{pf}

Now we proceed toward establishing Theorem \ref{thm:fp.m}.
The proof of the first claim of Theorem
\ref{thm:fp.m}(i) is just as for the single-hop
case. The second claim follows from the following lemma.\vspace*{-2pt}

\begin{lemma}
\label{lem:lyapunov.multihop}
For all $\mathbf{q}\in\RealsP^N$,
%
%
\begin{equation}
\label{eq:m.lyapunov.neg}
\blambda\bdot f(\mathbf{q}) - \max_{\bpi\in\sS} \bpi\bdot
(I-R)f(\mathbf{q})\leq0.
\end{equation}
Also, every fluid model solution satisfies
\[
\frac{d}{dt}L(\mathbf{q}(t)) = \blambda\bdot f(\mathbf{q}(t)) -
\max_{\bpi\in\sS} \bpi\bdot
(I-R)f(\mathbf{q}(t))
\leq0.\vspace*{-2pt}
\]
\end{lemma}

\begin{pf}
Since $\tblambda\in\Lambda$, $\tR\blambda\leq\bsigma$
componentwise for some $\bsigma\in\langle{\sS}\rangle$.
Because $\tR\geq0$ and $\blambda\geq\bZero$,
$\tR\blambda\geq\bZero$ componentwise. By Lemma \ref{lem:monotone},
$\blambda=(I-R^\T)\bsigma'$ for some $\bsigma'\in\langle{\sS
}\rangle$.
Hence
\[
\blambda\bdot f(\mathbf{q}) - \max_{\bpi\in\sS}\bpi\bdot
(I-R)f(\mathbf{q})
=
\bsigma'\bdot(I-R)f(\mathbf{q}) - \max_{\bpi\in\sS}\bpi\bdot
(I-R)f(\mathbf{q})
\leq0.
\]
For the claim about fluid model solutions,
\begin{eqnarray*}
\frac{d}{dt} L(\mathbf{q}(t))
&=&
\dot{\mathbf{q}}(t) \bdot f(\mathbf{q}(t))\\
&=&
\biggl(
\blambda-(I-R^\T)\biggl[\sum_\bpi\dot{s}_\bpi(t)\bpi+\dot
{\mathbf{y}}(t)\biggr]\biggr)
\bdot f(\mathbf{q}(t))\\
&&\mbox{by differentiating (\ref{eq:multihop.fluid.queue})}\\
&=&
\blambda\bdot f(\mathbf{q}(t)) - \sum_\bpi\dot{s}_\bpi(t)
\bpi\bdot(I-R)f(\mathbf{q}(t)) + \dot{\mathbf{y}}(t)\bdot
(I-R)f(\mathbf{q}(t)).
\end{eqnarray*}
For the middle term,
\begin{eqnarray*}
\sum_\bpi\dot{s}_\bpi(t)\bpi\bdot(I-R)f(\mathbf{q}(t))
&=&
\max_\brho\brho\bdot(I-R)f(\mathbf{q}(t)) \sum_\bpi\dot{s}_\bpi(t)
\qquad\mbox{by (\ref{eq:multihop.fluid.mwm})}\\
&=&
\max_\brho\brho\bdot(I-R)f(\mathbf{q}(t))
\qquad\mbox{by (\ref{eq:fluid.busy})}.
\end{eqnarray*}
For the last term, we claim that
%
%
\begin{equation}
\label{eq:multihop.noidle}
\dot{\mathbf{y}}(t)\bdot(I-R)f(\mathbf{q}(t))=0.
\end{equation}
To see this, consider first a queue $n$ with
$[(I-R)f(\mathbf{q}(t))]_n>0$. As noted in (\ref{eq:backpressure}), this
implies $f(q_n(t))>f([R\mathbf{q}(t)]_n)$. By Assumption
\ref{multihop.cond.f} it must be that $q_n(t)>0$, hence
$\dot{y}_n(t)=0$ by (\ref{eq:fluid.idling}).
Second, consider a queue $n$ with $[(I-R)f(\mathbf{q}(t))]_n<0$.
It must be that all of the active schedules do not serve this
queue, that is, $\dot{s}_\bpi(t)>0 \implies\pi_n=0$, since
otherwise by
Assumption \ref{cond.monotone} there is another schedule that has
bigger weight than $\bpi$, contradicting
(\ref{eq:multihop.fluid.mwm}). Third, if $[(I-R)f(\mathbf
{q}(t))]_n=0$ then
obviously $\dot{y}_n(t)[(I-R)f(\mathbf{q}(t))]_n=0$. Putting these three
together proves (\ref{eq:multihop.noidle}).

Putting together these findings for the middle and last terms,
\[
\frac{d}{dt} L(\mathbf{q}(t))
=
\blambda\bdot f(\mathbf{q}(t))
- \max_\brho\brho\bdot(I-R)f(\mathbf{q}(t)).
\]
Applying (\ref{eq:m.lyapunov.neg}) this is $\leq0$.\vadjust{\goodbreak}
\end{pf}

The proof of Theorem \ref{thm:fp.m}(ii) is broadly
similar to the
single-hop case, Lem\-ma~\ref{lem:lift.continuous}, but the formulae all
have to be adjusted to deal with multi-hop.

\begin{lemma}
$\LIFTW$ is continuous.
\end{lemma}
\begin{pf}
If $\CLVR(\tblambda)$ is empty, then the lifting map is trivial, and
the result is trivial. In what follows, we shall assume that
$\CLVR(\tblambda)$ is nonempty, and we will abbreviate it to
$\CLVR$. Furthermore\vspace*{1pt} note that for every
$\bxi\in\CLVR$ we know $\bxi\bdot\tblambda=1$
by definition of $\CLVR$,
and hence there is some
queue $n$ such that $\xi_n>0$ and $\tlambda_n>0$.

Pick any sequence $\mathbf{q}^k\tend\mathbf{q}$, and let $\br
^k=\LIFTW(\mathbf{q}^k)$ and
$\br=\LIFTW(\mathbf{q})$. We want to prove that $\br^k\tend\br$.
We shall
first prove that there is a compact set $[0,h]^N$ such that
$\br^k\in[0,h]^N$ for all $k$. We shall then prove that any
convergent subsequence of $\br^k$ converges to $\br$; this
establishes continuity of~$\LIFTW$.

First---compactness. A suitable value for $h$ is
\[
h = \max_{\bxi\in\CLVR} \max_{n\dvtx\xi_n>0} \sup_k \frac{\bxi
\bdot\tbq^k}{\xi_n}.\vadjust{\goodbreak}
\]
Note that the maximums are over a nonempty set, as noted at the
beginning of the proof. Note also that $h$ is finite because $\mathbf
{q}$ is
finite. Now, suppose that $\br^k\notin[0,h]^N$ for some $k$,
that is, that there is some queue $n$ for which $r^k_n>h$, and let
$\br'=\br^k$ in every coordinate except for $r'_n=h$. We claim that
$\br'$ satisfies the two constraints of $\ALGDUAL(\mathbf{q}^k)$. To
see that
it satisfies the second constraint, note that $\br'\leq\br^k$, and
hence if $\tlambda_n=0$, then $\vec{r} '_n\leq\vec{r} {}^k_n\leq
\tq_n$. To see that
it satisfies the first constraint, pick any $\bxi\in\CLVR$.
Either\vspace*{1pt}
$\xi_m=0$ for all queues~$m$ that are downstream of $n$, that is, for
which $\tR_{m n}=1$; if this is so, then\looseness=-1
\[
\bxi\bdot\tbr' = \bxi\bdot\tbr^k + \bxi\bdot(\tbr'-\tbr^k)
= \bxi\bdot\tbr^k + \sum_l(r'_l-r^k_l)\sum_m\xi_m\tR_{m l} =
\bxi\bdot\tbr^k.
\]\looseness=0
Or $\xi_m>0$ for some queue $m$ that is downstream of $n$; if this is
so, then
\[
\bxi\bdot\tbr' \geq\xi_m\vec{r}'_m \geq\xi_m h \geq\bxi\bdot
\tbq^k\qquad
\mbox{by construction of $h$}.\vadjust{\goodbreak}
\]
Applying this repeatedly, if
$\br^k\notin[0,h]^N$, then we can reduce it to a queue size vector in
$[0,h]^N$, thereby improving on $L(\br^k)$, yet still meeting the
constraints of $\ALGDUAL(\mathbf{q}^k)$; this contradicts the
optimality of
$\br^k$. Hence $\br^k\in[0,h]^N$.

Next---convergence on subsequences. With a slight abuse of notation,
let $\LIFTW(\mathbf{q}^k)=\br^k\tend\bs$ be a convergent
subsequence, and
recall that $\LIFTW(\mathbf{q})=\br$ and $\mathbf{q}^k\tend\mathbf
{q}$. By continuity of
the constraints of $\ALGDUAL$, $\bs$ is feasible for $\ALGDUAL
(\mathbf{q})$;
we shall next show that $L(\bs)\leq L(\br)$. Since $\br$ is the unique
optimum, it must be that $\bs=\br$.

It remains to show that $L(\bs)\leq L(\br)$. We will construct a
sequence
$\br-\bdelta^k+\varepsilon^k\bP$ of candidate solutions to
$\ALGDUAL(\mathbf{q}^k)$, choosing\vadjust{\goodbreak} $\bdelta^k\geq\bZero$ and
$\varepsilon^k\bP\geq\bZero$ to ensure
that the candidate solutions are feasible. Specifically, we define
\[
\delta^k_n = \cases{
0, &\quad if $\tlambda_n>0$,\cr
(\tq_n-\tqsup^k_n)^+, &\quad if $\tlambda_n=0$,}
\]
and $P_n=1_{\tlambda_n>0}$, and
\[
\varepsilon^k = \max_{\bxi\in\CLVR}
\frac{(\bxi\bdot\tbq^k-\bxi\bdot\tbq)^++\bxi\bdot
\tbdelta^k}{\bxi\bdot\tbP}.
\]
We will first deal with the
feasibility constraint that pertains when $\tlambda_n=0$. Note that
this implies $\lambda_m=0$ for all queues $m$ that are upstream of
$n$, since $\tlambda_n=\sum_m \tR_{n m}\lambda_m$, and hence that
$\tlambda_m=0$ for all upstream queues.
Using this we find
\begin{eqnarray*}
&&[\tR(\br-\bdelta^k+\varepsilon^k\bP)]_n\\
&&\qquad=
\sum_m \tR_{n m}[\br-\bdelta^k+\varepsilon^k\bP]_m\\
&&\qquad=
\sum_m \tR_{n m}\bigl(r_m - (\tq_m-\tqsup^k_m)^+ \bigr)
\qquad\mbox{since
$\tlambda_m=0$ when $\tR_{n m}=1$}\\
&&\qquad=
\biggl(\sum_m \tR_{n m}r_m \biggr) - \biggl(\sum_m \tR_{n m} (\tq
_m-\tqsup^k_m)^+ \biggr) \\
&&\qquad\leq
\biggl(\sum_m \tR_{n m}r_m\biggr) - (\tq_n-\tqsup^k_n)^+ \qquad\mbox{as
$\tR_{nn}=1$, $\tR_{n m} \geq0$ for all $m$}\\
&&\qquad=
\vec{r}_n - (\tq_n-\tqsup^k_n)^+\\
&&\qquad\leq
\tq_n - (\tq_n-\tqsup^k_n)^+ \qquad\mbox{since $\br$ is feasible for
$\ALGDUAL(\mathbf{q})$}\\
&&\qquad= \min(\tq_n,\tqsup^k_n) \leq\tqsup^k_n.
\end{eqnarray*}
Hence $\br-\bdelta^k+\varepsilon^k\bP$ satisfies the second
feasibility constraint of
$\ALGDUAL(\mathbf{q}^k)$.
For the other feasibility constraint of $\ALGDUAL(\mathbf{q}^k)$,
pick any
$\bxi\in\CLVR$. Then
\begin{eqnarray*}
&&\bxi\bdot\tR(\br-\bdelta^k+\varepsilon^k\bP)\\
&&\qquad=
\bxi\bdot(\tbr-\tbdelta^k) + \varepsilon^k\bxi\bdot\tbP\\
&&\qquad\geq
\bxi\bdot(\tbr-\tbdelta^k) +
( \bxi\bdot\tbq^k-\bxi\bdot\tbq)^++\bxi\bdot
\tbdelta^k\qquad
\mbox{by construction of $\varepsilon^k$}\\
&&\qquad\geq
\bxi\bdot\tbq- (\bxi\bdot\tbq^k-\bxi\bdot\tbq)^+\qquad
\mbox{since $\tbr$ is feasible for $\ALGDUAL(\tbq)$}\\
&&\qquad=
\max(\bxi\bdot\tbq,\bxi\bdot\tbq^k) \geq\bxi\bdot\tbq^k.
\end{eqnarray*}
Since the candidates are feasible solutions to $\ALGDUAL(\mathbf
{q}^k)$, and
$\br^k$ is an optimal solution, it must be that
\[
L(\br^k) \leq L(\br-\bdelta^k+\varepsilon^k\bP).
\]
Taking the limit as $k\tend\infty$, and noting that $L$ is
continuous
and $\bdelta^k\tend\bZero$ and \mbox{$\varepsilon^k\tend0$}, we find
\[
L(\bs) \leq L(\br)
\]
as required. This completes the proof.
\end{pf}

For the proof of Theorem \ref{thm:fp.m}(iii), it is useful
to work with a different representation of $\LIFTW$, provided by the
following lemma, which draws on monotonicity of~$\sS$.
\begin{lemma}
\label{lem:monotoneDW}
For any $\mathbf{q}\in\RealsP^N$,
$\LIFTW(\mathbf{q})$ can be written
\[
\LIFTW(\mathbf{q}) = \mathbf{q}+ t\bigl(\blambda-(I-R^\T)\bsigma
\bigr)\qquad
\mbox{for some } t\geq0, \bsigma\in\langle{\sS}\rangle.
\]
\end{lemma}
\begin{pf}
we will choose $\bsigma$ simply by multiplying each side of the
desired equation by $\tR$
\[
\tbr= \tbq+ t(\tblambda-\bsigma) \qquad\mbox{where }
\br=\LIFTW(\mathbf{q})\vadjust{\goodbreak}
\]
or, rearranging,
\[
\bsigma= \tblambda-(\tbr-\tbq)/t.
\]
We will\vspace*{1pt} show that $\bZero\leq\bsigma\leq\brho$ for some
$\brho\in\langle{\sS}\rangle$, hence by Lemma
\ref{lem:monotone} $\bsigma\in\langle{\sS}\rangle$.

First, we show $\bZero\leq\bsigma$. If $\tlambda_n>0$ this can be
achieved by choosing $t$ sufficiently large. If $\tlambda_n=0$, then by
the second constraint of $\ALGDUAL(\mathbf{q})$ we know that $\vec
{r}_n\leq\tq_n$
so $\sigma_n\geq0$.

Second, we show $\bsigma\bdot\bxi\leq1$ for all $\bxi$ that are
feasible for $\DUAL(\tblambda)$.
Either \mbox{$\bxi\bdot\tblambda=1$}, in which case
$\bxi\in\langle{\CLVR(\tblambda)}\rangle$ and so by\vspace*{2pt} the first constraint
of $\ALGDUAL$
we know that $\tbr\bdot\bxi\geq\tbq\bdot\bxi$.
Or $\bxi\bdot\tblambda<1$, in which case we simply need to choose $t$
sufficiently large. Either way, $\bsigma\bdot\bxi\leq1$ for all
dual-feasible $\bxi$, hence $\DUAL(\bsigma)\leq1$, hence
$\PRIMAL(\bsigma)\leq1$, hence $\bsigma\leq\brho$ for some
$\brho\in\langle{\sS}\rangle$ by the definition of $\PRIMAL
(\bsigma)$.
\end{pf}

The proof of Theorem \ref{thm:fp.m}(iii) is given by the
following lemma.
This proof is similar to the single-hop case,
Lemma \ref{lem:lift.scaleinvariant}, but it is much shorter
because the monotonicity assumption gives us a stronger
representation of the lifting map, Lem\-ma~\ref{lem:monotoneDW}.
Also, this version makes a weaker claim, namely that the lifting map
is scale-invariant at invariant states, whereas the single-hop version
shows that the lifting map is invariant everywhere.
\begin{lemma}[(Scale-invariance of the lifting map)]
\label{lem:multihop.scale}
If $\mathbf{q}=\LIFTW(\mathbf{q})$ then $\kappa\mathbf{q}=\LIFTW
(\kappa\mathbf{q})$
for all $\kappa>0$.
\end{lemma}
\begin{pf}
Suppose that $\mathbf{q}=\LIFTW(\mathbf{q})$, and let
$\kappa\br=\LIFTW(\kappa\mathbf{q})$. Clearly $\kappa\mathbf{q}$
is feasible for
$\ALGDUAL(\kappa\mathbf{q})$; we shall show that $L(\kappa\br)\geq
L(\kappa\mathbf{q})$, whence $\kappa\mathbf{q}$ is also optimal for
$\ALGDUAL(\kappa\mathbf{q})$, whence $\kappa\mathbf{q}=\kappa\br
$ by uniqueness of
the optimum.

It remains to prove that $L(\kappa\br)\,{\geq}\,L(\kappa\mathbf{q})$. By
Lemma \ref{lem:monotoneDW}, we can write $\kappa\br$~as
\[
\kappa\br= \kappa\mathbf{q}+ t\bigl(\blambda-(I-R^\T)\bsigma
\bigr)\vadjust{\goodbreak}
\]
for some $t\geq0$ and some $\bsigma\in\Sigma$. Now consider the value
of $L$ along a straight-line trajectory from $\kappa\mathbf{q}$ to
$\kappa\br$
\begin{eqnarray*}
&&\frac{d}{du} L\bigl(\kappa\mathbf{q}+(\blambda-(I-R^\T)\bsigma
)u\bigr) \bigg|_{ u=0}\\
&&\qquad=
\bigl(\blambda-(I-R^\T)\bsigma\bigr)\bdot f(\kappa\mathbf{q})\\
&&\qquad=
\blambda\bdot f(\kappa\mathbf{q}) - \bsigma\bdot(I-R)f(\kappa
\mathbf{q})\\
&&\qquad\geq
\blambda\bdot f(\kappa\mathbf{q}) - \max_{\rho\in\sS}\brho\bdot
(I-R)f(\kappa\mathbf{q})\qquad
\mbox{for any $\sigma\in\langle{\sS}\rangle$}\\
&&\qquad=0.
\end{eqnarray*}
The final equality is because
\[
\blambda\bdot f(\mathbf{q})
=\tblambda\bdot(I-R)f(\mathbf{q})=
\max_{\bpi\in\sS}\bpi\bdot(I-R)f(\mathbf{q})\vadjust{\goodbreak}
\]
by Lemma \ref{lem:multihop.fixed}(iv)
below (the proof of which does not assume the result of this lemma).
Hence
\[
\blambda\bdot f(\kappa\mathbf{q})
=
\tblambda\bdot(I-R)f(\kappa\mathbf{q})
=
\max_{\bpi\in\sS}\bpi\bdot(I-R)f(\kappa\mathbf{q})
\]
using Assumption \ref{multihop.cond.f} and the fact that
$\tblambda\in\langle{\sS}\rangle$ by Lemma \ref{lem:monotone}.
\end{pf}

The proof of Theorem \ref{thm:fp.m}(iv) relies on the
following lemma.
\begin{lemma}[(Fluid model trajectories preserve $\ALGDUAL$ feasibility)]
\label{lem:multihop.feasible}
Consider any fluid model solution, for any scheduling policy, with
initial queue size $\mathbf{q}(0)$. Then $\mathbf{q}(t)$ is feasible for
$\ALGDUAL(\mathbf{q}(0))$ for all $t\geq0$.
\end{lemma}
\begin{pf}
Feasibility\vspace*{1pt} for $\ALGDUAL(\mathbf{q}(0))$ has two parts. For the
first part,
pick any critically loaded virtual resource\vspace*{1pt}
$\bxi\in\CLVR(\tblambda)$, and multiply each side of
(\ref{eq:multihop.fluid.queue}) by $\tR=(I-R^T)^{-1}$ and then by
$\bxi$ to get
\[
\bxi\bdot\tbq(t)
=
\bxi\bdot\tbq(0) + \bxi\bdot\tR\mathbf{a}(t) - \bxi\cdot
\biggl(\sum_\bpi
s_{\bpi}(t)\bpi- \mathbf{y}(t)\biggr).
\]
Defining $\bsigma(t)=\sum\bpi s_\bpi(t)/t$, which is in $\langle
{\sS}\rangle$
by (\ref{eq:fluid.busy}),
\begin{eqnarray*}
\bxi\bdot\tbq(t)
&\geq&
\bxi\bdot\tbq(0) + t\bigl(\bxi\bdot\tblambda-
\bxi\bdot\bsigma(t)\bigr)
\qquad\mbox{by (\ref{eq:fluid.arrivals}) and because $\mathbf
{y}(t)\geq\bZero$}\\[-2pt]
&=&
\bxi\bdot\tbq(0) + t\bigl(1 -
\bxi\bdot\bsigma(t)\bigr)
\qquad\mbox{since $\bxi\in\CLVR(\tblambda)$}\\[-2pt]
&\geq&
\bxi\bdot\tbq(0) + t(1 -
1)
\qquad\mbox{since $\bxi$ is a virtual resource and $\bsigma\in
\langle{\sS}\rangle$}\\[-2pt]
&=&
\bxi\bdot\tbq(0)
\end{eqnarray*}
as required for the first part of $\ALGDUAL$-feasibility. For the
second part,
suppose that $\tlambda_n=0$ for some queue $n$. Multiply
each side of (\ref{eq:multihop.fluid.queue}) by $\tR$ to get
\[
\tbq(t)
=
\tbq(0) + \tblambda t
- \sum_\bpi\bpi s_\bpi(t) + \mathbf{y}(t)
\leq
\tbq(0) + \tblambda t,\vadjust{\goodbreak}
\]
where the inequality is by (\ref{ineq:fluid.littleidle}). Since we
assumed $\tlambda_n=0$, $\tq_n(t)\leq\tq_n(0)$.
This completes the proof that $\mathbf{q}(t)$ is feasible for
$\ALGDUAL(\mathbf{q}(0))$.
\end{pf}

The proof of Theorem \ref{thm:fp.m}(iv) is implied by
parts (i) and (ii) of the following lemma.
\begin{lemma}[(Characterization of invariant states of MW-$f$ backpressure)]
\label{lem:multihop.fixed}
The following are equivalent, for $\mathbf{q}^0\in\RealsP^N$:
\begin{longlist}
\item
$\mathbf{q}^0 = \LIFTW(\mathbf{q}^0)$;\vspace*{1pt}
\item
$\mathbf{q}^0$ is an invariant state;
\item
there exists a fluid model solution with $\mathbf{q}(t)=\mathbf{q}^0$
for all $t$;
\item
$\blambda\bdot f(\mathbf{q}^0) = \max_{\bpi\in\sS} \bpi\bdot
(I-R)f(\mathbf{q}^0)$.\vadjust{\goodbreak}
\end{longlist}
\end{lemma}
\begin{pf}
That (i)${}\implies{}$(ii)${}\implies{}$(iii)${}\implies{}$(iv) is proved in the
same way as in the single-hop case.
We just need to appeal to Lemma \ref{lem:lyapunov.multihop} rather
than \ref{lem:lyapunov} for the fact that $L(\mathbf{q}(t))$ is
decreasing, and to
Lemma \ref{lem:multihop.feasible} rather than \ref{lem:feasible} for
the fact that $\mathbf{q}(t)$ remains feasible.\vspace*{8pt}

\textit{Proof that} (iv)${}\implies{}$(i).\quad
Let $\br=\LIFTW(\mathbf{q}^0)$. By Lemma \ref{lem:monotoneDW},
$\br=\mathbf{q}^0+t(\blambda-(I-R^\T)\bsigma)$ for some $t\geq0$ and
$\bsigma\in\langle{\sS}\rangle$.
By considering the value of $L(\cdot)$
along the trajectory from $\mathbf{q}^0$ to $\br$,
and using (iv), we conclude that $L(\br)\geq L(\mathbf{q}^0)$.
By the same argument as in the single-hop case,
$\mathbf{q}^0=\LIFTW(\mathbf{q}^0)$.
\end{pf}

The proof of Theorem \ref{thm:fp.m}(v) is given by the
following lemma.
\begin{lemma}
\label{lem:multihop.conv}
Given $\tblambda\in\Lambda$,
for any $\varepsilon>0$ there exists an $H_\varepsilon>0$ such
that for every fluid model solution with arrival rate $\blambda$,
for which $|\mathbf{q}(0)|\leq1$,
$|\mathbf{q}(t)-\LIFTW(\mathbf{q}(t))|<\varepsilon$ for all $t\geq
H_\varepsilon$.
\end{lemma}
\begin{pf}
The proof of Lemma \ref{lem:criticalfluid.time} goes through almost
verbatim. The only changes are in the penultimate paragraph, which
should be replaced by the following:\vspace*{8pt}

\textit{Time to hit $\cK_\delta$.}\quad
Consider first the rate of change of $K(\cdot)$ while the process is
in $\cD\setminus\cK_\delta$
\begin{eqnarray*}
\dot{K}(\mathbf{q}(t))
&\leq&
\dot{L}(\mathbf{q}(t))
=
\blambda\bdot f(\mathbf{q}(t)) - \max_{\bpi\in\sS} \bpi\bdot
(I-R)f(\mathbf{q}(t))\\
&\leq&
\sup_{\br\in\cD\setminus\cK_\delta} \Bigl[\blambda\bdot f(\br
) -
\max_{\bpi\in\sS} \bpi\bdot(I-R)f(\br) \Bigr]\\
&\leq&0\qquad
\mbox{by Lemma \ref{lem:lyapunov.multihop}.}
\end{eqnarray*}
This supremum is of a continuous function of $\br$,
taken over a closed
and bounded set, hence the supremum is attained at some
$\hbr\in\cD\setminus\cK_\delta$.
If the supremum were equal to 0, then
$\blambda\bdot f(\hbr)=\max_\bpi\bpi\bdot(I-R)f(\hbr)$
so
$\hbr\in\INV$
by Lemma \ref{lem:multihop.fixed};
but $\hbr\in\cD\setminus\cK_\delta$ and we just proved that
$\INV\subset\cK_\delta$;
hence the supremum is some $-\eta_\delta<0$.
\end{pf}

\section{\texorpdfstring{Multiplicative state-space
collapse.}{Multiplicative state-space
collapse}}
\label{sec:ht}

This section establishes multiplicative state space collapse of queue
size. It shows that under the MW-$f$ policy, and with suitable
initial conditions when the network is not overloaded (i.e., when
$\blambda\in\Lambda$), the appropriately normalized queue size
vector is constrained to lie in or close to the \textit{set of invariant states}
\[
\INV= \{\mathbf{q}\in\RealsP^\N\dvtx\mathbf{q}=\LIFT\WORK
(\mathbf{q})\}.
\]
We assume that arrivals satisfy Assumption \ref{cond:mssc.stoch}, and
let the arrival rate vector $\blambda$ be as specified in that
assumption. The function $\LIFT\WORK$ depends on $\blambda$ and $f$, as\vadjust{\goodbreak}
specified in Sections \ref{sec:fluid.behavior} and
\ref{sec:fluid.behavior.m} for single-hop and multi-hop networks, respectively,
and the interesting case is where $\blambda\in\dLambda$ (since
otherwise $\LIFT\WORK$ is trivial).

This section
mostly follows the method developed by Bramson \cite{bramson},
except that our proof avoids the need for regenerative assumptions on
the arrival process by imposing slightly tighter bounds on the
uniformity of
their convergence, as expressed by Assumption
\ref{cond:mssc.stoch}.

Consider a sequence of systems of the type described in Section
\ref{sec:model.queue} running a scheduling policy of the type
described in Section \ref{sec:model.alg}.
Let the systems all have the same number of queues $N$, the same set
of allowed schedules $\sS$, the same routing matrix $R$ and the same
scheduling policy.
Let the sequence of
systems be indexed by $r\in\Naturals$.
Write
\[
X^r(\tau)=(\bQ^r(\tau),\bA^r(\tau),\bZ^r(\tau),S^r(\tau)),\qquad
\tau\in\IntegersP,
\]
for the $r$th system.
Define the scaled system $\hat{x}^r(t)=(\hbq^r(t),\hba^r(t),\hbz
^r(t),\hs^r(t))$
for $t\in\RealsP$ by
\begin{eqnarray*}
\hbq^r(t) &=& \bQ^r(r^2t)/r,\qquad
\hba^r(t) = \bA^r(r^2t)/r,\\
\hbz^r(t) &=& \bZ^r(r^2t)/r,\qquad
\hs^r_\bpi(t) = S^r_\bpi(r^2t)/r
\end{eqnarray*}
after extending the domain of $X^r(\cdot)$ to
$\RealsP$ by linear interpolation in each interval $(\tau,\tau+1)$.
Note that each sample path of a scaled system $\hat{x}^r(t)$ over
the interval $t\in[0,T]$ lies in $C^I(T)$ with $I=3N+|\sS|$. $T > 0$ will
be fixed for the remainder of this section. Recall the norm
$\|x\|={\sup_{0\leq t\leq T}}|x(t)|$. The main result of this paper is
the following.
\begin{theorem}[(Multiplicative state-space collapse)]
\label{thm:heavytraffic}
Consider a sequence of (single-hop or multi-hop) switched networks
indexed by $r \in\Naturals$, operating under the MW-$f$ policy
(with $f$ satisfying Assumptions \ref{cond.f} or \ref{multihop.cond.f}
and $\sS$ with Assumption~\ref{cond.monotone}),
as described above. Assume that the arrival processes satisfy
Assumption~\ref{cond:mssc.stoch} with $\blambda\in\Lambda$. Also
assume that the initial
queue sizes are nonrandom, and satisfy
$\lim_{r\tend\infty} \hbq^r(0)=\hbq_0$ for some $\hbq_0\in\INV$.
Then\vadjust{\goodbreak} for any $\delta>0$,
%
%
\begin{equation}
\label{eq:thm.mssc}
\Prob\biggl(
\frac{ \| \hbq^r(\cdot) - \LIFT\WORK(\hbq^r(\cdot))
\| }
{\|\hbq^r(\cdot)\| \mmax1}
< \delta
\biggr)
\tend1\qquad
\mbox{as $r\tend\infty$}.
\end{equation}
\end{theorem}

Simulations suggest that a stronger result holds in the widely-studied
\textit{diffusion} or \textit{heavy traffic} scaling, $\blambda^r=\blambda
-\Gamma/r$
for some nontrivial $\Gamma\in\RealsP^N$ and $\blambda\in\dLambda$.
We conjecture the following.
\begin{conjecture}
\label{conj:heavytraffic.ssc}
Under the assumptions of Theorem \ref{thm:heavytraffic} and
the additional assumption that increments in the arrival process
are i.i.d. and uniformly bounded, under the diffusion scaling for
any $\delta> 0$
%
%
\begin{equation}
\label{eq:thm.ssc}
\Prob\bigl( \| \hbq^r(\cdot) - \LIFT\WORK(\hbq^r(\cdot))
\|
< \delta
\bigr)
\tend1\qquad
\mbox{as $r\tend\infty$}.
\end{equation}
\end{conjecture}

%
%
\begin{figure}

\includegraphics{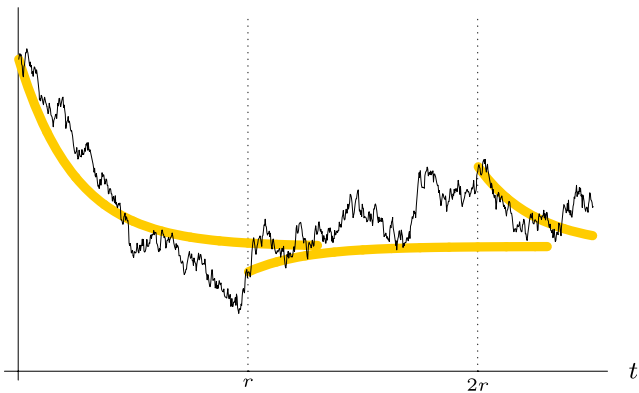}

\caption{Splitting the process into fluid-scaled parts, starting at
$0,r,2r,\ldots.$}
\label{fig:htscale}
\end{figure}

\subsection{\texorpdfstring{Outline of the proof of Theorem \protect\ref
{thm:heavytraffic}.}{Outline of the proof of Theorem 7.1}}

The outline of the proof of Theorem~\ref{thm:heavytraffic}
is as follows. We are interested in the dynamics of
$\hbq^r(t)$ over $t\in[0,T]$, that is, of $\bQ^r(\tau)$ over
$\tau\in[0,r^2T]$. We will split this time interval into
$\lfloor rT\rfloor+1$ pieces starting at $0,r,2r,\ldots,$ and look\vspace*{1pt} at each
piece under a~fluid scaling. We will define
a ``good event'' $\hE_r$ under
which the arrivals in all of the pieces are well behaved (Section \ref
{sec:ht.good}). We then
apply Theorem \ref{thm:fluid} to deduce that, under this event, the queue
size process in each of the pieces can be (uniformly) approximated by a fluid
model solution (Lemma \ref{lem:ht.fms}). We then use the properties of
the fluid model solution
stated in Theorem~\ref{thm:fp} to show that in each of the pieces, the
queue size is (uniformly) close to the set of invariant states (Lemmas
\ref{lem:ht.m} and \ref{lem:ssc.pathwise}).
Figure \ref{fig:htscale} depicts the idea. Finally we show that
$\Prob(\hE_r)\tend1$ (Lemma \ref{lem:ht.prob}). The formal proof is
given in Section~\ref{sec:ht.proof}.

Note that Lemmas \ref{lem:ht.fms}--\ref{lem:ssc.pathwise} are all
sample path-wise results that hold for every $\omega\in\hE_r$, and so
questions of independence etc. do not arise. The only part of the proof where
probability comes in is Lemma \ref{lem:ht.prob}.

The proof is written out for a single-hop switched network. For the
multi-hop case, the argument holds verbatim; simply replace all
references to the single-hop fluid limit Theorem \ref{thm:fluid} by
references to the equivalent multi-hop result,
and replace all references to the description of
single-hop fluid model solutions in Theorem \ref{thm:fp} by references
to the multi-hop version Theorem~\ref{thm:fp.m}.

\subsubsection{\texorpdfstring{The good event and the fluid-scaled
pieces.}{The good event and the fluid-scaled
pieces}}
\label{sec:ht.good}
Define the fluid-scaled pieces $\fx^{r,m,z}(u) = (\fbq^{r,m,z}(u),
\fba^{r,m,z}(u),\fby^{r,m,z}(u), \fls^{r,m,z}(u))$ of the original
process by
\begin{eqnarray*}
\fbq^{r,m,z}(u) &=& \bQ^r(rm+zu)/z,\\
\fba^{r,m,z}(u) &=& \bigl(\bA^r(rm+zu)-\bA^r(rm)\bigr)/z,\\
\fby^{r,m,z}(u) &=& \bigl(\bY^r(rm+zu)-\bY^r(rm)\bigr)/z,\\
\fls^{r,m,z}(u) &=& \bigl( S^r(rm+zu)- S^r(rm)\bigr)/z
\end{eqnarray*}
for $0\leq m\leq\lfloor rT\rfloor$, $z\geq r$ and $u \geq0$. Here $r$
indicates which process we are considering, $m$ indicates the piece
and $z$ indicates the fluid-scaling parameter. The scaling parameter
$z_{r,m}=|\bQ^r(rm)|\mmax r$ is particularly important, and for
convenience we will define $\fx^{r,m}(u) =
(\fbq^{r,m}(u), \fba^{r,m}(u), \fby^{r,m}(u), \fls^{r,m}(u))$
by\looseness=-1
\begin{eqnarray*}
\fbq^{r,m}(u) &=& \fbq^{r,m,z_{r,m}}(u),\qquad
\fba^{r,m}(u) = \fba^{r,m,z_{r,m}}(u), \\
\fby^{r,m}(u) &=& \fby^{r,m,z_{r,m}}(u),\qquad
\fls^{r,m}(u) = \fls^{r,m,z_{r,m}}(u).
\end{eqnarray*}\looseness=0
The good event is defined to be
%
%
\begin{eqnarray}
\label{eq:ht.good}
\hE_r &=&
\Bigl\{
\sup_{u\in[0,\Tfluid]}| \fba
^{r,m,w_{r,k}}(u)-\blambda^r
u|<\eta_r
\mbox{ for all $0\leq m\leq\lfloor rT\rfloor$}\nonumber\\[-9pt]\\[-9pt]
&&\hspace*{19pt}\mbox{and $0\leq k\leq\lfloor Lr\log r\rfloor$, where
$w_{r,k}=r(1+k/\log r)$}
\Bigr\}.
\nonumber
\end{eqnarray}
By this, we mean that $\hE_r$ is a subset of the sample space for the
$r$th system, and we write $\fx^{r,m,w_{r,k}}(\cdot)(\omega)$ etc. for
$\omega\in\hE_r$ when we wish to emphasize the dependence on~$\hE_r$.
The constants here are
$\Tfluid=(2+\lambda^{\max}+N S^{\max})(H_\zeta+1)$,
$\lambda^{\max}=\sup_r|\blambda^r|$,
$\zeta>0$ is chosen as specified in Section~\ref{sec:ht.proof} below,
$H_\zeta$ is chosen as in Theorem \ref{thm:fp}(v),
$L=1+T(1+\lambda^{\max}+N S^{\max})$,
$S^{\max}=\max_{\bpi\in\sS}|\bpi|$
and
the sequence of deviation terms $\eta_r\in[0,1]$ is chosen as
specified in Lemma \ref{lem:ht.prob}
such that $\eta_r\tend0$ as $r\tend\infty$.
\begin{lemma}
\label{lem:ht.fms}
Let $\FMS$ be the set of fluid model solutions
over time horizon $[0,\Tfluid]$ for arrival
rate vector $\blambda$, and let $\FMS(\mathbf{q}_0)$ and $\FMS_1$
be as
specified in Definition \ref{def:fms}.
Then
%
%
\begin{equation}\label{eq:ht.fms.1}
\sup_{\omega\in\hE_r} \max_{0\leq m\leq\lfloor rT\rfloor} d
(\fx
^{r,m}(\cdot)(\omega),\FMS_1)
\tend0 \qquad\mbox{as $r\tend\infty$}\vadjust{\goodbreak}
\end{equation}
and
%
%
\begin{equation}\label{eq:ht.fms.2}
\sup_{\omega\in\hE_r} d(\fx^{r,0}(\cdot)(\omega),\FMS
(\mathbf{q}_0))
\tend0 \qquad\mbox{as $r\tend\infty$},
\end{equation}
where $\mathbf{q}_0=\hbq_0/(|\hbq_0|\mmax1)$.
\end{lemma}
\begin{pf}
The proof of each equation will use Theorem \ref{thm:fluid}. We start
with (\ref{eq:ht.fms.1}). The theorem requires the use of an index $j$
in some totally ordered countable set; here we shall use the pair
$j\equiv(r,m)$ ordered lexicographically, where $r\in\Naturals$ and
$0\leq
m\leq\lfloor rT\rfloor$. Lexicographic ordering means $(r,m)\geq(r',m')$
iff either
$r>r'$ or both $r=r'$ and $m\geq m'$. Note that $j\tend\infty$ implies
$r\tend\infty$ (and vice versa).

To apply the theorem, we first need to pick constants.
Let $K=1$, let $\blambda$ and $\blambda^j=\blambda^r\tend\blambda$
as per Assumption \ref{cond:mssc.stoch} and let\vadjust{\goodbreak} $\varepsilon_j=\eta
_r(1+1/\log r)$ so
that $\varepsilon_j\tend0$. Thus condition (\ref{cond.rateconv}) of Theorem
\ref{thm:fluid} is satisfied. Now let
\[
G_j
\equiv
G_{r,m}
=
\{(\fx^{r,m}(\cdot)(\omega),z_{r,m}(\omega)) \dvtx
\omega\in\hat{E}_r\}.
\]
It is worth stressing that $G_j\equiv G_{r,m}$ is a set of sample
paths and associated scaling parameters, not a probabilistic event,
and so any questions about the lack of independence between
$\fx^{r,m}(\cdot)(\omega)$ and $z_{r,m}(\omega)$ are void.
Note also that although the events $\hE_r$ lie in different
probability spaces
for each $r$, this has no bearing on the definition of $G_j$ nor on
the application of Theorem~\ref{thm:fluid}.

We next show that $G_j$ satisfies conditions (\ref
{eq:fluid.scale})--(\ref{eq:fluid.qbound}) of
Theorem \ref{thm:fluid}, for~$j$ sufficiently large.
Equation (\ref{eq:fluid.scale}) follows straightforwardly from the
fact that $z_{r,m}\geq r$, hence $\inf\{z\dvtx(\fx,z)\in G_{r,m}\}\geq r$,
hence $\inf\{z\dvtx(\fx,z)\in G_j\}\tend\infty$.
For~(\ref{eq:fluid.arr}), later in the proof we will
establish that, under $\hE_r$ for $r$ large enough,
%
%
\begin{equation}
\label{eq:ht.fms.arr}
\quad{\sup_{t \in[0,\Tfluid]}}
|\fba^{r,m}(t) - \blambda^r t|
<
\eta_r\biggl(1+\frac{1}{\log r}\biggr)\qquad
\mbox{for all } 0\leq m\leq\lfloor rT \rfloor,
\end{equation}
which implies that for all $(\fx,z)\in G_j\equiv G_{r,m}$,
${\sup_{t\in[0,\Tfluid]}} |\fba^j(t)-\blambda^j t|<\varepsilon_j$
as required.
Equation (\ref{eq:fluid.qbound}) follows straightforwardly\vspace*{1pt} from the
scaling used to define $\fbq^j(0)\equiv\fbq^{r,m}(0)$: for every
$\omega$, not merely $\omega\in\hE_r$,
\[
| \fbq^{r,m}(0)|
=
\biggl| \frac{\bQ^r(rm)}{z_{r,m}} \biggr|
=
\biggl| \frac{\bQ^r(rm)}{|\bQ^r(rm)|\mmax r} \biggr|
\leq1.
\]

Since $G_j$ satisfies the conditions of Theorem \ref{thm:fluid}
for sufficiently large $j$, we can apply that theorem to deduce
\[
\sup_{(\fx,z)\in G_j} d(\fx,\FMS_1) \tend0 \qquad\mbox{as $j\tend
\infty$}.
\]
Rewriting $j$ as $(r,m)$, and turning the limit statement into a
$\limsup$ statement,
\[
\sup_{(r',m)\geq(r,0)} \sup_{(\fx,z)\in G_{r',m}}
d(\fx,\FMS_1)
\tend0 \qquad\mbox{as $r\tend\infty$}
\]
and in particular
\[
\max_{0\leq m\leq\lfloor rT\rfloor} \sup_{(\fx,z)\in
G_{r,m}}
d(\fx,\FMS_1) \tend0
\qquad\mbox{as $r\tend\infty$.}
\]
Rewriting $(\fx,z)\in G_{r,m}$ in terms of $\omega\in\hE_r$, as per
the definition of $G_{r,m}$,
\[
\max_{0\leq m\leq\lfloor rT\rfloor} \sup_{\omega\in\hE_r}
d(\fx^{r,m}(\cdot)(\omega),\FMS_1) \tend0
\qquad\mbox{as $r\tend\infty$.}
\]
Interchanging the $\max$ and the $\sup$ gives (\ref{eq:ht.fms.1}).

To establish (\ref{eq:ht.fms.2}), we will again apply Theorem
\ref{thm:fluid} but this time using the index $j\equiv r$,
$\blambda^j=\blambda^r\tend\blambda$ and
$\varepsilon_j=\eta_r(1+1/\log r)$ as above, and define
\[
G_j
\equiv
G_r
=
\{(\fx^{r,0}(\cdot)(\omega),z_{r,0}(\omega))\dvtx
\omega\in\hE_r\}.\vadjust{\goodbreak}
\]
Equations (\ref{cond.rateconv})--(\ref{eq:fluid.qbound}) hold just as
before. For (\ref{eq:fluid.initq}), we will use $\mathbf
{q}_0$ as in
the statement of this lemma, and
$\varepsilon'_j=|\fbq^{r,0}(0)-\mathbf{q}_0|$. This is a
well-defined constant (i.e., it does not depend on the randomness
$\omega$), because we assumed in Theorem \ref{thm:heavytraffic}
that the initial queue sizes $\bQ^r(0)$ are nonrandom, and by definition
$\fbq^{r,0}(0)=\bQ^r(0)/(|\bQ^r(0)|\mmax r)$. Furthermore,\vspace*{1pt} Theorem
\ref{thm:heavytraffic}
assumes $\hbq^r(0)\tend\hbq_0$, which implies
$\fbq^{r,0}(0)\tend\mathbf{q}_0$ hence $\varepsilon'_j\tend0$.
Equation (\ref{eq:fluid.initq}) then follows straightforwardly, for
every $\omega$ not merely $\omega\in\hE_r$. Applying Theorem~\ref{thm:fluid}, we deduce that
\[
\sup_{(\fx,z)\in G_j} d(\fx,\FMS(\mathbf{q}_0))\tend0
\qquad\mbox{as $j\tend\infty$}.
\]
Equivalently,
\[
\sup_{\omega\in\hE_r}
d(\fx^{r,0}(\cdot)(\omega),\FMS(\mathbf{q}_0))\tend0
\qquad\mbox{as $r\tend\infty$}
\]
as required.

To complete the proof of Lemma \ref{lem:ht.fms}, the only remaining
claim that
needs to be established is (\ref{eq:ht.fms.arr}). We will proceed in two
steps. First we prove that $|\bQ^r(rm)|\leq Lr^2$ under $\hE_r$, for
$r$ sufficiently
large and for all $0\leq m\leq\lfloor rT\rfloor$.
To see this, note from (\ref{eq:multihop.discrete.queue.bound}) that
\[
\bQ^r(rm) \leq\bQ^r(0) + \bA^r(rm) + N r m S^{\max}.
\]
Now $\hE_r$ gives a suitable bound on arrivals: for all $0\leq
m'\leq\lfloor rT\rfloor$, and using the fact that $1\leq\Tfluid$,
\[
\biggl| \frac{\bA^r(rm'+r)-\bA^r(rm')}{r} - \blambda^r \biggr|
=
| \fba^{r,m',w_{r,0}}(1) - \blambda^r| < \eta_r,
\]
and by applying this from $m'=0$ to $m'=m-1$
we find $|\bA^r(rm)|\leq
rm(\lambda^{\max}+\eta_r)$. The assumptions of Theorem
\ref{thm:heavytraffic} tell us that $\bQ^r(0)/r\tend\hbq_0$ for some
$\hbq_0\in\RealsP^N$. Putting all this together, we find that
for sufficiently large $r$
\[
|\bQ^r(rm)| \leq r^2\bigl(1+T(1+\lambda^{\max}+N
S^{\max})\bigr)
= Lr^2 \qquad\mbox{for any } 0\leq m\leq\lfloor rT\rfloor.
\]
Now we proceed to prove (\ref{eq:ht.fms.arr}), under $\hE_r$ for $r$
sufficiently large.
Observe that (for $r>2$) there exists $k\in\{1,\ldots,\lfloor Lr\log
r\rfloor\}$ such that
$w_{r,k-1}\leq z_{r,m}\leq w_{r,k}$; this follows from $r\leq
z_{r,m}=|\bQ^r(rm)|\mmax r\leq Lr^2$ and the definition of $w_{r,k}$
in~(\ref{eq:ht.good}).
Hence for any $t\in[0,\Tfluid]$,
\begin{eqnarray*}
&&|\fba^{r,m}(t)-\blambda^rt|\\[-2pt]
&&\qquad=
\biggl|\frac{\bA^r(rm+z_{r,m}t)-\bA^r(rm)}{z_{r,m}} -
\blambda^rt\biggr|\\[-2pt]
&&\qquad=
\biggl|\frac{\bA^r(rm+w_{r,k}u)-\bA^r(rm)}{w_{r,k}} - \blambda
^ru\biggr|
\biggl(\frac{w_{r,k}}{z_{r,m}}\biggr)
\qquad\mbox{where $u=t z_{r,m}/w_{r,k}$}\\[-2pt]
&&\qquad=
| \fba^{r,m,w_{r,k}}(u)-\blambda^r u| \biggl(\frac
{w_{r,k}}{z_{r,m}}\biggr)\\[-2pt]
&&\qquad<
\eta_r \frac{w_{r,k}}{z_{r,m}} \qquad\mbox{since $\hE_r$
holds and $u\leq t\leq\Tfluid$}\\[-2pt]
&&\qquad\leq
\eta_r \frac{w_{r,k}}{w_{r,k-1}} \qquad\mbox{since $z_{r,m}\geq
w_{r,k-1}$}\\[-2pt]
&&\qquad=
\eta_r \biggl(1+\frac{1}{k-1+\log r}\biggr)
\leq
\eta_r\biggl(1+\frac{1}{\log r}\biggr).
\end{eqnarray*}
This establishes (\ref{eq:ht.fms.arr}) and completes the proof.
\end{pf}
\begin{lemma}[(Choice of approximating piece)]
\label{lem:ht.m}
Given $t\in[0,T]$ and $r\in\Naturals$,
define $m^*=m^*(r,t)$ and $u^*=u^*(r,t)$ by
\[
m^* = \min\{ m\in\IntegersP\dvtx rm \leq r^2t \leq rm + \Tfluid
z_{r,m} \} ,\qquad
u^* = \frac{r^2t - r m^*}{z_{r,m^*}}.
\]
This is a sound definition (i.e., the set for $m^*$ is nonempty).
Further, under event~$\hE_r$, either $m^*=0$
and $0\leq u^*\leq\Tfluid$,
or $0<m^*\leq\lfloor rT\rfloor$ and $H_\zeta<u^*\leq\Tfluid$.
\end{lemma}
%
%
\begin{pf}
The set for $m^*$ is nonempty because $z_{r,m}\geq r$
and $\Tfluid\geq1$. The upper bound for $m^*$ is trivial.
The upper bound for $u^*$ in either case is trivial.
To prove the lower bound for $u^*$ when $m^*>0$,
$r^2t>r(m^*-1)+\Tfluid z_{r,m^*-1}$ due to the minimality
of $m^*$. Hence
\[
u^* = \frac{r^2t-rm^*}{z_{r,m^*}}
>
\frac{\Tfluid z_{r,m^*-1} - r}{z_{r,m^*}}
\geq
\Tfluid\frac{z_{r,m^*-1}}{z_{r,m^*}} - 1.
\]
To bound $z_{r,m^*-1}/z_{r,m^*}$, we can use
(\ref{eq:multihop.discrete.queue.bound}) and the bound on
$\mathbf{a}^{r,m,w_{r,0}}(1)$ provided by $\hE_r$ to show that for
any $m$
\begin{eqnarray*}
z_{r,m} &=& |\bQ^r(rm)|\mmax r\\[-2pt]
&\leq&
\bigl( |\bQ^r(rm-r)| +
r(\lambda^{\max}+\eta_r+N S^{\max})
\bigr)\mmax r\\[-2pt]
&\leq&
|\bQ^r(rm-r)|\mmax r + r(\lambda^{\max}+1+N
S^{\max})
\qquad\mbox{since $\eta_r\leq1$}\\[-2pt]
&\leq&
z_{r,m-1} (2+\lambda^{\max}+N S^{\max})
\qquad\mbox{since $z_{r,m-1}\geq r$}.
\end{eqnarray*}
Substituting this back into the earlier bound for $u^*$,
\[
u^*
>
\frac{\Tfluid}{2+\lambda^{\max}+N S^{\max}} - 1,
\]
and this is equal to $H_\zeta$ by choice of $\Tfluid$.
\end{pf}
\begin{lemma}[(Pathwise multiplicative state space collapse)]
\label{lem:ssc.pathwise}
Let \mbox{$0<\zeta<1$}, $t\in[0,T]$ and $r\in\Naturals$ be given.
Suppose there exist $m\in\{0,\ldots,\lfloor rT\rfloor\}$,
$u\in[0,\Tfluid]$ and $x\in\FMS$ such that $r^2t=rm+z_{r,m}u$ and
$\|\fx^{r,m}-x\|<\zeta$, and
furthermore\vadjust{\goodbreak}
either \textup{(i)} $m>0$ and $u>H_\zeta$ and $x\in\FMS_1$, or \textup{(ii)} $m=0$ and
$x\in\FMS(\mathbf{q}_0)$ where $\mathbf{q}_0$ is as defined in
Lemma \ref{lem:ht.fms}.
Then
%
%
\begin{equation}
\label{eq:ht.pathwise}
\quad\frac{ | \hbq^r(t) - \LIFT\WORK(\hbq^r(t)) | }
{\|\hbq^r(\cdot)\|\mmax1}
\leq
\frac{ | \hbq^r(t) - \LIFT\WORK(\hbq^r(t)) | }{z_{r,m}/r}
<
2\zeta+
\modcont_\zeta(\LIFT\WORK),
\end{equation}
where $\modcont_\zeta(\LIFT\WORK)$ is the modulus of continuity of
the map
$\mathbf{q}\mapsto\LIFT\WORK(\mathbf{q})$ over
\[
\cD= \{ \mathbf{q}'\in\RealsP^N \dvtx |\mathbf{q}'-\mathbf
{q}|\leq1 \mbox{ for some $\mathbf{q}$
such that } L(\mathbf{q})\leq L(\bOne) \}
\]
for $L(\cdot)$ as in Definition \ref{def:lyapunov}.
\end{lemma}
\begin{pf}
The first inequality is trivially true because
\[
\frac{z_{r,m}}{r} = \frac{|\bQ^r(rm)|\mmax r}{r}
\leq
\Bigl(\sup_{u\in[0,T]} \hbq
^r(u)\Bigr) \mmax1.
\]
For the second inequality, note that after unwrapping the
$\hbq^r(\cdot)$ scaling and wrapping it up again in the $\fbq^{r,m}$ scaling,
the middle term in (\ref{eq:ht.pathwise}) is
\[
\mbox{MT}=
| \fbq^{r,m}(u) - \LIFT\WORK(\fbq^{r,m}(u)) |.
\]
Writing
$\mathbf{q}$ for the queue component of $x$,
\begin{eqnarray*}
\mbox{MT}
&\leq&
| \fbq^{r,m}(u) - \mathbf{q}(u) |
+
| \mathbf{q}(u)-\LIFT\WORK(\mathbf{q}(u)) |
+
| \LIFT\WORK(\mathbf{q}(u)) - \LIFT\WORK(\fbq^{r,m}(u))
|\\
&=&
\mbox{(\ref{equ52a})}
+
\mbox{(\ref{equ52b})}
+
\mbox{(\ref{equ52c})}\qquad\mbox{respectively}.
\end{eqnarray*}
We can bound each term as follows:
%
%
\begin{subequation}
\begin{eqnarray}
\label{equ52a}
\begin{tabular}{p{317pt}}
is $<\zeta$ since
$\|\fx^{r,m}-x\|<\zeta$ by an assumption of the lemma.
\end{tabular}
\hspace*{-42pt}
\\
\label{equ52b}
\\[-18pt]
\begin{tabular}{p{317pt}}
is $<\zeta$ in the case $m>0$: by the assumptions of the lemma,
$x\in\FMS_1$ so $|\mathbf{q}(0)|\leq1$, and also $u>H_\zeta$.
The requirements of Theorem~\ref{thm:fp}(v)
are met; hence we obtain the inequality.
\end{tabular}\hspace*{-42pt}\nonumber\\
\begin{tabular}{p{317pt}}
is $=0$ in the case $m^*=0$:
In this case, by assumption of the lemma $x\in\FMS(\mathbf{q}_0)$,
that is, $\mathbf{q}(0)=\mathbf{q}_0$. By assumption of Theorem \ref
{thm:heavytraffic}
$\hbq_0\in\INV$, that is, $\hbq_0=\LIFT\WORK(\hbq_0)$,
therefore by Theorem \ref{thm:fp}(iii) $\mathbf{q}_0\in
\INV$,
therefore by Theorem \ref{thm:fp}(iv)
the fluid model solution $\mathbf{q}(\cdot)$ stays constant at
$\mathbf{q}_0$ and so $\mbox{(\ref{equ52b})}=0$.
\end{tabular}\nonumber
\hspace*{-42pt}
\\
\label{equ52c}
\\[-18pt]
\begin{tabular}{p{317pt}}
is $\leq\modcont_\zeta(\LIFT\WORK)$: by assumption of the
lemma, either $m>0$ and $x\in\FMS_1$, or $m=0$ and
$x\in\FMS(\mathbf{q}_0)$ where $\mathbf{q}_0\leq\bOne$
componentwise; either
way $\mathbf{q}(0)\leq\bOne$ componentwise, so $L(\mathbf
{q}(0))\leq L(\bOne)$.
By Theorem \ref{thm:fp}(i) $L(\mathbf{q}(u))\leq
L(\bOne)$
so $\mathbf{q}(u)\in\cD$. Furthermore, since $\|\fx^{r,m}-x\|<\zeta$
by assumption of the lemma, $|\fbq^{r,m}(u)-\mathbf{q}(u)|<\zeta<1$
and so
$\fbq^{r,m}(u)\in\cD$. The inequality then follows from the
definition of the modulus of continuity.
\end{tabular}\nonumber
\hspace*{-42pt}
\end{eqnarray}
\end{subequation}
\upqed\end{pf}
\begin{lemma}[(The good event has high probability)]
\label{lem:ht.prob}
Under\vspace*{1pt} the assumptions of Theorem \ref{thm:heavytraffic},
$\Prob(\hE_r)\tend1$ as $r\tend\infty$.
The deviation terms are given by
$\eta_r=\min(1,\sup_{z\geq r}\Tfluid\delta_{\lfloor z
\Tfluid\rfloor})$ and
$\eta_r\tend0$ as $r\tend\infty$.\vadjust{\goodbreak}
\end{lemma}
\begin{pf}
By a simple union bound, and then using the fact that the arrival
process has stationary increments,
\begin{eqnarray*}
\Prob(\hE_r)
&\geq&
1 - \sum_{m=0}^{\lfloor rT\rfloor} \sum_{k=0}^{\lfloor Lr\log
r\rfloor}
\Prob\Bigl( {\sup_{u\in[0,\Tfluid]}} |
\fba^{r,m,w_{r,k}}(u) - \blambda^r u| \geq\eta_r\Bigr)\\[-2pt]
&=&
1 - (1+\lfloor rT\rfloor) \sum_{k=0}^{\lfloor Lr\log r\rfloor}
\Prob\Bigl( {\sup_{u\in[0,\Tfluid]}} |
\fba^{r,0,w_{r,k}}(u) - \blambda^r u| \geq\eta_r\Bigr)\\[-2pt]
&=&
1 - (1+rT) \sum_{k=0}^{\lfloor Lr\log r\rfloor}
\Prob\biggl( \sup_{u\in[0,\Tfluid]} \biggl|
\frac{A^r(w_{r,k}u)}{w_{r,k}} - \blambda^r u\biggr| \geq\eta
_r\biggr).
\end{eqnarray*}
To bound this we will use Assumption \ref{cond:mssc.stoch}, which says that
\[
z (\log z)^2 \Prob \biggl(\sup_{\tau\leq z} \frac{1}{z}|
\bA^r(\tau)-\blambda^r\tau| \geq\delta_z \biggr)
\tend0
\qquad\mbox{as $z\tend\infty$}
\]
uniformly in $r$.
After extending the domain of $\bA^r$ to $\RealsP$
by linear interpolation in each interval $(\tau,\tau+1)$, and
extending the domain of $\delta_z$ to $z\in\RealsP$ by
$\delta(z)=\delta_{\lfloor z\rfloor}\mmax\delta_{\lceil z\rceil}$,
and rescaling $z$ by $\Tfluid$,
\[
z(\log z)^2 \Prob\biggl( \sup_{u\in[0,\Tfluid]}
\biggl| \frac{\bA^r(zu)}{z} - \blambda^ru\biggr| \geq\Tfluid
\delta(z \Tfluid) \biggr)
\tend0 \qquad\mbox{as $z\tend\infty$}
\]
uniformly
in $r$. In other words, for any $\phi>0$ there exists
$z_0$ such that for all $z\geq z_0$ and all $r$,
\[
\Prob\biggl( \sup_{u\in[0,\Tfluid]}
\biggl| \frac{\bA^r(zu)}{z} - \blambda^ru\biggr| \geq\Tfluid
\delta(z \Tfluid) \biggr) < \frac{\phi}{z(\log z)^2}.
\]
Now pick $r_0$ sufficiently large that $r_0\geq z_0$ and
$\sup_{z\geq r_0}\Tfluid\delta(z\Tfluid)<1$, which we can do since
$\delta(z)\tend0$ as $z\tend\infty$ by Assumption \ref{cond:mssc.stoch}.
This choice implies that for any $r\geq r_0$ and $z\geq r$,
$\Tfluid\delta(z\Tfluid)\leq\eta_r$ [recall that
$\eta_r =\min(1,\sup_{z\geq r}\Tfluid\delta_{\lfloor z
\Tfluid\rfloor})$].
Hence, for any $r\geq r_0$ and $z\geq r$,
\[
\Prob\biggl( \sup_{u\in[0,\Tfluid]}
\biggl| \frac{\bA^r(zu)}{z} - \blambda^ru\biggr| \geq\eta_r
\biggr) < \frac{\phi}{z(\log z)^2}.
\]
Applying this bound to $\Prob(\hE_r)$, and using the facts that
$w_{r,k}\geq r$ and\break
$w_{r,k}(\log w_{r,k})^2\geq r(1+k/\log
r)(\log r)^2$,
\begin{eqnarray*}
1-\Prob(\hE_r)
&<&
\phi\biggl(\frac{1+\lfloor rT\rfloor}{r(\log r)^2}\biggr) \sum
_{k=0}^{\lfloor Lr\log r\rfloor}
\frac{1}{1+k/\log r}\\[-2pt]
&=&
\phi\biggl(\frac{1+\lfloor rT\rfloor}{r\log r}\biggr) \sum
_{k=0}^{\lfloor Lr\log r\rfloor}\frac{1}{k+\log r}\\[-2pt]
&\leq&
\phi\biggl(\frac{1+\lfloor rT\rfloor}{r\log r}\biggr)
\int_{\ell=0}^{\lfloor Lr\log r\rfloor} \frac{1}{\ell-1+\log r}\,
d\ell
\\
&=& \phi\biggl( \frac{1+\lfloor rT\rfloor}{r\log r}\biggr)
\log\biggl(1+\frac{\lfloor Lr\log r\rfloor}{\log r-1}\biggr).
\end{eqnarray*}
The final expression converges to $\phi T$ as $r\tend\infty$. Since
$\phi$ can
be chosen arbitrarily small, $\Prob(\hE_r)\tend1$ as $r\tend\infty$.
\end{pf}

\subsubsection{\texorpdfstring{Proof of Theorem \protect\ref
{thm:heavytraffic}.}{Proof of Theorem 7.1}}
\label{sec:ht.proof}

Given $\delta>0$, pick $\zeta>0$ such that
$2\zeta+\modcont_\zeta(\LIFT\WORK) < \delta\mmin1$
where $\modcont_\zeta(\LIFT\WORK)$ is the modulus of continuity
of $\LIFT\WORK$ over the set $\cD$ specified in Lemma
\ref{lem:ssc.pathwise}. We can achieve the desired bound
by making $\zeta$ sufficiently small; this is because
$\LIFT\WORK$ is continuous, hence uniformly continuous on compact
sets, and $\cD$ is compact as a consequence of
Theorem \ref{thm:fp}(i), hence
$\modcont_\zeta(\LIFT\WORK)\tend0$ as $\zeta\tend0$.
With this choice of $\zeta$, define the good sets $\hE_r$ and
the constants $\Tfluid$ and $H_\zeta$ as specified by
(\ref{eq:ht.good}).

By Lemma \ref{lem:ht.fms}, there exists $r_0$ such that for $r\geq
r_0$ and for all $\omega\in\hE_r$ and all $m$,
$d(\fx^{r,m},\FMS_1) < \zeta$ and
$d(\fx^{r,0},\FMS(\mathbf{q}_0)) < \zeta$, where $\mathbf{q}_0$ is
defined in
the statement of the lemma.

Now, pick any $t\in[0,\Tfluid]$ and $r\geq r_0$, and assume $\hE_r$
holds. Lemma~\ref{lem:ht.m} says that
we can choose $m\in\{0,\ldots,\lfloor rT\rfloor\}$ and $u\in
[0,\Tfluid]$ such
that $r^2t=rm+z_{r,m}u$, and furthermore either (i) $m>0$ and
$u>H_\zeta$ or (ii) $m=0$.
By Lem\-ma~\ref{lem:ht.fms}, we can pick $x\in\FMS$ (depending on $r$,
$t$ and the $\omega$) such that $\|\fx^{r,m}-x\|<\zeta$ and
furthermore either
(i) $m>0$ and $x\in\FMS_1$ or (ii) $m=0$ and $x\in\FMS(\mathbf{q}_0)$.
Then, by Lemma \ref{lem:ssc.pathwise},
\[
\frac{ | \hbq^r(t) - \LIFT\WORK(\hbq^r(t)) | }
{\|\hbq^r(\cdot)\|\mmax1}
<
\delta.
\]
This bound holds for every $t\in[0,T]$ and $r\geq r_0$, in a
sample path-wise sense, whenever $\omega\in\hE_r$.

Finally, Lemma \ref{lem:ht.prob} says that $\Prob(\hE_r)\tend1$. This
completes the proof.
\hfill$\Box$

\section{\texorpdfstring{An optimal policy?}{An optimal policy}}
\label{sec:optimal}

Our motivation for this work was Conjecture \ref{conj:keslassy}, which
says that for an input-queued switch the performance of MW-$\alpha$
improves as $\tendsdown$. We have not been able to prove this.
However, under a~condition on the arrival rate, we can show (i) that
the critically-loaded fluid model solutions for a single-hop switched
network approach optimal (in the sense of minimizing
total amount of work in the network) as $\tendsdown$; and (ii)
that for an input-queued switch the set of invariant states $\INV$
defined in
Section \ref{sec:ht} is sensitive to $\alpha$. We speculate that
these findings might eventually form part of a proof of a heavy
traffic limit theorem supporting Conjecture~\ref{conj:keslassy}, given
that critically loaded fluid models and invariant states play an
important role in heavy traffic theorems.

In this section we state the condition on the arrival rates,
and give the results (i) and (ii). Motivated by these results, we make
a conjecture about an optimal scheduling policy.\vadjust{\goodbreak}
\begin{definition}[(Complete loading)]
Consider a switched network with arrival rate vector
$\blambda$. Say that $\blambda$ satisfies the
\textit{complete loading condition} if
$\blambda\in\Lambda$, and
there is a convex combination of
critically loaded virtual resources that gives equal weight to each
queue; in other words if
\[
\frac{\bOne}{\max_{\bpi\in\sS} \bOne\bdot\bpi} \in\langle
{\CLVR(\blambda)}\rangle.
\]
\end{definition}
\begin{theorem}[(Near-optimality of fluid models under complete
loading)]\label{thm:opt.dynamic}
Consider a single-hop switched network with arrival rate vector
$\blambda\in\Lambda$.
\begin{longlist}
\item
For any fluid model solution for the MW-$\alpha$ policy,
$\bOne\bdot\mathbf{q}(t)\leq\break N^{\alpha/(1+\alpha)}\bOne\bdot
\mathbf{q}(0)$.
\item
For any fluid model solution for any scheduling policy,
if $\blambda$ satisfies the complete loading condition, then
$\bOne\bdot\mathbf{q}(t)\geq\bOne\bdot\mathbf{q}(0)$.
\end{longlist}
\end{theorem}
\begin{pf}
The first claim relies on the
standard result that for any $\bx\in\RealsP^\N$ and
$\beta>1$,
%
%
\begin{equation}
\label{eq:holder2}
\frac{1}{\N^{1-1/\beta}} \sum_{n}x_{n}
\leq
\biggl( \sum_{n}x_{n}^\beta\biggr)^{1/\beta}
\leq
\sum_{n}x_{n}.
\end{equation}
Using the Lyapunov function from Definition \ref{def:lyapunov},
\begin{eqnarray*}
\bOne\bdot\mathbf{q}(t)
&\leq&
\N^{1-1/(1+\alpha)} \biggl( \sum_{n}q_{n}(t)^{1+\alpha}
\biggr)^{1/(1+\alpha)}\qquad
\mbox{by the first inequality in (\ref{eq:holder2})}\\
&=&
\N^{\alpha/(1+\alpha)} L(\mathbf{q}(t))^{1/(1+\alpha)}
\qquad\mbox{by definition of $L(\cdot)$}\\
&\leq&
\N^{\alpha/(1+\alpha)} L(\mathbf{q}(0))^{1/(1+\alpha)}
\qquad\mbox{since $\dot{L}(\mathbf{q}(t))\leq0$ by
Theorem \ref{thm:fp}(i)}\\
&\leq&
N^{\alpha/(1+\alpha)} \bOne\bdot\mathbf{q}(0)\qquad
\mbox{by the second inequality in (\ref{eq:holder2}).}
\end{eqnarray*}
The second claim is a simple
consequence of Lemma \ref{lem:feasible}.
(This lemma is for a single-hop network. The multi-hop version, Lemma
\ref{lem:multihop.feasible}, does not have such a simple interpretation.)
\end{pf}
\begin{theorem}[($\INV$ is sensitive to $\alpha$ for an input-queued switch)]
\label{thm:mwma}
Consider an $M\times M$ input-queued switch running MW-$\alpha$,
as introduced in Section \ref{ssec:iq}.
Let $\lambda_{i j}$ be the arrival rate at the queue at input port $i$ of
packets destined for output port $j$, $\blambda\in\RealsP^{M\times
M}$.
Suppose that $\blambda>\bZero$
componentwise, and furthermore that every input port and every output
port is critically loaded,\vadjust{\goodbreak} that is,
%
%
\begin{equation}
\label{eq:iq.completeload}
\sum_{j=1}^M \lambda_{\hi j}=1
\quad\mbox{and}\quad
\sum_{i=1}^M \lambda_{i \hj}=1
\qquad\mbox{for every }
1\leq\hi,\hj\leq M.
\end{equation}
Then $\blambda$ satisfies the complete
loading condition, and the critically loaded virtual resources are
\[
\CLVR(\blambda) = \{\br_\hi\mbox{ for all }1\leq\hi\leq
M\}
\cup
\{\bc_\hj\mbox{ for all }1\leq\hj\leq M\},
\]
where\vspace*{-1pt} $\br_\hi$ and $\bc_\hj$ are the row and column indicator matrices,
$(\br_\hi)_{i,j}=1_{i=\hi}$ and $(\bc_\hj)_{i,j}=1_{j=\hj}$.
Define the workload map $W\dvtx\RealsP^{M\times M}\to\RealsP^{2M}$ by
$W(\mathbf{q})=[\bxi\bdot\mathbf{q}]_{\bxi\in\CLVR(\blambda)}$.
Denoting the
invariant set by $\INV(\alpha)$:
\begin{longlist}
\item
if $w$ is in the relative interior of $\{W(\mathbf{q})\dvtx\mathbf{q}\in
\RealsP^{M\times
M}\}$, then $w$ is in
$\{W(\mathbf{q})\dvtx\mathbf{q}\in\INV(\alpha)\}$ for sufficiently
small $\alpha>0$;
\item
for a $2\times2$ input-queued switch, $\{W(\mathbf{q})\dvtx\mathbf{q}\in
\INV(\alpha)\}$ is strictly
increasing as $\tendsdown$.
\end{longlist}
\end{theorem}

Item (i) essentially says that $W(\INV(\alpha))$ becomes
as large as possible as $\tendsdown$, except for some possible
issues at the boundary. We have only been able to prove (ii)
for a $2\times2$ switch, but we believe it holds for any $M\times M$ switch.
The proofs are rather long, and depend on the specific structure of
the input-queued switch, so they are left to the \hyperref[app]{Appendix}.

Conjecture \ref{conj:keslassy} claims that, for an input-queued
switch, performance improves as $\tendsdown$. Examples due to
Ji, Athanasopoulou and Srikant \cite{srikantsmallswitch}
and Stolyar (personal communication) show
that this is not true for general switched networks. However, Theorem
\ref{thm:opt.dynamic} suggests that the conjecture might apply
not just to input-queued switches but also to generalized switches
under the complete loading condition; the examples of
Ji, Athanasopoulou and Srikant \cite{srikantsmallswitch} and Stolyar
do not satisfy this
condition. We therefore extend Conjecture \ref{conj:keslassy} as follows.
\begin{conjecture}
Consider a general single-hop switched network as described in Section
\ref{sec:model}, running MW-$\alpha$.
Consider the diffusion scaling limit (described in Conjecture
\ref{conj:heavytraffic.ssc}), and let $\blambda$ be the limiting
arrival rates; assume~$\blambda$ satisfies the complete loading condition.
For every $\alpha>0$ there is a~limiting stationary queue size
distribution. The expected value of the sum
of queue sizes under this distribution is nonincreasing as
$\tendsdown$.
\end{conjecture}

Theorem \ref{thm:opt.dynamic} says that MW-$\alpha$ approaches optimal
as $\tendsdown$, under the complete loading condition. It is
natural to ask if there is a policy that is optimal, rather than just a
sequence of policies that approach optimal.
Given that MW-$\alpha$ chooses\vadjust{\goodbreak} a schedule $\bpi$ to maximize
$\bpi\bdot\mathbf{q}^\alpha$ (where the exponent is taken
componentwise), and since
\[
x^\alpha= \cases{
1+\alpha\log x + O(\alpha^2), &\quad if $x>0$,\cr
0, &\quad if $x=0$,}
\]
we
propose the following formal limit policy, which we call MSMW-log:
\textit{at each time step}, \textit{look at all maximum-size schedules},
\textit{that is}, \textit{those
$\bpi\in\sS$ for which $\sum_n\pi_n 1_{Q_n>0}$ is maximal. Among
these}, \textit{pick one which has maximal log-weight}, \textit{that is}, \textit{for which
$\sum_{n\dvtx Q_n>0}\pi_n\log Q_n$ is maximal}, \textit{breaking ties randomly.}
\begin{conjecture}
\label{conj:sw}
Consider a general single-hop switched network running MSMW-log.
Consider the diffusion scaling limit, and let $\blambda$ be the
limiting arrival rates; assume $\blambda$ satisfies the complete
loading condition. There is
a limiting stationary queue size
distribution. This distribution minimizes the expected value of the
sum of the queue sizes, over all scheduling policies for which this
expected value is defined.
\end{conjecture}

Scheduling policies based on MW are computationally difficult to
implement because there are so many comparisons to be made. In future
work we plan to investigate whether the techniques described in this
paper can be applied to policies that may have worse
performance but simpler implementation.

\begin{appendix}\label{app}

\section*{Appendix: Results for input-queued switches}

In this section we prove Theorem \ref{thm:mwma}. Throughout this
Appendix we are considering a $M\times M$ input-queued switch.
The set of schedules $\sS$ consists of all $M\times M$ permutation
matrices. We assume the arrival rate matrix $\blambda$ satisfies the
complete loading condition (\ref{eq:iq.completeload})
and that $\blambda>\bZero$ componentwise.
We let the
scheduling algorithm be MW-$\alpha$, and define $\INV(\alpha)$ to be
the set of invariant states.

\subsection{\texorpdfstring{Identifying
$\Lambda$, $\sS^*$, $\CLVR$ and $\CLVR^+$.}{Identifying
Lambda, S^*, Xi and Xi^+}}
The Birkhoff--von Neumann decomposition
result says that a matrix is doubly substochastic
if and only if it is less than or
equal to a convex combination of permutation matrices, which yields
\[
\Lambda= \Biggl\{ \blambda\in[0,1]^{M\times M} \dvtx
\sum_{j=1}^M \lambda_{\hi,j}\leq1 \mbox{ and }
\sum_{i=1}^M \lambda_{i,\hj}\leq1 \mbox{ for all $\hi, \hj$}
\Biggr\}.
\]
Since $\blambda$ satisfies the complete loading condition
(\ref{eq:iq.completeload}), $\blambda\in\Lambda$.

Lemma \ref{lem:iq.vr} below gives $\sS^*$, the set of virtual
resources, that is, maximal extreme points
of the set of feasible solutions to $\DUAL(\blambda)$.
From the complete loading condition,
it is clear that $\CLVR(\blambda)=\sS^*$ as claimed in the theorem.
It will also be useful, for the proof of Theorem
\ref{thm:mwma}(i), to identify $\CLVR^+(\blambda)$\vadjust{\goodbreak}
as defined by (\ref{eq:clvrplus}). We claim that
%
%
\begin{equation}
\label{eq:mwma.clvrp}
\CLVR^+(\blambda)=\CLVR(\blambda).
\end{equation}
To see this,
suppose $\bxi$ is a nonmaximal
extreme point of the set of feasible solutions to $\DUAL(\blambda)$;
then there exists some other extreme point
$\bzeta$ such that $\bxi\leq\bzeta$ and $\bxi\neq\bzeta$;
but because $\blambda>\bZero$ componentwise it must be that
$\bxi\bdot\blambda<\bzeta\bdot\blambda$. We have found that
$\blambda\in\Lambda$,
so the solution to $\DUAL(\blambda)$ is $\leq1$, hence
$\bxi\bdot\blambda<1$. Therefore $\bxi\notin\CLVR^+(\blambda)$,
that is, $\CLVR^+(\blambda)$ consists only of maximal extreme points.
\setcounter{assumption}{0}
\begin{lemma}
\label{lem:iq.vr}
The set of maximal extreme points of the set
\[
F = \{ \bxi\in\RealsP^{M\times M} \dvtx \bxi\bdot\bpi\leq1
\mbox{ for all }\bpi\in\sS\}
\]
is given by
\[
\sS^*
=
\{\br_\hi\mbox{ for all $1\leq\hi\leq M$}\}
\union
\{\bc_\hj\mbox{ for all $1\leq\hj\leq M$}\},
\]
where the row and column indicator matrices $\br_\hi$ and $\bc_\hj$
are defined by
$(\br_\hi)_{i,j}=1_{i=\hi}$ and $(\bc_\hj)_{i,j}=1_{j=\hj}$.
\end{lemma}
\begin{pf}
First we argue that every $\bxi\in\sS^*$ is a maximal extreme point of~$F$. It is simple to check that $\bxi\in F$. Also, $\bxi$ is extreme
because $F\subset[0,1]^{M\times M}$. Finally, $\bxi$ is maximal
because if it were not then there would be some
$\bvarepsilon\geq\bZero$, $\bvarepsilon\neq\bZero$, such that
$\bxi+\bvarepsilon\in F$; but for any such $\bvarepsilon$ there is a~matching~$\bpi$ such that $\bvarepsilon\bdot\bpi>0$ hence
$\bxi+\bvarepsilon\notin F$.

Next we argue the converse, that all maximal extreme points
of $F$ are in $\sS^*$. The first step is to characterize the
extreme points of $F$. We claim that if $\bzeta\in F$, then it can be
written
$\bzeta\leq\bxi$ for some $\bxi\in\langle{\sS^*}\rangle$.
Consider the optimization problem
%
%
\begin{eqnarray}
\label{eq:iqd.primal}
&&\mbox{minimize }
\sum_{\hi=1}^Mx_\hi+ \sum_{\hj=1}^My_\hj
\quad\mbox{over}\quad
x_\hi\geq0, y_\hj\geq0 \qquad\mbox{for all }\hi,\hj\nonumber\\[-10pt]\\[-10pt]
&&\mbox{such that }
\bzeta\leq\sum_\hi x_\hi\br_\hi+ \sum_\hj y_\hj\bc_\hj.\nonumber
\end{eqnarray}
The dual of this problem is
%
%
\begin{eqnarray}
\label{eq:iqd.dual}
&&\mbox{maximize }
\mathbf{a}\bdot\bzeta
\quad\mbox{over}\quad
\mathbf{a}\in\RealsP^{M\times M}\nonumber\\[-8pt]\\[-8pt]
&&\mbox{such that }
\mathbf{a}\bdot\br_\hi\leq1,\qquad \mathbf{a}\bdot\bc_\hj\leq1\qquad
\mbox{for all }\hi,\hj.\nonumber
\end{eqnarray}
[These problems are just $\PRIMAL(\bzeta)$ and $\DUAL(\bzeta)$,
resp., but with respect to the ``virtual'' schedule set $\sS^*$
rather than the actual schedule set $\sS$.] By Slater's condition,
strong duality holds. Now, any matrix $\mathbf{a}$ that is feasible for
(\ref{eq:iqd.dual}) is nonnegative and doubly substochastic, hence by
the Birkhoff--von Neumann decomposition result it can be written as a
$\mathbf{a}\leq\mathbf{b}$ where $\mathbf{b}$ is a~convex combination of permutation matrices, that is, $\mathbf{b}\in
\langle{\sS}\rangle$.
But by the assumption\vadjust{\goodbreak} that $\bzeta\in F$, $\bzeta\bdot\bpi\leq1$ for
all $\bpi\in\sS$, hence $\mathbf{b}\bdot\bzeta\leq1$,
hence $\mathbf{a}\bdot\bzeta\leq1$,
so the value of the
optimization problem (\ref{eq:iqd.dual}) is $\leq1$. By strong
duality, the value of the optimization problem (\ref{eq:iqd.primal})
is also $\leq1$. Therefore $\bzeta\leq\bxi$ for some
$\bxi\in\langle{\sS}\rangle$.

We now claim that if $\bzeta\in F$ then it can be written
$\bzeta= \sum_\bxi a_\bxi\bxi$
where the sum is over some finite collection of values
drawn from the set
\[
E =\{\bxi\in\RealsP^{M\times M} \dvtx
\bxi\leq\br_\hi\mbox{ or }\bxi\leq\bc_\hj\mbox{ for some }\hi
, \hj
\},
\]
and all $a_\bxi$ are $\geq0$.
We have just shown that
%
%
\begin{equation}
\label{eq:iq.induction}
\bzeta= \sum_\bxi a_\bxi\bxi- \mathbf{z}
\qquad\mbox{for some }\mathbf{z}\geq\bZero
\mbox{ and }a_\bxi\geq0,\qquad \bxi\in E.
\end{equation}
If $\mathbf{z}=\bZero$ we are
done. Otherwise, pick some $(k,l)$ such that $z_{k,l}>0$ and define~$\bn^{k,l}$ by $(\bn^{k,l})_{i,j}=1-1_{i=k\mbox{ and }j=l}$. Noting
that $\sum a_\bxi\xi_{k,l}\geq z_{k,l}>0$, we can rewrite
$\bzeta$ as
\[
\bzeta=
\frac{z_{k,l}}{\sum a_\bxi\xi_{k,l}} \sum a_\bxi\bxi
+
\biggl( 1- \frac{z_{k,l}}{\sum a_\bxi\xi_{k,l}} \biggr)
\sum a_\bxi\bxi\bn^{k,l}
-
\mathbf{z}\bn^{k,l},
\]
where matrix multiplication is componentwise as per the notation
specified in
Section \ref{sec:intro}. We have thus rewritten $\bzeta$ in the form
(\ref{eq:iq.induction}), but now $\mathbf{z}$ has one fewer nonzero element.
Continuing in this way we can remove all nonzero elements of $\mathbf{z}$,
until we are left with $\bzeta\in\langle{E}\rangle$.

We have therefore shown that all extreme points of $F$ are
in $E$. Clearly, all maximal points of $E$ are in $\sS^*$. Therefore,
all the maximal extreme points of $F$ are in $\sS^*$ as claimed.
\end{pf}

\subsection{\texorpdfstring{Proof of Theorem \protect\ref
{thm:mwma}(i).}{Proof of Theorem 8.3(i)}}

We first state two lemmas which will be needed in the proof. The
first is a general closure property of permutation matrices, and the
second is a property of invariant states of MW-$\alpha$. We then
prove the theorem and the two lemmas.
\begin{lemma}
\label{lem:mwm.closure}
Let $\bx\in\RealsP^{M\times M}$, and define $\bA\in\RealsP
^{M\times
M}$ by $A_{i,j}=1$ if there is some matching $\brho\in\sS$ whose
weight $\brho\bdot\bx$ is maximal
(i.e., $\brho\bdot\bx=\max_{\bsigma\in\sS}\bsigma\bdot\bx$)
and for which $\rho_{i,j}=1$; and $A_{i,j}=0$ otherwise. Then, for any
matching $\bpi$ such that
\[
\pi_{i,j}=1 \quad\implies\quad A_{i,j}=1,
\]
$\bpi$ is itself a maximum weight matching.
\end{lemma}
\begin{lemma}
\label{lem:mwm.sameweight}
Fix any $\blambda\in\Lambda$, and any $\mathbf{q}\in\INV(\alpha)$.
For every $1\leq i,j\leq M$ such that $\lambda_{i,j}>0$, there exists
a matching $\bpi\in\sS$ whose weight $\bpi\bdot\mathbf{q}^\alpha
$ is maximal
(i.e., $\bpi\bdot\mathbf{q}^\alpha=\max_{\bsigma\in\sS}\bsigma
\bdot\mathbf{q}^\alpha$,
where the exponent is taken componentwise)
and for which $\pi_{i,j}=1$.
\end{lemma}

\begin{pf*}{Proof of Theorem \ref{thm:mwma}\textup{(i)}}
Suppose the claim of the theorem is not true, that is, that there
exists a
sequence $\tendsdown$ such that
$w\notin\SSC(\alpha)=\{W(\mathbf{q})\dvtx\mathbf{q}\in\INV(\alpha
)\}$ for each $\alpha$
in the sequence.\vadjust{\goodbreak}

Write $w=(w_{1\cdot},\ldots,w_{M\cdot}, w_{\cdot1},\ldots
,w_{\cdot
M})$,
and define the function $\LIFT^\alpha\dvtx\RealsP^{2M}\to\RealsP
^{M\times
M}$ to give the (unique) solution to the optimization problem
\begin{eqnarray*}
&&\mbox{minimize }
\frac{1}{1+\alpha}\sum_{i,j} q_{i,j}^{1+\alpha}
\quad\mbox{over}\quad
\mathbf{q}\in\RealsP^{M\times M}\\
&&\mbox{such that }
\br_\hi\bdot\mathbf{q}\geq w_{\hi\cdot}
\quad\mbox{and}\quad
\bc_\hj\bdot\mathbf{q}\geq w_{\cdot\hj}
\qquad\mbox{for all }1\leq\hi,\hj\leq M.
\end{eqnarray*}
This is\vspace*{1pt} the optimization problem defined in (\ref{eq:liftwplus}); we have
simply written out $\CLVR^+(\blambda)$ explicitly using
(\ref{eq:mwma.clvrp}).
By Lemma \ref{lem:liftwplus}, the map
$\LIFTW\dvtx\RealsP^{M\times M}\to\RealsP^{M\times M}$
that defines $\INV(\alpha)$
is simply the
composition of $W\dvtx\RealsP^{M\times M}\to\RealsP^{2M}$ and
$\LIFT^\alpha\dvtx\RealsP^{2M}\to\RealsP^{M\times M}$.
Let $\mathbf{q}(\alpha)=\LIFT^\alpha(w)$.
Note that $\mathbf{q}(\alpha)\in\INV(\alpha)$; this is because
$\mathbf{q}(\alpha)$ is optimal for $\LIFT^\alpha(w)$, therefore it is
optimal for $\LIFT^\alpha(W(\mathbf{q}(\alpha)))$ which has a smaller
feasible region, therefore
$\mathbf{q}(\alpha)=\LIFT^\alpha(W(\mathbf{q}(\alpha)))=\break\LIFTW
(\mathbf{q}(\alpha))$, that is, $\mathbf{q}(\alpha)\in\INV(\alpha)$.

We next establish this claim: that \textit{for each $\alpha$ in
the sequence}, \textit{there exists}~$i$, $j$, $i'$ \textit{and} $j'$ \textit{such
that} $q_{i,j}(\alpha)=0$, $q_{i',j}(\alpha)\geq w_{\cdot j}/M$
\textit{and} $q_{i,j'}(\alpha)\geq w_{i\cdot}/M$.
To prove this claim,
observe that $W(\mathbf{q}(\alpha))\geq w$ by the constraints of the
optimization problem $\LIFT^\alpha$; and that
$\mathbf{q}(\alpha)\in\INV(\alpha)$ hence $W(\mathbf{q}(\alpha
))\in\SSC(\alpha)$;
hence $W(\mathbf{q}(\alpha))\neq w$ by assumption that $w\notin
\SSC(\alpha)$.
Therefore $W(\mathbf{q}(\alpha))>w$ in some component,
that is, there is some $i$ or $j$ such that
$\br_i\bdot\mathbf{q}(\alpha)>w_{i\cdot}$ or $\bc_j\bdot\mathbf
{q}(\alpha)>w_{\cdot j}$. Indeed,
there must be both such an $i$ and $j$, since otherwise the sum of row
workloads and column workloads would not be equal. Therefore
$q_{i,j}(\alpha)=0$, since if $q_{i,j}(\alpha)>0$, then we could reduce
$q_{i,j}(\alpha)$ and
still have a feasible solution to the problem that defines $\LIFT
^\alpha(w)$ but
with a smaller value of the objective function, which contradicts
optimality of $\mathbf{q}(\alpha)$. There must also be a $j'$ such that
$q_{i,j'}(\alpha)\geq w_{i\cdot}/M$ since otherwise the workload
constraint $\br_i\bdot\mathbf{q}(\alpha)\geq w_{i\cdot}$ would not be
met. Likewise for $i'$.

We assumed that $w\in\RealsP^{2M}$ is in the relative interior of
$\SSCm= \{W(\mathbf{q})\dvtx\mathbf{q}\in\RealsP^{M\times M}\}$. This
set is clearly convex,
and from the characterization of relative interior for convex sets,
for all $x\in\SSCm$ there exists $y\in\SSCm$ and $0<a<1$
such that $w=ax+(1-a)y$. In particular, by choosing $x=W(\bOne)$, we
find that $w>0$ componentwise.

For each $\alpha$ in the sequence we can find indices
$(i(\alpha),j(\alpha),i'(\alpha),j'(\alpha))$ as above. Some set of
indices $(i,j,i',j')$ must be repeated infinitely often, since there
are only finitely many choices. Restrict attention to the subsequence
of $\alpha$
for which $i(\alpha)=i$, $j(\alpha)=j$, $i'(\alpha)=i'$, $j'(\alpha)=j'$.

Now, consider the submatrix
\[
\pmatrix{
q_{i,j}(\alpha) & q_{i,j'}(\alpha)\cr
q_{i',j}(\alpha) & q_{i',j'}(\alpha)}.
\]
Recall that $\mathbf{q}(\alpha)\in\INV(\alpha)$.
By Lemma \ref{lem:mwm.sameweight}, and the assumption that
$\blambda>\bZero$ componentwise, every queue is involved in some
maximum weight matching. By Lem\-ma~\ref{lem:mwm.closure}, every
matching is a maximum weight matching.
Let $\bpi$ be any matching with\vadjust{\goodbreak} $\pi_{i,j}=\pi_{i',j'}=1$, and let
$\brho$ be like $\bpi$ except with $(i,j)$ and $(i',j')$ flipped,
that is,
$\rho_{i,j}=\rho_{i',j'}=0$ and
$\rho_{i,j'}=\rho_{i',j}=1$. We can write out explicitly the
difference in weight between these two matchings
\[
\brho\bdot\mathbf{q}(\alpha)^\alpha- \bpi\bdot\mathbf{q}(\alpha
)^\alpha
=
q_{i',j}(\alpha)^\alpha+ q_{i,j'}(\alpha)^\alpha-
q_{i',j'}(\alpha)^\alpha- q_{i,j}(\alpha)^\alpha.\vspace*{-2pt}
\]
Here, $\mathbf{q}(\alpha)^\alpha$ denotes component-wise exponentiation.
Recall that along the subsequence we have chosen,
$q_{i,j}(\alpha)=0$, $q_{i,j'}(\alpha)\geq w_{i\cdot}/M$
and $q_{i',j}(\alpha)\geq w_{\cdot j}/M$.
Therefore $\liminf_{\alpha\tend0}
\brho\bdot\mathbf{q}(\alpha)^\alpha-\bpi\bdot\mathbf{q}(\alpha
)^\alpha>0$.
This contradicts the finding that every matching is a maximum weight
matching for $\mathbf{q}(\alpha)$.

Thus, we have contradicted our original assumption that there exists a~sequence $\tendsdown$ with $w\notin\SSC(\alpha)$.
This completes the proof.\vspace*{-2pt}
\end{pf*}
\begin{pf*}{Proof of Lemma \ref{lem:mwm.closure}}
Let $a=\max_\brho\brho\bdot\bx$
be the weight of the maximum weight matching, let
$\cA=\{\brho\dvtx \brho\bdot\bx=a \}$ and
let $\bsigma= \sum_{\brho\in\cA} \brho$. Observe that
$\bsigma\in\IntegersP^{M \times M}$
and $A_{i,j} = 1_{\{\sigma_{i,j}>0\}}$. Therefore,
$\bsigma\geq\bA$ componentwise. Further, by definition
$\bA\geq\bpi$. Therefore, the matrix $\bsigma-\bpi$ is
nonnegative. Since~$\bsigma$ is sum of $|\cA|$
permutation matrices, all its row sums and column sums
are equal to $|\cA|$. And since~$\bpi$ is a permutation
matrix as well, the matrix~$\bsigma-\bpi$ has all
its row sums and column sums equal to $|\cA|-1$;
therefore by the Birkhoff--von Neumann decomposition
\[
\bsigma= \bpi+ \sum_{\brho\in\sS} \alpha_\brho\brho,\vspace*{-2pt}
\]
where each $\alpha_\brho\geq0$ and $\sum\alpha_\brho=|\cA|-1$.
Now, $\bsigma\bdot\bx=|\cA| a$ by construction of $\bsigma$.
Therefore
\begin{eqnarray*}
|\cA| a
&=&
\bsigma\bdot\bx
=
\bpi\bdot\bx+ \sum_{\brho\in\sS} \alpha_\brho\brho\bdot\bx\\[-2pt]
&\leq&
\bpi\bdot\bx+ (|\cA| - 1) a
\qquad\mbox{because } a=\max_{\brho\in\sS}\brho\bdot\bx.\vspace*{-2pt}
\end{eqnarray*}
Rearranging, $\bpi\bdot\bx\geq a$. But $a$ is the weight of the
maximum weight matching, thus $\bpi\bdot\bx=a$.\vspace*{-2pt}
\end{pf*}
\begin{pf*}{Proof of Lemma \ref{lem:mwm.sameweight}}
Since $\blambda\in\Lambda$, $\blambda\leq\bsigma$ for some
$\bsigma\in\langle{\sS}\rangle$, that is, $\bsigma=\sum_{\bpi\in
\sS
}a_\bpi\bpi$
where $\sum a_\bpi=1$ and each $a_\bpi\geq0$. Therefore
\[
\blambda\bdot\mathbf{q}^\alpha
\leq
\bsigma\bdot\mathbf{q}^\alpha
=
\sum a_\bpi\bpi\bdot\mathbf{q}^\alpha
\leq\Bigl(\sum a_\bpi\Bigr) \max_{\bpi\in\sS}\bpi\bdot\mathbf
{q}^\alpha
=
\max_{\bpi\in\sS}\bpi\bdot\mathbf{q}^\alpha.\vspace*{-2pt}
\]
But by Lemma \ref{lem:fixedpoint}(iv),
$\blambda\bdot\mathbf{q}^\alpha=\max_\bpi\bpi\bdot\mathbf
{q}^\alpha$, therefore both
the inequalities in the above must be equalities.
In particular, all matchings $\bpi$ for which $a_\bpi>0$ are maximum
weight matchings.
If $\lambda_{i,j}>0$, then at least one of these matchings has
$\bpi_{i,j}=1$.\vspace*{-2pt}
\end{pf*}

\subsection{\texorpdfstring{Proof of Theorem \protect\ref
{thm:mwma}(ii).}{Proof of Theorem 8.3(ii)}}
Consider a $2\times2$ switch
with arrival rate matrix $\blambda$. Since $\blambda$ satisfies
(\ref{eq:iq.completeload}), we may write
\[
\blambda=
\pmatrix{
\lambda_{1,1} & 1-\lambda_{1,1}\vspace*{-1pt}\cr
1-\lambda_{1,1} & \lambda_{1,1}}\vspace*{-2pt}\vadjust{\goodbreak}
\]
for some $\lambda_{1,1}\in(0,1)$. To find $\INV(\alpha)$ we use the
characterization from Lem\-ma~\ref{lem:fixedpoint}(iv), which
says that $\mathbf{q}\in\INV(\alpha)$ if and only if
$\blambda\bdot\mathbf{q}^\alpha=\max_\bpi\bpi\bdot\mathbf
{q}^\alpha$, that is, if and
only if
\[
\lambda_{1,1} (q_{1,1}^\alpha+ q_{2,2}^\alpha)
+ (1-\lambda_{1,1}) (q_{1,2}^\alpha+ q_{2,1}^\alpha)
= ( q_{1,1}^\alpha+ q_{2,2}^\alpha)
\mmax
( q_{1,2}^\alpha+ q_{2,1}^\alpha).
\]
Now, the equation $\lambda_{1,1}x+(1-\lambda_{1,1})y=x\mmax y$ is satisfied
if and only if $x=y$, given $0<\lambda_{1,1}<1$. Therefore
\[
\INV(\alpha) = \{ \mathbf{q}\in\RealsP^{2\times2} \dvtx
q_{1,1}^\alpha+q_{2,2}^\alpha= q_{1,2}^\alpha+ q_{2,1}^\alpha
\}.
\]
We wish show that $\{W(\mathbf{q})\dvtx\mathbf{q}\in\INV(\alpha)\}$ is strictly
increasing as $\tendsdown$, where $\WORK(\mathbf{q})
=
(\br_1\bdot\mathbf{q}, \br_2\bdot\mathbf{q},
\bc_1\bdot\mathbf{q}, \bc_2\bdot\mathbf{q})$.
It suffices to show that $\hcW(\alpha) = \{\hWORK(\mathbf
{q})\dvtx\mathbf{q}\in\INV(\alpha)\}$ is
strictly increasing, where
$\hWORK(\mathbf{q})
=
(\br_1\bdot\mathbf{q}, \bc_1\bdot\mathbf{q}, \bOne
\bdot\mathbf{q})
$,
since\vspace*{1pt} there is a~straightforward bijection between $\WORK(\mathbf
{q})$ and $\hWORK(\mathbf{q})$.
Now, $(w_{1\cdot},w_{\cdot1},w_{\cdot\cdot})\in\RealsP^3$ is in
$\hcW(\alpha)$ iff there
exists $\mathbf{q}\in\RealsP^{2\times2}$ such that
\[
q_{1,1}^\alpha+q_{2,2}^\alpha=q_{1,2}^\alpha+q_{2,1}^\alpha,\qquad
\br_1\bdot\mathbf{q}=w_{1\cdot},\qquad
\bc_1\bdot\mathbf{q}=w_{\cdot1},\qquad
\bOne\bdot\mathbf{q}=w_{\cdot\cdot},
\]
that is, iff there exists $x\in\Reals$ such that
%
%
\begin{eqnarray}
\label{eq:mwma2.balance}
&\displaystyle x^\alpha+ (w_{\cdot\cdot}-w_{1\cdot}-w_{\cdot1}+x)^\alpha
- (w_{1\cdot}-x)^\alpha- (w_{\cdot1}-x)^\alpha= 0,&\\
\label{eq:mwma2.bounds}
&\displaystyle \max(0,w_{1\cdot}+w_{\cdot1}-w_{\cdot\cdot}) \leq x
\leq
\min(w_{1\cdot},w_{\cdot1}).&
\end{eqnarray}
Write $\theta(x)$ for the left-hand side of
(\ref{eq:mwma2.balance}). Since $\theta(x)$ is increasing in $x$,
there exists a solution to (\ref{eq:mwma2.balance}) iff
$\theta(x)\leq0$ at the lower bound in
(\ref{eq:mwma2.bounds}) and $\theta(x)\geq0$ at the upper bound.
By considering four separate cases of which of the bounds in
(\ref{eq:mwma2.bounds}) are tight, and after some algebra, we find
that there exists a~solution iff
%
%
\begin{equation}
\label{eq:ssc.w}
w_{i\cdot} + w_{\cdot j} +
(w_{i\cdot}^\alpha+ w_{\cdot j}^\alpha){}^{1/\alpha}
\geq
w_{\cdot\cdot}
\qquad\mbox{for each $i,j\in\{1,2\}$},
\end{equation}
where $w_{2\cdot}=w_{\cdot\cdot}-w_{1\cdot}$ and
$w_{\cdot2}=w_{\cdot\cdot}-w_{\cdot1}$.
Now, it is a standard inequality that for
any $x>0$ and $y>0$, and any $0<\alpha<\beta$,
\[
(x^\alpha+y^\alpha)^{1/\alpha} > (x^\beta+y^\beta)^{1/\beta}.
\]
Applying this inequality to (\ref{eq:ssc.w}), it follows that
$\hcW(\alpha)$ is strictly increasing as \mbox{$\tendsdown$}.
\end{appendix}

\section*{\texorpdfstring{Acknowledgments.}{Acknowledgments}}

We would like to thank Mike Harrison, Frank Kelly, Balaji Prabhakar
and Ruth Williams for helpful discussions. We thank Ciamac Moallemi
for suggesting Lemma \ref{lem:equiv}. We also thank the anonymous
reviewers for their help in improving readability of the manuscript.


%

%
\printaddresses

\end{document}